\documentclass[10pt]{article}
\usepackage{amsmath}
\usepackage{amsfonts}
\usepackage{graphicx}
\usepackage{amssymb}
\usepackage{leftidx}
\usepackage{color}
\allowdisplaybreaks

\newtheorem{condition**}{A*}
\newtheorem{condition***}{C*}
\newtheorem{condition*}{C}

\newtheorem{example}{Example}[section]
\newtheorem{proposition}{Proposition}[section]
\newtheorem{corollary}{Corollary}[section]
\newtheorem{definition}{Definition}[section]
\newtheorem{theorem}{Theorem}[section]
\newtheorem{lemma}{Lemma}[section]
\newtheorem{remark}{Remark}[section]

\begin{document}

\title{BSDE Approach for  $\alpha $-Potential Stochastic Differential Games}
\author{Xin Guo$^{1}$\thanks{Department of Industrial Engineering and Operations Research, University of California, Berkeley, CA, USA. Email: xinguo@berkeley.edu. X. Guo acknowledges the general support from Coleman Fung endowment fund.}, \ 
Xun Li$^{2}$\thanks{Department of Applied Mathematics, The Hong Kong Polytechnic University, Hong Kong, China. E-mail: li.xun@polyu.edu.hk.}, \ Liangquan Zhang$^{3}$\thanks{School of Mathematics, Renmin University of China, Beijing 100872, China. Email: xiaoquan51011@163.com.}} 

\maketitle

\begin{abstract}
In this paper, we examine a class of $\alpha$-potential stochastic
differential games with random coefficients via the  backward stochastic
differential equations (BSDEs) approach. Specifically, we show that the first and second order linear derivatives of the objective function for each player can be expressed through the corresponding first and second-order adjoint equations, which leads to rigorous estimates for $\alpha$. We illustrate the dependence of $\alpha$ on game characteristics through detailed  analysis of  linear-quadratic games, and with common noise. 
\end{abstract}


\tableofcontents

\section{Introduction\label{sect1}}


The birth of game theory as a mathematical discipline is attributed to Von Neumann and Morgenstern from their seminal work \cite{NM44}, with some foundational ideas and concepts rooted in military strategy  and ancient philosophical discussions.
One of the most influential notions in game theory is   Nash Equilibrium (NE), introduced  by Nash in 1950s \cite{Nash50} for competitive games. A Nash equilibrium is a set of strategies where no player has incentive to deviate from the chosen strategy, given strategies of other players.
The inception of NE has significantly expanded the scope and influence of game theory beyond economics, with a wide arrange of applications including political science, biology, engineering, and social psychology \cite{ms68}. 
 
Meanwhile, recent rapid development in machine learning theory and in particular in multi-agent reinforcement learning leads to renewed interests in  studies of  competitive games in dynamic systems,   including  autonomous driving, smart grid management, and crowd funding. In these systems, players interact with random environment and among themselves via strategic plays.
 
However, NEs for competitive games are notoriously difficult to derive or approximate  due to interdependence among players with computational complexity increasing exponentially with respect to the number of players, see for instance  \cite{Papadimitriou2007}.

\paragraph{Potential games.}
Potential games, introduced by Monderer and Shapley in \cite{ms96},  are a class of static games  where the incentive of all players to change their strategy can be expressed through a single function called  ``potential function''. Potential games  simplifies the otherwise complicated task of identifying a Nash equilibrium for $N$-player competitive  games to  optimizing  this single function.  This alignment between changes in individual player payoffs and the changes of the potential function  implies that individual rational behavior leads towards a collectively  stable state that can be analyzed from a global perspective. 
It enables modeling and analyzing diverse real-world problems across various fields, including market competitions  and resource management in an economic system \cite{nlks22}, traffic routing and energy distribution in a network systems \cite{lv2024,mas2009,zmz2021}, and multi-agent games \cite{dwzj,lopp21,mzz18,mwps22} and controls of multi-robot \cite{swhw2024} in engineering systems. 

Despite the ingenious idea of potential function and the promising power of the potential game paradigm, most dynamic games are not Markov
potential games. It is well known  that even a Markov game where the game at each state is a static potential game
may not be a dynamic potential game \cite{lopp21}. In practice, dynamic potential game framework imposes restrictive assumptions that are infeasible and inappropriate for modeling dynamic games with complex game interactions and dynamics \cite{gz23}.

\paragraph{$\alpha$-potential games.} 
$\alpha$-potential games are a class of  dynamic games,  proposed first   for $N$-player non-cooperative discrete-time Markov games with finite states and actions  \cite{glmsw03}, and later  for  general  non-cooperative games  \cite{GLZ1}.
This dynamic framework generalizes the traditional potential game:  the change in a player's objective function upon her unilateral deviation is equal to the change of a single function called ``$\alpha$-potential function,''  up to an error $\alpha$. When $\alpha=0$, the game is simply a potential game and the auxiliary function is now a potential function. 

Unlike the the restrictive framework of  potential games, this simple modification  by the parameter $\alpha$ provides a  tractable and powerful way to analyze and approximate equilibria in  complex dynamics and multi-agent systems. It is shown in  \cite{glmsw03}  that a general finite-state and finite-action discrete-time dynamic game can be formulated as an $\alpha$-potential game through a semi-infinite program which identifies the precise $\alpha$. Moreover, a number of popular dynamic games which are not potential games can now be analyzed in this framework with convergent computational algorithms and NE-regret bound analysis.

In an $\alpha$-potential games, maximizing a single $\alpha$-potential function yields an $\alpha$-Nash equilibrium (NE), i.e., an NE up to an error $\alpha$.
Meanwhile, the  optimization problem for the associated $\alpha$-potential function is equivalent to  a conditional MKV control problem, as shown and analyzed in details in  \cite{GLZ1}.
Moreover, the parameter $\alpha$ is shown (\cite{glmsw03} and \cite{GLZ1})  to be capable of capturing the degree of heterogeneity and asymmetry among players' objective functions and interactions for sequential decision makings. 
These results are particularly relevant for practical applications, such as network routing problems or resource allocation in dynamic systems, where exact solutions may be computationally infeasible and where approximate equilibria and in particular the quantity of $\alpha$ nonetheless provide valuable insights for strategic interactions among players and for game design.

\paragraph{$\alpha$-NE and estimation of $\alpha$.}
  $\alpha$-potential games exploit through  $\alpha$ to obtain $\alpha$-NE as well as crucial information and insight of the game strategy and its dependence on game characteristics.  Given aforementioned mathematical developments, the most crucial remaining piece is
   now a precise bound analysis of $\alpha$ for general nonzero-sum games.  
For discrete-time competitive  games, the  precise estimate of its $\alpha$ can be obtained by a semi-infinite linear program \cite{glmsw03}.
The analysis of  $\alpha$ for continuous-time differential games remains a challenging problem. 

 For a special class of stochastic differential games, \cite{GLZ1} shows intriguing  dependence of $\alpha$ on a number of factors, including 
the intensity of interactions among the players through either the state dynamics and/or through their value functions, and their strategy choices. In the particular case of distributed games where the strategy of players does not influence the 
dynamics of the underlying state with special structures of objective functions, $\alpha=0$; in a network game where players interact weakly, $\alpha$ is shown to decay to $0$ with varying degree depending on the intensity  of players' interaction. In games with mean-field type interactions, the order of $\alpha$ would differ for open-loop and closed-loop strategies. 
The key technical tool in \cite{GLZ1} is the first and the second order sensitivity
processes of the state dynamics with respect to controls. However, their approach for $\alpha$ 
seems to encounter a major technical barrier for
 more general form of stochastic differential games.
For instance, for the case with   diffusion controls, they are unable  to explicitly bound the second order  derivative of the objective functions.

Since the precise estimate of $\alpha$ is the key to evaluate the quality of the optimization solution for the associated $\alpha$-potential functions, i.e., the 
solution of the $\alpha$-NE, alternative approaches  are needed in order to obtain more precise $\alpha$ and for general classes of differential 
games.

\paragraph{Our approach.} In
this paper, we adopt a BSDE approach to analyze the bounds of $\alpha$. The  key idea is to revisit the first and the second order  sensitivity processes exploited in \cite{GLZ1} through their BSDE roots.
Observing that these processes  are crucial for  deriving stochastic maximum principle and for backward stochastic differential equations, we   reformulate them in BSDE forms.  Along with the associated duality principle from BSDE, this reformulation  enables us to eliminate the second-order sensitivity process through the duality principle, and to avoid  expensive and complex computations for the second order variational equation. Consequently, the problem of estimating $\alpha$ is transformed into analyzing a class of linear BSDE, even when the   diffusion term is controlled and includes both state and control variables with possibly random  coefficients. More importantly, it enables the analysis  of stochastic differential games, and with common noise. 
The effectiveness of this BSDE approach is demonstrated via several examples, including linear quadratic games on directed graphs,  large-population systems with $N$ weakly-coupled agents, as well as  a linear-quadratic mean-field game with common noise.

In contrast to the technique presented in \cite{GLZ1},  the methodology of (BSDEs)  enables more precise estimate of $\alpha$ within a general game framework. More importantly, our analysis also highlights  the potential of  this BSDE approach, which  may be extended to handle other forms of admissible policies beyond the open-loop case studied here.
Indeed, 
when the dynamic programming principle becomes intractable, for instance due to the lack of continuity or differentiability in  the value function or for a high dimensional state space where  the infinite dimensional partial differential equations in the Wasserstein space is  too challenging to solve,  BSDEs may offer a more feasible  approach to  handle non-smooth value functions and allow for connecting optimal control with the state process via the martingale representation theorem,
see \cite{GLZ15clo}.

\paragraph{Organization.} The remainder of this paper is organized as follows. Section \ref{sect2} presents some preliminary results concerning $\alpha$-potential games and the BSDE theory.  Section \ref{sect3} derives the BSDEs representation of $\frac{\protect\delta V_{i}}{\protect\delta %
u_{h}}$ and $\frac{\protect\delta ^{2}V_{i}}{\protect\delta u_{h}\protect%
\delta u_{\ell }}$. Subsequently, Section \ref{sect4} contains the estimate of $\alpha$ via the BSDE technique, with illustrations of several special classes of games.
   Proofs of all main results are in Section \ref{sect6}. After concluding remarks in Section \ref{sect7}, the proofs of technical lemmas are included in Appendix \ref{APP}. 


\section{Mathematical Set-up and Preliminaries}
\label{sect2}

\paragraph{Notation.}
Let us first introduce some notation.  Throughout this paper, we
denote the $k$ dimensional Euclidean space by $\mathbb{R}^{k}$ with the standard
Euclidean norm $|\cdot |$ and the standard Euclidean inner product $\langle
\cdot ,\cdot \rangle $, the transpose of a vector (or matrix) $x$ 
by $x^{\top }$,  the trace of a square matrix $A$ by $\text{Tr}(A)$. Let
$\mathbb{R}^{m\times n}$ be the Hilbert space consisting of all ($m\times n$%
)-matrices with the inner product $\langle A,B\rangle :=\text{Tr}(AB^{\top
}) $ and the norm $|A|:=\langle A,A\rangle ^{\frac{1}{2}}=\sqrt{\text{Tr}%
(AA^{\top })}$. $\mathbb{S}^{n}$ denotes the set of symmetric $n\times n$ matrices.
For any $A\in \mathbb{R}^{k\times m},$ we adopt the usual matrix norm in $\mathbb{R}^{k\times m}$ as follows:

\begin{eqnarray*}
\begin{aligned}
\left\Vert A\right\Vert _{2} & \triangleq \sqrt{\max \sigma \left( AA^{\top
}\right) }\leq \sqrt{\text{Tr} \left( AA^{\top }\right) }\equiv
\left\vert A\right\vert  \\
& \leq \sqrt{k\wedge m}\sqrt{\max \sigma \left( AA^{\top }\right) }\equiv
\sqrt{k\wedge m}\left\Vert A\right\Vert _{2},
\end{aligned}
\end{eqnarray*}%
where $\sigma \left( AA^{\top }\right) $ is the set of all eigenvalues of $AA^{\top}$. 

We will denote $(\Omega ,\mathcal{F},\mathbb{F}=\left( \mathcal{F}%
_{t}\right) _{0\leq t\leq T},\mathbb{P})$ as a complete filtered probability
space satisfying the usual condition, where $N$ independent standard
Brownian motion $W^{i}\left( \cdot \right) ,$ $i\in I_{N},$ are defined and
each takes value in $\mathbb{R}$ on $\left[ 0,T\right] $ for a fixed finite
time $T$. For notational simplicity, when $x$ belongs to an Euclidean space,
we will write $x_{i}$ for its $i$-th coordinate and $|x|$ for its Euclidean
norm. When there is no ambiguity, we denote $\mathbb{F}^{i}\mathbb{=}\left\{
\mathcal{F}_{t}^{i}\right\} _{0\leq t\leq T}$ for the natural filtration of $%
W^{i}$, augmented by all $\mathbb{P}$-null sets in $\mathcal{F}$. For each $%
p\geq 1$, let $\mathcal{S}^{p}(\mathbb{R}^{N})$ be the space of $\mathbb{R}%
^{N}$-valued $\mathbb{F}$-progressively measurable processes $X:\Omega
\times \left[ 0,T\right] \rightarrow \mathbb{R}^{N}$ satisfying $\left\Vert
X\right\Vert _{\mathcal{S}^{p}(\mathbb{R}^{N})}=\mathbb{E}\left[ \sup_{s\in %
\left[ 0,T\right] }\left\vert X_{s}\right\vert ^{p}\right] ^{1/p}<\infty $,
and let $\mathcal{H}^{p}(\mathbb{R}^{N})$ be the space of $\mathbb{R}^{N}$%
-valued $\mathcal{F}$-progressively measurable processes $X:\Omega \times %
\left[ 0,T\right] \rightarrow \mathbb{R}^{N}$ satisfying $\mathbb{E}\left[
\int_{0}^{T}\left\vert X_{s}\right\vert ^{p}\mathrm{d}s\right] ^{1/p}<\infty
$. Let $L^{2}\left( \Omega ;\mathbb{R}\right) $ denote the space of square
integrable $\mathcal{F}_{0}$-measurable random variables. $\delta _{i,j}$
denotes the Kronecker delta such that $\delta _{i,j}=1$ if $i=j$ and $0$
otherwise.  We
denote by $\mathcal{M}(0,T;\mathbb{R}^{n})$ the set of all $\mathbb{R}^{n}$%
-valued processes $\{\varphi _{t}\}_{0\leq t\leq T}$ that are $\mathcal{F}%
_{t}$-adapted and satisfy $\mathbb{E}\left[ \int_{0}^{T}|\varphi _{t}|^{2}%
\mathrm{d}t\right] <\infty $. Whenever there is no risk case of ambiguity, the statement
\textquotedblleft for a.a. $t\in \lbrack 0,T]$, a.s. $\omega \in \mathbb{P}$
($\mathbb{P}$-a.s.)\textquotedblright\ will be abbreviated to
\textquotedblleft for a.a. $\left( t,\omega \right) $\textquotedblright .

\paragraph{Mathematical setup and assumptions.}
The focus of this paper is the analysis of stochastic non-zero-sum games and the estimate of its $\alpha$. For sake of simplicity, we focus on the one-dimensional case.
Specifically, 
consider the differential game $\mathcal{G}^{\mathsf{op}}=\left(
I_{N},\left( \mathcal{A}_{i}\right) _{i\in I_{N}},\left( V_{i}\right) _{i\in
I_{N}}\right) $ defined as follows: let $I_{N}=\left\{ 1,\ldots ,N\right\} $%
, and for each $i\in I_{N}$, let $\mathcal{A}_{i}$ be a convex subset of $%
\mathcal{H}^{2}(\mathbb{R})$ representing the set of admissible open-loop
controls of player $i$. For each $\mathbf{u}=\left( u_{i}\right) _{i\in
I_{N}}$, $\mathbf{X}^{\mathbf{u}}=\left( X_{i}^{\mathbf{u}}\right)
_{i=1}^{N}\in \mathcal{S}^{2}(\mathbb{R}^{N})$ is the associated state
process governed by the following dynamics: For all $i\in I_{N}$ and $\xi
_{i}\in L^{2}\left( \Omega ;\mathbb{R}\right) $,
\begin{equation}
\left\{\begin{aligned}
\mathrm{d}X_{t,i} & = b_{i}\left( t,X_{t,i},\mathbf{X}_{t}^{\mathbf{u}%
},u_{t,i}\right) \mathrm{d}t+\sigma _{i}\left( t,X_{t,i},\mathbf{X}_{t}^{%
\mathbf{u}},u_{t,i}\right) \mathrm{d}W_{t}^{i}, \\
X_{0,i} & = \xi _{i},\text{ }\forall t\in \lbrack 0;T],%
\end{aligned}\right.  \label{sde1}
\end{equation}
where $b_{i}: [0,T] \times \mathbb{R\times R}^{N}\times \mathbb{R}\rightarrow \mathbb{R}$ and $\sigma _{i}:
[0,T] \times \mathbb{R\times R}^{N}\times \mathbb{R}\rightarrow
\mathbb{R}$
satisfy appropriate conditions to ensure the definedness of (\ref{sde1}), as below.
The cost functional $V_{i}:\mathcal{A}^{\left( N\right) }\subset \mathcal{H}^{2}\left( \mathbb{R}\right) ^{N}\rightarrow \mathbb{R}$ of player $i$ is
given by
\begin{equation}
V_{i}\left( \mathbf{u}\right) =\mathbb{E}\left[ \int_{0}^{T}f_{i}\left( t,%
\mathbf{X}_{t}^{\mathbf{u}},\mathbf{u}_{t}\right) \mathrm{d}t+g_{i}\left(
\mathbf{X}_{T}^{\mathbf{u}}\right) \right] ,  \label{cost}
\end{equation}%
where $f_{i}: [ 0,T] \times \mathbb{R}^{N}\times
\mathbb{R}^{N}\rightarrow \mathbb{R}$ and $g_{i}: \mathbb{R}^{N}\rightarrow \mathbb{R}$ satisfy appropriate conditions to  ensure the definedness of the game, as specified below.

\begin{enumerate}
\item[\textbf{(A1)}] \textbf{[Assumptions for the dynamics] }For each $i\in
I_{N},$ $\varphi _{i}=b_{i},$ or $\sigma _{i}$

(i) The maps $\varphi _{i}$ are $\mathcal{B}\left( \left[ 0,T\right] \times \mathbb{R\times
R}^{N}\times \mathbb{R}\right) \otimes \mathcal{F}_{T}^{i}$-measurable.

(ii) For all $\left( x,y,u\right) \in \mathbb{R\times R}^{N}\times \mathbb{R}
$, the process $\varphi _{i}\left(t,x,y,u\right) \in \mathbb{R}$ is $\mathbb{F}$-adapted for $t \in [0, T]$.

(iii) There exists some positive constants $L^{\varphi _{i}}$ and  
$L_{y}^{\varphi _{i}}$ such that for almost all $\left( t,\omega \right)
\in \left[ 0,T\right] \times \Omega $, $\left( t,x,y,u\right) \in \left[ 0,T%
\right] \times \mathbb{R\times R}^{N}\times \mathbb{R}$ and $i,j\in I_{N},$%
\begin{equation*}
\left\{\begin{array}{l}
\left\vert \varphi _{i}\left( t,0,0,u\right) \right\vert \leq L^{\varphi
_{i}}\left( 1+\left\vert u\right\vert \right), \\
\left\vert \left( \partial _{x}\varphi _{i}\right) \left( t,x,y,u\right)
\right\vert +\left\vert \left( \partial _{u}\varphi _{i}\right) \left(
t,x,y,u\right) \right\vert +\left\vert \left( \partial _{uu}\varphi
_{i}\right) \left( t,x,y,u\right) \right\vert \leq L^{\varphi _{i}}, \\
\left\vert \left( \partial _{xx}\varphi _{i}\right) \left( t,x,y,u\right)
\right\vert +\left\vert \left( \partial _{xu}\varphi _{i}\right) \left(
t,x,y,u\right) \right\vert \leq L^{\varphi _{i}}, \\
\left\vert \left( \partial _{y_{i}}\varphi _{i}\right) \left( t,x,y,u\right)
\right\vert \leq \frac{L_{y}^{\varphi _{i}}}{N}, \\
\left\vert \left( \partial _{xy_{i}}^{2}\varphi _{i}\right) \left(
t,x,y,u\right) \right\vert +\left\vert \left( \partial _{uy_{i}}^{2}\varphi
_{i}\right) \left( t,x,y,u\right) \right\vert \leq \frac{L_{y}^{\varphi _{i}}%
}{N}, \\
\left\vert \left( \partial _{y_{i}y_{j}}^{2}\varphi _{i}\right) \left(
t,x,y,u\right) \right\vert \leq \frac{L_{y}^{\varphi _{i}}}{N}1_{i=j}+\frac{%
L_{y}^{\varphi _{i}}}{N^{2}}1_{i\neq j}.%
\end{array}%
\right.
\end{equation*}
For any $\psi =\left( \psi _{i}\right) _{i\in I_{N}}$, where $\psi _{i}$
satisfies (A1), we define $L^{\psi }=\displaystyle\max_{i\in I_{N}}L^{\varphi _{i}}$
and $L_{y}^{\psi }=\displaystyle\max_{i\in I_{N}}L_{y}^{\varphi _{i}}$. 
\end{enumerate}
\begin{enumerate}
\item[\textbf{(A2)}] \textbf{[Assumptions for the costs]} (i) The maps $%
f_{i} $ and $g_{i}$ are respectively $\mathcal{B}\left( \left[ 0,T\right]
\times \mathbb{R}^{N}\times \mathbb{R}\right) \otimes \mathcal{F}_{T}^{i}$
and $\mathcal{B}\left( \mathbb{R}^{N}\right) \otimes \mathcal{F}_{T}^{i}$
measurable.

(ii) For all $\left( y,u\right) \in \mathbb{R}^{N}\times \mathbb{R}$, the
process $\left[ 0,T\right] \ni t\rightarrow f_{i}\left( t,y,u\right) \in
\mathbb{R}$ is $\mathbb{F}$-adapted.

(iii) For almost all $\left( t,\omega \right) \in \left[ 0,T\right] \times
\Omega,$ $f_{i}$ and $g_{i}$ are twice continuously
differentiable with respect to $\left( y,u\right) \in \mathbb{R}^{N}\times
\mathbb{R}^{N}$ such that $\sup_{t\in \left[ 0,T\right] }\left( \left\vert
f_{i}\left( t,0,0\right) \right\vert +\left\vert \left( \partial _{\left(
y,u\right) }f_{i}\right) \left( t,0,0\right) \right\vert \right) <\infty $
and second-order derivatives are uniformly bounded.
\end{enumerate}
\begin{remark}
In the special case of
\begin{eqnarray*}
\begin{aligned}
b_{i}\left( t,x,y,u\right) & = A_{i}x+\frac{1}{N}\sum_{k=1}^{N}y_{k}+B_{i}u+\mathfrak{b}_{t,i}, \\
\sigma _{i}\left( t,x,y,u\right) & = C_{i}x+\frac{1}{N}\sum_{k=1}^{N}y_{k}+D_{i}u+\mathfrak{l}_{t,i},
\end{aligned}
\end{eqnarray*}
where $\mathfrak{b}_{t,i}$ and $\mathfrak{l}_{t,i}$ are bounded processes. Then 
\begin{eqnarray*}
\begin{aligned}
\left\vert b_{i}\left( t,0,0,u\right) \right\vert +\left\vert \sigma
_{i}\left( t,0,0,u\right) \right\vert & \leq \left\vert B_{i}u+\mathfrak{b}%
_{t,i}\right\vert +\left\vert D_{i}u+\mathfrak{l}_{t,i}\right\vert \\
& \leq \left\vert B_{i}u\right\vert +\left\vert \mathfrak{b}%
_{t,i}\right\vert +\left\vert D_{i}u\right\vert +\left\vert \mathfrak{l}%
_{t,i}\right\vert \\
& \leq C\left( 1+\left\vert u\right\vert \right)
\end{aligned}
\end{eqnarray*}
for some constant $C>0$. 
\end{remark}

\begin{remark}
\label{re1}Under (A1), it is easy to check that by the mean value theorem,
\begin{eqnarray*}
\begin{aligned}
\left\vert \varphi _{i}\left( t,x,y,u\right) \right\vert 
& = \left\vert \varphi _{i}\left( t,x,y,u\right) -\varphi _{i}\left(
t,0,0,u\right) +\varphi _{i}\left( t,0,0,u\right) \right\vert \\
& \leq \left\vert \varphi _{i}\left( t,x,y,u\right) -\varphi _{i}\left(
t,0,0,u\right) \right\vert +\left\vert \varphi _{i}\left( t,0,0,u\right)
\right\vert \\
& \leq L^{\varphi _{i}}\left\vert x\right\vert +\frac{L_{y}^{\varphi _{i}}}{N%
}\left\vert y\right\vert +L^{\varphi _{i}}\left( 1+\left\vert u\right\vert
\right) \\
& \leq L^{\varphi _{i}}\left( 1+\left\vert x\right\vert +\left\vert
u\right\vert \right) +\frac{L_{y}^{\varphi _{i}}}{N}\sum_{i=1}^{N}\left\vert
y_{i}\right\vert.
\end{aligned}
\end{eqnarray*}
\end{remark}

\begin{remark}
In (A1), for every $i\in I_{N}$, the partial derivatives of $y\rightarrow
b_{i}\left( t,x,y,u\right) $ admit explicit decay rates in terms of $N$.
This assumption naturally holds if each player's state depends on the
empirical measure of the joint state process, i.e., the mean-field
interaction. For instance 
\begin{equation*}
b_{i}\left( t,x,y,u\right) =\psi _{i}\left( t,x,\frac{1}{N}%
\sum_{j=1}^{N}\delta _{y_{j}},u\right)
\end{equation*}
with $\left( t,x,y,u\right) \in \left[ 0,T\right] \times \mathbb{R}\times
\mathbb{R}^{N}$ for a sufficient smooth measurable function $\psi _{i}:\left[ 0,T\right]
\times \mathbb{R}\times \mathcal{P}_{2}\left( \mathbb{R}\right) \times
\mathbb{R}\rightarrow \mathbb{R}$, where $\mathcal{P}_{2}\left( \mathbb{R}%
\right) $ is the space of probability measures on $\mathbb{R}$ with second
moments, (see Propositions 5.35 and 5.91 in \cite{cd18}). Let $%
i,j\in I_{N}$ with $i<j,$ and define $\Delta _{i,j}^{f}=f_{i}-f_{j}$ and $%
\Delta _{i,j}^{g}=g_{i}-g_{j}.$  

In general, the dependence of the
constant $L_{y}^{b}$ reflects the degree of
coupling among all players' state dynamics. For instance, if $L_{y}^{b}$
remains bounded as $N\rightarrow \infty $, then the state dynamics can have
mean-field type interactions. Alternatively, if $L_{y}^{b}=0,$ then all
players' states are decoupled.
\end{remark}

\paragraph{$\alpha$-potential games and $\alpha$-Nash equilibrium.}
Next, recall the definition of
$\alpha $-potential games and some of its properties, introduced and derived in \cite{glmsw03} for the discrete-time setting and in \cite{GLZ1} for general games.

It starts with  a finite player  game $\mathcal{G}=(I_{N},\mathcal{S},\left( \mathcal{A}%
_{i}\right) _{i\in I_{N}},\left( V_{i}\right) _{i\in I_{N}}),$ with  $I_{N}=\left\{ 1,\ldots, N\right\}$, $N\in \mathbb{N}$ the number of player, $\mathcal{S}$  the state space of the underlying dynamics, $\mathcal{A}_{i}$ a subset of a real vector space
representing all admissible strategies of player $i$, and $\mathcal{A}%
^{\left( N\right) }=\prod_{i\in I_{N}}$ the set of strategy profiles for
all players. For each $i\in I_{N}$, $V_{i}:\mathcal{A}^{\left( N\right)
}\rightarrow \mathbb{R}$ is the cost function of player $i$, where $%
V_{i}\left( \mathbf{u}\right) $ is player $i$'s expected cost if the state
dynamics starts with a fixed initial state $s_{0}\in \mathcal{S}$ and all
players take the strategy profile $\mathbf{u\in }\mathcal{A}^{\left(
N\right) }$.  We
denote by $\mathcal{A}_{-i}^{\left( N\right) }=\prod_{j\in I_{N}\backslash
\left\{ i\right\} }$ the set of strategy profiles of all players except
player $i$, and by $\mathbf{u}$ and $u_{-i}$ a generic element of $\mathcal{A}^{\left( N\right) }$ and $\mathcal{A}_{-i}^{\left( N\right) }$,
respectively. For any $i\in I_{N}$, player $i$ aims to minimize $V_{i}$ over all admissible strategies in $\mathcal{A}_{i}$.

\begin{definition}\label{pf}
Given a game $\mathcal{G}=(I_{N},\mathcal{S},\left( \mathcal{A}_{i}\right)
_{i\in I_{N}},\left( V_{i}\right) _{i\in I_{N}})$, if there exists $\alpha
>0 $ and $\Phi :\mathcal{A}^{\left( N\right) }\rightarrow \mathbb{R}$ such
that for all $i\in I_{N}$, $a_{i},a_{i}^{\prime }\in \mathcal{A}_{i}$ and $a_{-i}\in \mathcal{A}_{-i}^{\left( N\right)}$,
\begin{equation*}
\left\vert V_{i}\left( \left( a_{i}^{\prime },a_{-i}\right) \right)
-V_{i}\left( \left( a_{i},a_{-i}\right) \right) -\left( \Phi \left( \left(
a_{i}^{\prime },a_{-i}\right) \right) -\Phi \left( \left(
a_{i},a_{-i}\right) \right) \right) \right\vert \leq \alpha ,
\end{equation*}
then  $\mathcal{G}$ is called an $\alpha $-potential game, and $\Phi $ is an $%
\alpha $-potential function for $\mathcal{G}$. In the case of $\alpha =0$,  the game $\mathcal{G}$ is simply a potential game with $\Phi $ a
potential function for $\mathcal{G}$.
\end{definition}

The following proposition \cite{glmsw03} underscores the importance of explicitly characterizing $\alpha $ for a given game. 

\begin{proposition}
\label{an}Let $\mathcal{G}$ denote an $\alpha $-potential game for some $%
\alpha $ with $\Phi$ the associated $\alpha $-potential function. For
each $\epsilon >0$, if there exists $\mathbf{a}^{\ast }\in \mathcal{A}%
^{\left( N\right) }$ such that $\Phi \left( \mathbf{a}^{\ast }\right) \leq
\inf_{\mathbf{a}\in \mathcal{A}^{\left( N\right) }}\Phi \left( \mathbf{a}%
\right) +\epsilon $, then $\mathbf{a}^{\ast }$ constitutes an $\left( \alpha
+\epsilon \right) $-Nash equilibrium of $\mathcal{G}$.
\end{proposition}

To analyze the change of strategies,  a scalar-valued
function with respect to unilateral deviations of strategies,  called linear derivative, is introduced \cite{glmsw03}.  For each $i\in
I_{N}$, denote by span($\mathcal{A}_{i}$) the vector space of all linear
combinations of strategies in $\mathcal{A}_{i}$, i.e.,
\begin{equation*}
\text{span}(\mathcal{A}_{i})=\left\{ \sum_{l=1}^{m}c_{l}a^{\left( l\right)
}\left\vert m\in N,\left( c_{l}\right) _{l=1}^{m}\subset \mathbb{R},\left(
a^{\left( l\right) }\right) _{l=1}^{m}\subset \mathcal{A}_{i}\right.
\right\} .
\end{equation*}

\begin{definition}
\label{ld}Let $\mathcal{A}^{\left( N\right) }=\prod_{i\in I_{N}}$ be a
convex set and $f:\mathcal{A}^{\left( N\right) }\rightarrow \mathbb{R}$. For
each $i\in I_{N}$, we say $f$ has a linear derivative with respect to $%
\mathcal{A}_{i}$, if there exists $\frac{\delta f}{\delta a_{i}}:\mathcal{A}%
^{\left( N\right) }\times $span$(\mathcal{A}_{i})\rightarrow \mathbb{R}$,
such that for all $\mathbf{a=}\left( a_{i},a_{-i}\right) \in \mathcal{A}%
^{\left( N\right) }$, $\frac{\delta f}{\delta a_{i}}\left( \mathbf{a};\cdot
\right) $ is linear and
\begin{equation}
\lim_{\varepsilon \searrow 0}\frac{f\left(a_{i}+\varepsilon \left(
a_{i}^{\prime }-a_{i}\right), a_{-i} \right) -f\left( \mathbf{a}\right) }{\varepsilon }=\frac{\delta f}{\delta a_{i}}\left( \mathbf{a};a_{i}^{\prime}-a_{i}\right),\text{ }\forall a_{i}^{\prime }\in \mathcal{A}_{i}.
\label{lfd}
\end{equation}
Moreover, for each $i,j\in I_{N}$, we say $f$ has second-order linear
derivatives with respect to $\mathcal{A}_{i}\times \mathcal{A}_{j}$, if (1)
for all $k\in \left\{ i,j\right\}$, $f$ has a linear derivative $\frac{%
\delta f}{\delta a_{k}}$ with respect to $\mathcal{A}_{k}$, and (2) for all $%
\left( k,l\right) \in \left\{ \left( i,j\right) ,\left( j,i\right) \right\}$, there exists $\frac{\delta ^{2}f}{\delta a_{k}\delta a_{l}}:\mathcal{A}^{\left( N\right) }\times $span$(\mathcal{A}_{i})\times $span$(\mathcal{A}%
_{j})\rightarrow \mathbb{R}$ such that for all $\mathbf{a\in }\mathcal{A}%
^{\left( N\right) },$ $\frac{\delta ^{2}f}{\delta a_{k}\delta a_{l}}\left(
\mathbf{a};\cdot ,\cdot \right) $ is bilinear and for all $a_{k}^{\prime
}\in $span$(\mathcal{A}_{k}),$ $\frac{\delta ^{2}f}{\delta a_{k}\delta a_{l}}%
\left( \cdot ;a_{k}^{\prime },\cdot \right) $ is a linear derivative of $%
\frac{\delta f}{\delta a_{k}}\left( \cdot ;a_{k}^{\prime }\right) $ with
respect to $\mathcal{A}_{l}$. We refer to $\frac{\delta ^{2}f}{\delta
a_{i}\delta a_{j}}$ as the
second-order linear derivative of f with respect to $\mathcal{A}_{i}\times
\mathcal{A}_{j}$.
\end{definition}

For simplicity, we set $\frac{\delta f}{\delta a_{i}} = \partial_{a_i}f$ and $\frac{\delta ^{2}f}{\delta a_{i}\delta a_{j}} = \partial_{a_ia_j}^2f$. 

With this notion of linear derivatives, one can construct (see Theorem 2.4 in \cite{GLZ1}) an $\alpha$-potential function and provide an upper bound for the associated parameter $\alpha$ for a given game $\mathcal{G}$, where all players' strategy sets are convex. 
\begin{lemma} [\cite{GLZ1}]
    For each $i \in I_N$, let span($\mathcal{A}_i$) denote the vector space consisting of all linear combinations of strategies in $\mathcal{A}_i$. If the expected cost functions of the game $\mathcal{G}$ admit second-order partial derivatives, then under certain mild regularity conditions, for any fixed $%
\mathbf{z}\in \mathcal{A}^{(N)}$, the function%
\begin{equation}
\Phi \left( \mathbf{a}\right) =\int_{0}^{1}\sum_{j=1}^{N}\frac{\partial V_{j}%
}{\partial a_{j}}\left( z+r\left( \mathbf{a-z}\right) ;a_{j}-z_{j}\right)
\mathrm{d}r  \label{apf}
\end{equation}%
is an $\alpha $-potential function of $\mathcal{G},$ with
\begin{equation}
\alpha \leq 2\sup_{i\in I_{N},a_{i}^{\prime }\in \mathcal{A}_{i},\mathbf{a,a}%
^{\prime \prime }\in \mathcal{A}^{(N)}}\sum_{j=1}^{N}\left\vert \frac{\delta
^{2}V_{i}}{\delta a_{i}\delta a_{j}}\left( \mathbf{a};a_{i}^{\prime
},a_{j}^{\prime \prime }\right) -\frac{\delta ^{2}V_{j}}{\delta a_{j}\delta
a_{i}}\left( \mathbf{a};a_{j}^{\prime \prime },a_{i}^{\prime }\right)
\right\vert .  \label{aph}
\end{equation}
\end{lemma}

It is therefore critical to estimate the linear derivatives of $V_i$, in order to get a sharp bound of $\alpha$.

In \cite{GLZ1}, for stochastic differential games with open-loop admissible strategies,   the linear derivative $\frac{\partial V_{i}}{\partial u_{i}}$ of $V_{i}$ with respect to $\mathcal{A}_{i}$ is represented by the sensitivity process $Y^{u;u_{i}^{\prime }}$, as the derivative (in the $L^{2}$ sense) of the state $\mathbf{X}^{u}$ when player $i$ perturbs her control in the direction $u_{i}^{\prime }$. In the particular case where the diffusion term (as in SDE (3.1) in \cite{GLZ1}) is independent of both the state and control variables,  $Y^{u;u_{i}^{\prime }}$
satisfies a linear SDE:
\begin{equation}
\left\{\begin{aligned}
\mathrm{d}Y_{t,i}^{u;u_{i}^{\prime }} 
= ~ & \Big [\left( \partial
_{x}b_{i}\right) \left( t,X_{t,i}^{\mathbf{u}},\mathbf{X}_{t}^{\mathbf{u}%
},u_{t,i}\right) Y_{t,i}^{u;u_{i}^{\prime }}+\sum_{j=1}^{N}\left( \partial
_{y_{j}}b_{i}\right) \left( t,X_{t,i}^{\mathbf{u}},\mathbf{X}_{t}^{\mathbf{u}%
},u_{t,i}\right) Y_{t,i}^{u;u_{i}^{\prime }} \\
& +\left( \partial _{u}b_{i}\right) \left( t,X_{t,i}^{\mathbf{u}},\mathbf{%
X}_{t}^{\mathbf{u}},u_{t,i}\right) \delta _{h,i}u_{t,h}^{\prime }\Big ]%
\mathrm{d}t, \\
Y_{0,i}^{u;u_{i}^{\prime }} = ~ & 0,\text{ }\forall i\in I_{N}.
\end{aligned}\right.   \label{ry}
\end{equation}
Meanwhile,  the second-order linear derivative of $V_{i}$ can be characterized by a second-order sensitivity processes
of $Z^{u;u_{i}^{\prime },u_{j}^{\prime \prime }}$. This derivative also satisfies a linear SDE without the diffusion control:
\begin{equation}
\left\{\begin{aligned}
\mathrm{d}Z_{t,i}^{u;u_{i}^{\prime },u_{j}^{\prime \prime }} 
= ~ & \Big [\left( \partial _{x}b_{i}\right) \left( t,X_{t,i}^{\mathbf{u}},\mathbf{X}%
_{t}^{\mathbf{u}},u_{t,i}\right) Z_{t,i}^{u;u_{i}^{\prime },u_{j}^{\prime
\prime }}+\sum_{j=1}^{N}\left( \partial _{y_{j}}b_{i}\right) \left(
t,X_{t,i}^{\mathbf{u}},\mathbf{X}_{t}^{\mathbf{u}},u_{t,i}\right)
Z_{t,j}^{u;u_{i}^{\prime },u_{j}^{\prime \prime }} \\
& +\mathfrak{f}_{t,i}^{\mathbf{u},u_{h}^{\prime },u_{\ell }^{\prime
\prime }}\Big ]\mathrm{d}t, \\
Z_{0,i}^{u;u_{i}^{\prime },u_{j}^{\prime \prime }} 
= ~ & 0,\text{ }\forall i\in I_{N}, 
\end{aligned}\right.   \label{rz}
\end{equation}
where 
\begin{equation*}
\mathfrak{f}_{t,i}^{\mathbf{u,}u_{h}^{\prime },u_{\ell }^{\prime \prime }}:=\left(
\begin{array}{c}
Y_{t,i}^{\mathbf{u,}u_{h}^{\prime }} \\
\mathbf{Y}_{t}^{\mathbf{u,}u_{h}^{\prime }}%
\end{array}%
\right) ^{\top }\left(
\begin{array}{cc}
\partial _{xx}^{2}b_{i} & \partial _{xy}^{2}b_{i} \\
\partial _{yx}^{2}b_{i} & \partial _{yy}^{2}b_{i}%
\end{array}%
\right) \left( t,X_{t,i}^{\mathbf{u}},\mathbf{X}_{t}^{\mathbf{u}}\right)
\left(
\begin{array}{c}
Y_{t,i}^{\mathbf{u,}u_{\ell }^{\prime \prime }} \\
\mathbf{Y}_{t}^{\mathbf{u,}u_{\ell }^{\prime \prime }}%
\end{array}%
\right).
\end{equation*}
Subsequently, the linear derivatives of $V_{i}$ in (\ref{cost}) can be represented by means of
the sensitivity processes satisfying (\ref{v1}) (see (5.3) in \cite{GLZ1}).
\begin{lemma}
For all $\forall i,h\in I_{N}$, define the map $\frac{\partial V_{i}%
}{\partial u_{h}}:\mathcal{A}^{\left( N\right) }\times \mathcal{H}^{2}\left(
\mathbb{R}\right) \rightarrow \mathbb{R}$ such that for all $\mathbf{u\in }%
\mathcal{A}^{\left( N\right) }$ and $u_{h}^{\prime }\in \mathcal{H}%
^{2}\left( \mathbb{R}\right) $, then
\begin{equation}
\frac{\delta V_{i}}{\delta u_{h}}\left( \mathbf{u;}u_{h}^{\prime }\right) =%
\mathbb{E}\left[ \int_{0}^{T}\left(
\begin{array}{c}
\mathbf{Y}_{t}^{\mathbf{u},u_{h}^{\prime }} \\
u_{t,h}^{\prime }%
\end{array}%
\right) ^{\top }\left(
\begin{array}{c}
\partial _{x}f_{i} \\
\partial _{u_{h}}f_{i}%
\end{array}%
\right) \left( t,\mathbf{X}_{t}^{\mathbf{u}},\mathbf{u}_{t}\right) \mathrm{d}%
t+\left( \partial _{x}g_{i}\right) ^{\top }\left( \mathbf{X}_{T}^{\mathbf{u}}\right) \mathbf{Y}_{T}^{\mathbf{u},u_{h}^{\prime }}\right]. 
\label{linear1}
\end{equation}
\end{lemma}

The expression in (\ref{linear1}) resembles the variational inequality in stochastic  control problems, which in the classical control theory has been utilized to derive the maximum principle. 
Indeed, these processes are also known as the variational equations in control theory \cite{YZ} to characterize the variation of system trajectories under small disturbances of the state dynamics with respect to controls. 
This observation inspires us to reformulate the linear derivatives of \( V_i \) using BSDEs (such as Eq. (\ref{potenfv})), after obtaining a BSDE representation of $Y$
as in Theorem \ref{Thm:Y}.  


\paragraph{Basics of BSDE.} To this end, let us briefly review the most relevant results on BSDE. Consider the  following nonlinear BSDE introduced in \cite{PP1}

\begin{equation}
y_{t}=\xi +\int_{t}^{T}\mathsf{g}_{s}\left( y_{s},z_{s}\right) \mathrm{d}%
s-\int_{t}^{T}z_{s}\mathrm{d}W_{s},  \label{bsde0}
\end{equation}%
where a $d$-dimensional standard Brownian motion $W(\cdot )$ is defined on $(\Omega ,\mathcal{F},\mathbb{F}=\left( \mathcal{F}%
_{t}\right) _{0\leq t\leq T},\mathbb{P})$.  $\mathsf{g}_{t}\left(y,z\right) : [ 0,T]\times \mathbb{R}^{m}\times \mathbb{R}^{m\times d}\rightarrow \mathbb{R}^{m}$, and for every $(y,z)\in \mathbb{R}^{m}\times \mathbb{R}^{m\times d}$,
$\mathsf{g}_{t}(y,z)$ is an $\mathcal{F}_{t}$-adapted $\mathbb{R}^{m}$-valued process
with
\begin{equation}
\int_{0}^{T}\left\vert \mathsf{g}_{s}\left(0,0\right) \right\vert \mathrm{d}%
s\in L^{2}\left( \Omega ,\mathcal{F}_{T},P;\mathbb{R}\right) .  \label{bsde1}
\end{equation}%
In addition, $\mathsf{g}$ is Lipschitz continuous with
respect to $(y,z)$: there exists a constant $C>0$ such that, for all $%
y,y^{\prime }\in \mathbb{R}^{m},$ $z,z^{\prime }\in \mathbb{R}^{m\times d}$,%
\begin{equation}
\left\vert \mathsf{g}_{t}\left(y,z\right) -\mathsf{g}_{t}\left(y^{\prime },z^{\prime }\right)
\right\vert \leq C\left( \left\vert y-y^{\prime }\right\vert +\left\vert
z-z^{\prime }\right\vert \right) .  \label{bsde2}
\end{equation}

\begin{proposition}
\label{bs1}Suppose that $\mathsf{g}$ satisfies the conditions (\ref{bsde1}) and (\ref{bsde2}). Then,
for any given terminal condition $\xi \in L_{\mathcal{F}%
_{T}}^{2}\left( \mathbb{R}^{m}\right) $, there exists a unique pair of $\mathcal{F}_{t}$%
-adapted processes $(y,z)\in \mathcal{M}(0,T;\mathbb{R}^{m}\times \mathbb{R}%
^{m\times d})$ satisfying BSDE (\ref{bsde0}).
\end{proposition}

The fundamental difference between ODE  and  BSDE with similar Lipschitz condition and  the initial value is that  a solution for BSDE needs to be adapted to the given filtration with a terminal constraint: simply reversing time does not yield a solution for the terminal value problem of SDEs, as it would violate adaptiveness.  The (adapted) solution of a BSDE consists of a pair of adapted stochastic processes. The second component corrects the potential ``non-adaptiveness" arising from the backward nature of the equations, including the prescribed terminal value of the first component. Intuitively speaking,  the first component represents the ``mean evolution" of the dynamics, while the second component signifies the uncertainty/risk between the current time and the terminal time, while ensuring its adaptiveness.

For the purpose of our analysis, we consider a linear BSDE, with its estimates detailed in Appendix \ref{APP}.

\begin{lemma}
\label{bs2}Consider the following linear BSDE:
\begin{equation}
y_{t}=\xi +\int_{t}^{T}\left(
A_{s}y_{s}+\sum_{j=1}^{d}B_{s}^{j}z_{s}^{j}+f_{s}\right) \mathrm{d}%
s-\int_{t}^{T}z_{s}\mathrm{d}W_{s},  \label{lbsde3}
\end{equation}%
where $A,$ $B^{1},\ldots ,B^{d}$ are bounded $\mathbb{R}^{m\times m}$-valued
$\{\mathcal{F}_{t}\}_{t>0}$-adapted processes, $f_{s}\in \mathcal{M}\left(
\left( 0,T\right) ;\mathbb{R}^{m}\right) $ and $\xi \in L_{\mathcal{F}%
_{T}}^{2}\left( \mathbb{R}^{m}\right) $. Then, there exists some constant $%
C>0$ such that%
\begin{equation*}
\mathbb{E}\left[ \sup_{0\leq t\leq T}\left\vert y_{t}\right\vert
^{2}+\sum_{j=1}^{d}\int_{0}^{T}\left\vert z_{s}^{j}\right\vert ^{2}\mathrm{d}%
s\right] \leq {C}\mathbb{E}\Bigg [\left\vert \xi \right\vert
^{2}+\int_{t}^{T}\left\vert f_{s}\right\vert ^{2}\mathrm{d}s\Bigg ],
\end{equation*}%
where the constant $C$ depends on the bounds of $\left\vert A_{t}\right\vert
,\left\vert B_{t}^{1}\right\vert ,\ldots ,\left\vert B_{t}^{d}\right\vert $,
$T,$ and $d.$
\end{lemma}

The following lemma is necessary for our subsequent analysis and 
follows directly from It\^{o}'s formula.

\begin{lemma}
\label{l1}Let $\mathcal{Y},$ $\mathcal{P}\in \mathcal{M}(0,T;\mathbb{R}%
^{n\times n})$ $s$atisfy the following:%
\begin{equation*}
\left\{ \begin{aligned} \mathrm{d}\mathcal{Y}_{t} & = \Phi
_{t}\mathrm{dt+}\sum_{j=1}^{N}\Psi _{t}^{j}\mathrm{d}W_{t}^{j}, \\
\mathrm{d}\mathcal{P}_{t} & = \Theta _{t}\mathrm{dt+}
\sum_{j=1}^{N}\mathbf{Q}_{t}^{j}\mathrm{d}W_{t}^{j}, \end{aligned}\right.
\end{equation*}%
with $\Phi ,$ $\Psi ^{j},$ $\Theta $, $\mathbf{Q}^{j}$ all being elements in
$\mathcal{M}(0,T;R^{n\times n})$. Then
\begin{equation*}
\mathbb{E}\left\{ \text{\emph{tr}}\left[ \mathcal{P}_{t}\mathcal{Y}_{t}%
\right] -\text{\emph{tr}}\left[ \mathcal{P}_{0}\mathcal{Y}_{0}\right]
\right\} =\mathbb{E}\left[ \int_{0}^{T}\left\{ \text{\emph{tr}}\left[ \Theta
_{t}\mathcal{Y}_{t}+\mathcal{P}_{t}\Phi _{t}+\sum_{j=1}^{N}\mathbf{Q}%
_{t}^{j}\Psi _{t}^{j}\right] \right\} \mathrm{d}t\right] .
\end{equation*}
\end{lemma}

\section{BSDEs Representation of $\frac{\protect\delta V_{i}}{\protect\delta 
u_{h}}$ and $\frac{\protect\delta ^{2}V_{i}}{\protect\delta u_{h}\protect
\delta u_{\ell }}$}
\label{sect3} 
In this section, we investigate the general framework of
stochastic differential games as described by equations (\ref{sde1})-(\ref{cost}). We will revisit its associated sensitivity processes  in a BSDE framework and rewrite the derivative of $V_{i}$. As we will see, this BSDE approach avoids the estimation of the second order linear derivative of the $V_i$ via a duality representation of BSDE, and enables the estimation of $\alpha$ (additional) diffusion control, which is not feasible using the technique of \cite{GLZ1}.

\subsection{Sensitivity process $Y$}
Let us first recall the sensitivity of the controlled state with respect to a
single player's control. For each $\mathbf{u}\in \mathcal{H}^{2}\left(
\mathbb{R}\right) ^{N}$, let $\mathbf{X}^{\mathbf{u}}=\left( X_{i}^{\mathbf{u%
}}\right) _{i\in I_{N}}$ be the state process satisfying (\ref{sde1}). For
each $h\in I_{N}$ and $u_{h}^{\prime }\in \mathcal{H}^{2}\left( \mathbb{R}%
\right) $, define $\mathbf{Y}^{\mathbf{u},u_{h}^{\prime }}$ as the solution
of the following dynamics: for all $t\in \lbrack 0,T]$,
\begin{equation}
\left\{\begin{aligned} \mathrm{d}Y_{t,i}^{h} = ~ & \Big [\left( \partial
_{x}b_{i}\right) \left(
t,X_{t,i}^{\mathbf{u}},\mathbf{X}_{t}^{\mathbf{u}},u_{t,i}\right)
Y_{t,i}^{h}+\displaystyle\sum_{j=1}^{N}\left( \partial _{y_{j}}b_{i}\right)
\left( t,X_{t,i}^{\mathbf{u}},\mathbf{X}_{t}^{\mathbf{u}},u_{t,i}\right)
Y_{t,j}^{h} \\ & +\left( \partial _{u}b_{i}\right) \left(
t,X_{t,i}^{\mathbf{u}},\mathbf{X}_{t}^{\mathbf{u}},u_{t,i}\right) \delta
_{h,i}u_{t,h}^{\prime }\Big ]\mathrm{d}t \\ & +\Big [\left( \partial
_{x}\sigma _{i}\right) \left(
t,X_{t,i}^{\mathbf{u}},\mathbf{X}_{t}^{\mathbf{u}},u_{t,i}\right)
Y_{t,i}^{h}+\displaystyle\sum_{j=1}^{N}\left( \partial _{y_{j}}\sigma
_{i}\right) \left(
t,X_{t,i}^{\mathbf{u}},\mathbf{X}_{t}^{\mathbf{u}},u_{t,i}\right)
Y_{t,j}^{h} \\ & +\left( \partial _{u}\sigma _{i}\right) \left(
t,X_{t,i}^{\mathbf{u}},\mathbf{X}_{t}^{\mathbf{u}},u_{t,i}\right) \delta
_{h,i}u_{t,h}^{\prime }\Big]\mathrm{d}W_{t}^{i}, \\ Y_{0,i}^{h} = ~ & 0,
\quad \forall i\in I_{N}. \end{aligned}\right.  \label{v1}
\end{equation}
We point it out that (\ref{v1}) is also known as the first order variational equation in control theory
(see Yong and Zhou \cite{YZ}). That is, in the $L^{2}$ sense, $\mathbf{Y}^{%
\mathbf{u},u_{h}^{\prime }}$ is the derivative of the controlled state $%
\mathbf{X}_{t}^{\mathbf{u}}$ when player $h$ varies her control in the
direction of $u_{h}^{\prime}$. 

In the special case of no diffusion control, this
SDE becomes a linear SDE (\ref{ry}), and it is possible to estimate directly $\sup_{0\leq t\leq T}\mathbb{E%
}\left[ \left\vert Y_{t,i}^{u;u_{i}^{\prime }}\right\vert ^{p}\right] $ (see
Proposition 7.2 in \cite{GLZ1}). 
In the more general case  with a controlled diffusion term,   the sensitive process $Y_{t,i}^{h}$ for player $i$ depends on the actions of the other players through $\mathbf{X}_{t}^{\mathbf{u}}$, and the analysis becomes more challenging. 
We will instead adopt the technique of BSDE,
and study the adjoint equation for player $i$. 
For saving notions, let $b_{i}\left( t,\cdot \right) =b_{i}\left( t,X_{t,i}^{\mathbf{u}},\mathbf{%
X}_{t}^{\mathbf{u}},u_{t,i}\right) $, $\sigma _{i}\left( t,\cdot \right)
=\sigma _{i}\left( t,X_{t,i}^{\mathbf{u}},\mathbf{X}_{t}^{\mathbf{u}%
},u_{t,i}\right) $.

\begin{theorem}\label{Thm:Y}
The sensitivity process given in Equation (\ref{v1}) can be represented in the following form:
\begin{equation}
\left\{\begin{aligned}
\mathrm{d}\mathbf{Y}_{t}^{\mathbf{u},u_{h}^{\prime}} = & \left( \mathbf{B}%
_{0,x,y}\left( t,\mathbf{X}_{t}^{\mathbf{u}},\mathbf{u}_{t}\right) \mathbf{Y}%
_{t}^{\mathbf{u},u_{h}^{\prime }}+\mathbf{B}_{1,u}^{h}\left( t,\mathbf{X}%
_{t}^{\mathbf{u}},\mathbf{u}_{t}\right) u_{t,h}^{\prime }\right) \mathrm{d}t
\\
& +\displaystyle\sum_{j=1}^{N}\left( \mathbf{\Pi }^{j}_{0,x,y}\left( t,%
\mathbf{X}_{t}^{\mathbf{u}},\mathbf{u}_{t}\right) \mathbf{Y}_{t}^{\mathbf{u}%
,u_{h}^{\prime }}+\mathbf{\Pi }_{1,u}^{j}\left( t,\mathbf{X}_{t}^{\mathbf{u}%
},\mathbf{u}_{t}\right) \delta _{j,h}u_{t,h}^{\prime }\right) \mathrm{d}%
W_{t}^{j}\mathrm{,} \\
\mathbf{Y}_{0}^{\mathbf{u},u_{h}^{\prime }} = & \, 0, 
\end{aligned}\right.  \label{var1}
\end{equation}
where
\begin{eqnarray}
&&\mathbf{B}_{0,x,y}\left( t,\mathbf{X}_{t}^{\mathbf{u}},\mathbf{u}%
_{t}\right) =\mathbf{\bar{B}}_{x}\left( t,\mathbf{X}_{t}^{\mathbf{u}},%
\mathbf{u}_{t}\right) +\mathbf{B}_{0,y}\left( t,\mathbf{X}_{t}^{\mathbf{u}},%
\mathbf{u}_{t}\right) ,  \label{bxy} \\
&&\mathbf{\Pi }_{0,x,y}^{j}\left( t,\mathbf{X}_{t}^{\mathbf{u}},\mathbf{u}%
_{t}\right) =\mathbf{\bar{\Pi}}_{x}^{j}\left( t,\mathbf{X}_{t}^{\mathbf{u}},%
\mathbf{u}_{t}\right) +\mathbf{\Pi }_{0,y}^{j}\left( t,\mathbf{X}_{t}^{%
\mathbf{u}},\mathbf{u}_{t}\right) ,\quad \quad j\in I_{N},  \label{pixy}
\end{eqnarray}
with 
\begin{eqnarray}
\mathbf{\bar{B}}_{x}\left( t,\mathbf{X}_{t}^{\mathbf{u}},\mathbf{u}%
_{t}\right)  &=&\text{diag}\left( \left( \partial _{x}b_{1}\right) \left(
t,\cdot \right) ,\ldots ,\left( \partial _{x}b_{N}\right) \left( t,\cdot
\right) \right) _{N\times N},  \label{bbarx} \\
\mathbf{B}_{1,u}^{i}\left( t,\mathbf{X}_{t}^{\mathbf{u}},\mathbf{u}%
_{t}\right)  &=&\left[
\begin{array}{ccccc}
0, & \cdots , & \left( \partial _{u}b_{i}\right) \left( t,\cdot \right), &
\cdots , & 0\end{array}\right] _{N\times 1}^\top,\quad i\in I_{N},  \label{bi1u} \\
\mathbf{B}_{0,y}\left( t,\mathbf{X}_{t}^{\mathbf{u}},\mathbf{u}_{t}\right)
&=&\left[
\begin{array}{cccc}
\left( \partial _{y_{1}}b_{1}\right) \left( t,\cdot \right)  & \left(
\partial _{y_{2}}b_{1}\right) \left( t,\cdot \right)  & \cdots  & \left(
\partial _{y_{N}}b_{1}\right) \left( t,\cdot \right)  \\
\left( \partial _{y_{1}}b_{2}\right) \left( t,\cdot \right)  & \left(
\partial _{y_{2}}b_{2}\right) \left( t,\cdot \right)  & \cdots  & \left(
\partial _{y_{N}}b_{2}\right) \left( t,\cdot \right)  \\
\vdots  & \vdots  & \vdots  & \vdots  \\
\left( \partial _{y_{1}}b_{h}\right) \left( t,\cdot \right)  & \left(
\partial _{y_{2}}b_{h}\right) \left( t,\cdot \right)  & \cdots  & \left(
\partial _{y_{N}}b_{h}\right) \left( t,\cdot \right)  \\
\vdots  & \vdots  & \vdots  & \vdots  \\
\left( \partial _{y_{1}}b_{N}\right) \left( t,\cdot \right)  & \left(
\partial _{y_{2}}b_{N}\right) \left( t,\cdot \right)  & \cdots  & \left(
\partial _{y_{N}}b_{N}\right) \left( t,\cdot \right)
\end{array}%
\right] _{N\times N},  \label{b0y}
\end{eqnarray}
and
\begin{eqnarray}
\mathbf{\bar{\Pi}}_{x}^{i}\left( t,\mathbf{X}_{t}^{\mathbf{u}},\mathbf{u}%
_{t}\right)  &=&\left( \left( \partial _{x}\sigma _{i}\right) \left( t,\cdot
\right) \delta _{ij}\right) _{N\times N}^\top,  \label{pibarix} \\
\mathbf{\Pi }_{0,y}^{i}\left( t,\mathbf{X}_{t}^{\mathbf{u}},\mathbf{u}%
_{t}\right)  &=&\left[
\begin{array}{cccc}
0 & 0 & \cdots  & 0 \\
\vdots  & \vdots  & \vdots  & \vdots  \\
\left( \partial _{y_{1}}\sigma _{i}\right) \left( t,\cdot \right)  & \left(
\partial _{y_{2}}\sigma _{i}\right) \left( t,\cdot \right)  & \cdots  &
\left( \partial _{y_{N}}\sigma _{i}\right) \left( t,\cdot \right)  \\
\vdots  & \vdots  & \vdots  & \vdots  \\
0 & 0 & \cdots  & 0,%
\end{array}%
\right] _{N\times N},  \label{pii0y} \\
\mathbf{\Pi }_{1,u}^{j}\left( t,\mathbf{X}_{t}^{\mathbf{u}},\mathbf{u}%
_{t}\right)  &=&\left[
\begin{array}{ccccc}
0, & \cdots , & \left( \partial _{u}\sigma _{j}\right) \left( t,X_{t,j}^{%
\mathbf{u}},\mathbf{X}_{t}^{\mathbf{u}},u_{t,j}\right) , & \cdots , & 0%
\end{array}%
\right] _{N\times 1}^{\top },\quad j\in I_{N}.  \label{pij1u}
\end{eqnarray}

\end{theorem}

With this BSDE representation of $Y$, one can rewrite the first order derivative of $V_i$ via an associated BSDE.

\begin{proposition}
\label{pro1} Assumptions \emph{(A1)} and \emph{(A2)} hold. Then 
\begin{eqnarray}
&&\frac{\delta V_{i}}{\delta u_{h}}\left( \mathbf{u};u_{h}^{\prime }\right)
=\mathbb{E}\left[ \int_{0}^{T}\left(
\begin{array}{c}
\mathbf{Y}_{t}^{\mathbf{u},u_{h}^{\prime }} \\
u_{t,h}^{\prime}%
\end{array}%
\right) ^{\top }\left(
\begin{array}{c}
\partial _{x}f_{i} \\
\partial _{u_{h}}f_{i}%
\end{array}\right) \left( t,\mathbf{X}_{t}^{\mathbf{u}},\mathbf{u}_{t}\right) \mathrm{d}%
t+\left( \partial _{x}g_{i}\right) ^{\top }\left( \mathbf{X}_{T}^{\mathbf{u}%
}\right) \mathbf{Y}_{T}^{\mathbf{u},u_{h}^{\prime }}\right]  \notag \\
& & = \mathbb{E}\Bigg [\int_{0}^{T}\Big (P_{t,i}\mathbf{B}_{1,u}^{h}\left( t,%
\mathbf{X}_{t}^{\mathbf{u}},\mathbf{u}_{t}\right) +\left( \partial
_{u}\sigma _{h}\right) \left( t,X_{t,h}^{\mathbf{u}},\mathbf{X}_{t}^{\mathbf{%
u}},u_{t,h}\right) \left( Q_{t,i,h}\right) _{h}  +\partial _{u_{h}}f_{i}\Big )u_{t,h}^{\prime }\mathrm{d}t\Bigg ],
\label{potenfv}
\end{eqnarray}
where $\left( P_{t,i},\left\{ Q_{t,i,j}\right\} _{j=1,\ldots ,N}\right)_{0\leq t\leq T}$ is the unique adapted solution to the following BSDE
\begin{equation}
\left\{\begin{aligned}
-\mathrm{d}P_{t,i} = & \Big [\mathbf{B}_{0,x,y}\left( t,\mathbf{X}_{t}^{\mathbf{u}},\mathbf{u}_{t}\right) ^{\top }P_{t,i}+\sum_{j=1}^{N}\mathbf{\Pi }_{0,x,y}^{j}\left( t,\mathbf{X}_{t}^{\mathbf{u}},\mathbf{u}_{t}\right)^{\top }Q_{t,i,j}+\left( \partial _{x}f_{i}\right) \left( t,\mathbf{X}_{t}^{\mathbf{u}},\mathbf{u}_{t}\right) \Big ]\mathrm{d}t \\
& -\displaystyle\sum_{j=1}^{N}Q_{t,i,j}\mathrm{d}W_{t}^{j}, \\
P_{T,i} = & \, \partial _{x}g_{i}\left( \mathbf{X}_{T}^{\mathbf{u}}\right), \quad \quad \forall i\in I_{N}, 
\end{aligned}\right.  \label{adj1}
\end{equation}
where $\left( Q_{t,i,h}\right) _{h}$ represents the $h$-th component of the
column vector $Q_{t,i,h}$, while $\mathbf{B}_{0,x,y},\mathbf{\Pi}_{0,x,y}^{j}$, and $\mathbf{B}_{1,u}^{h}$ are given in (\ref{bxy}), (\ref{pixy}) and (\ref{bi1u}) in Theorem \ref{Thm:Y}.

\end{proposition}

\begin{remark}
In stochastic control theory, the dimensionality of the adjoint equation
(which is a BSDE) corresponds to that of the forward stochastic system.
However, in the context of potential games, the value function for the $i$%
-th player depends on $g_{i}\left( \mathbf{X}_{T}^{\mathbf{u}}\right)$,
indicating that the gradient of $g_{i}$ is represented as a column vector.
As a result, the adjoint equation associated with the $i$-th player is $N$ dimensional instead of one dimensional.
\end{remark}

\begin{remark}
Note that the adjoint equation (\ref{adj1}) remains
independent of $h$, since neither $B_{0,y}\left( t,\mathbf{X}_{t}^{\mathbf{u}},\mathbf{u}_{t}\right)$, nor $\Pi _{0,y}^{i}\left( t,\mathbf{X}_{t}^{\mathbf{u}},\mathbf{u}_{t}\right)$, nor $\left( \partial _{x}f_{i}\right)\left( t,\mathbf{X}_{t}^{\mathbf{u}},\mathbf{u}_{t}\right)$ depends on $h$.
\end{remark}

\subsection{Sensitivity Process $Z$}

We next proceed to quantify the sensitivity of the state process with respect to the changes in two players' controls. Let $u\in H^{2}\left(\mathbb{R}\right)^{N}$. For each $h$, $\ell
\in I_{N}$ with $h\neq \ell$, and each $u_{h}^{\prime }$, $u_{\ell}^{\prime \prime }\in H^{4}\left( \mathbb{R}\right) ^{N}$, define $Z^{\mathbf{u},u_{h}^{\prime },u_{\ell }^{\prime \prime }}\in S^{2}\left(\mathbb{R}^{N}\right) $ as the solution of the following dynamics: for all $i\in I_{N}$ and $t\in \lbrack 0,T]$,
\begin{equation}
\left\{\begin{aligned}
\mathrm{d}Z_{t,i}^{h,\ell } = & \Big [\left( \partial _{x}b_{i}\right)
\left( t,X_{t,i}^{\mathbf{u}},\mathbf{X}_{t}^{\mathbf{u}},u_{t,i}\right)
Z_{t,i}^{h,\ell }+\sum_{j=1}^{N}\left( \partial _{y_{j}}b_{i}\right) \left(
t,X_{t,i}^{\mathbf{u}},\mathbf{X}_{t}^{\mathbf{u}},u_{t,i}\right)
Z_{t,j}^{h,\ell } \\
&  +\mathfrak{f}_{1,t,i}^{\mathbf{u},u_{h}^{\prime },u_{\ell }^{\prime
\prime }}+\mathfrak{f}_{2,t,i}^{\mathbf{u},u_{h}^{\prime },u_{\ell }^{\prime
\prime }}\Big ]\mathrm{d}t \\
&  +\Big [\left( \partial _{x}\sigma _{i}\right) \left( t,X_{t,i}^{\mathbf{%
u}},\mathbf{X}_{t}^{\mathbf{u}},u_{t,i}\right) Z_{t,i}^{h,\ell
}+\sum_{j=1}^{N}\left( \partial _{y_{j}}\sigma _{i}\right) \left( t,X_{t,i}^{%
\mathbf{u}},\mathbf{X}_{t}^{\mathbf{u}},u_{t,i}\right) Z_{t,j}^{h,\ell } \\
&  +\mathfrak{g}_{1,t,i}^{\mathbf{u},u_{h}^{\prime },u_{\ell }^{\prime
\prime }}+\mathfrak{g}_{2,t,i}^{\mathbf{u},u_{h}^{\prime },u_{\ell }^{\prime
\prime }}\Big ]\mathrm{d}W_{t}^{i}, \\
Z_{0,i}^{h,\ell} = & \, 0, \quad \quad \forall i\in I_{N},%
\end{aligned}\right.  \label{v3}
\end{equation}
where $\mathfrak{f}_{1,t,i}^{\mathbf{u},u_{h}^{\prime },u_{\ell }^{\prime
\prime }}:\Omega \times \left[ 0,T\right] \rightarrow \mathbb{R}$ and $\mathfrak{g}_{1,t,i}^{\mathbf{u},u_{h}^{\prime },u_{\ell }^{\prime \prime}}:\Omega \times \left[ 0,T\right] \rightarrow \mathbb{R}$ are defined by
%
\begin{eqnarray}
\mathfrak{f}_{1,t,i}^{\mathbf{u},u_{h}^{\prime },u_{\ell }^{\prime \prime }}
& \!\!\!=\!\!\! & \left(\begin{array}{c}
Y_{t,i}^{\mathbf{u},u_{h}^{\prime }} \\
\mathbf{Y}_{t}^{\mathbf{u},u_{h}^{\prime }}%
\end{array}%
\right) ^{\top }\left(
\begin{array}{cc}
\partial _{xx}^{2}b_{i} & \partial _{xy}^{2}b_{i} \\
\partial _{yx}^{2}b_{i} & \partial _{yy}^{2}b_{i}%
\end{array}%
\right) \left( t,\cdot \right) \left(
\begin{array}{c}
Y_{t,i}^{\mathbf{u},u_{\ell }^{\prime \prime }} \\
\mathbf{Y}_{t}^{\mathbf{u},u_{\ell }^{\prime \prime }}%
\end{array}%
\right) \notag\\
& \!\!\!=\!\!\! & \text{tr}\left[ \partial _{yy}^{2}b_{i}\left( t,\cdot \right) \mathcal{Y}%
_{t}^{\mathbf{u},u_{h}^{\prime },u_{\ell }^{\prime \prime }}\right] +\Gamma
_{t,i}^{\mathbf{u},u_{h}^{\prime },u_{\ell }^{\prime \prime }}\left( \bar{b}%
\right) , \\
\mathfrak{g}_{1,t,i}^{\mathbf{u},u_{h}^{\prime },u_{\ell }^{\prime \prime }}
& \!\!\!=\!\!\! & \left(
\begin{array}{c}
Y_{t,i}^{\mathbf{u},u_{h}^{\prime }} \\
\mathbf{Y}_{t}^{\mathbf{u},u_{h}^{\prime }}%
\end{array}%
\right) ^{\top }\left(
\begin{array}{cc}
\partial _{xx}^{2}\sigma _{i} & \partial _{xy}^{2}\sigma _{i} \\
\partial _{yx}^{2}\sigma _{i} & \partial _{yy}^{2}\sigma _{i}%
\end{array}%
\right) \left( t,\cdot \right) \left(
\begin{array}{c}
Y_{t,i}^{\mathbf{u},u_{\ell }^{\prime \prime }} \\
\mathbf{Y}_{t}^{\mathbf{u},u_{\ell }^{\prime \prime }}%
\end{array}\right) \notag\\
& \!\!\!=\!\!\! & \text{tr}\left[ \partial _{yy}^{2}\sigma _{i}\left( t,\cdot \right)
\mathcal{Y}_{t}^{\mathbf{u},u_{h}^{\prime },u_{\ell }^{\prime \prime }}%
\right] +\Xi _{t,i}^{\mathbf{u},u_{h}^{\prime },u_{\ell }^{\prime \prime
}}\left( \bar{\sigma}\right),
\end{eqnarray}
where, for simplicity, we put $b_{i}\left( t,\cdot \right) =b_{i}\left( t,X_{t,i}^{\mathbf{u}},\mathbf{X}_{t}^{\mathbf{u}},u_{t,i}\right)$
and $\sigma_{i}\left( t,\cdot \right) = \sigma_{i}\left(t,X_{t,i}^{\mathbf{u}},\mathbf{X}_{t}^{\mathbf{u}},u_{t,i}\right)$, 
\begin{eqnarray}
\mathcal{Y}_{t}^{\mathbf{u},u_{h}^{\prime },u_{\ell }^{\prime \prime }} & \!\!\!=\!\!\! & 
\mathbf{Y}_{t}^{\mathbf{u},u_{\ell }^{\prime \prime }}\left( \mathbf{Y}_{t}^{%
\mathbf{u},u_{h}^{\prime }}\right) ^{\top },  \label{fy} \\
\bar{b}\left( t,\mathbf{X},\mathbf{u}\right) & \!\!\!=\!\!\! & \left( b_{1},\ldots,b_{N}\right) ^{\top }\left( t,\mathbf{X},\mathbf{u}\right),\text{ }\bar{%
\sigma}_{j}\left( t,\mathbf{X},\mathbf{u}\right) =\left( 0,\ldots \sigma
_{j},\ldots 0\right) ^{\top }\left( t,\mathbf{X},\mathbf{u}\right),
\label{bsg} \\
\Gamma_{t,i}^{\mathbf{u},u_{h}^{\prime },u_{\ell }^{\prime \prime }}\left(
\bar{b}\right) & \!\!\!=\!\!\! & Y_{t,i}^{\mathbf{u},u_{h}^{\prime }}\partial
_{xx}^{2}b_{i}\left( t,\cdot \right) Y_{t,i}^{\mathbf{u},u_{\ell }^{\prime
\prime }}+\left( \mathbf{Y}_{t}^{\mathbf{u},u_{h}^{\prime }}\right) ^{\top
}\partial _{yx}^{2}b_{i}\left( t,\cdot \right) Y_{t,i}^{\mathbf{u},u_{\ell
}^{\prime \prime }}  \notag \\
&&+Y_{t,i}^{\mathbf{u},u_{h}^{\prime }}\partial _{xy}^{2}b_{i}\left( t,\cdot
\right) \mathbf{Y}_{t}^{\mathbf{u},u_{\ell }^{\prime \prime }},  \label{gama}
\\
\Xi _{t,i}^{\mathbf{u},u_{h}^{\prime },u_{\ell }^{\prime \prime }}\left(
\bar{\sigma}\right) & \!\!\!=\!\!\! & Y_{t,i}^{\mathbf{u},u_{h}^{\prime }}\partial
_{xx}^{2}\sigma _{i}\left( t,\cdot \right) Y_{t,i}^{\mathbf{u},u_{\ell
}^{\prime \prime }}+\left( \mathbf{Y}_{t}^{\mathbf{u},u_{h}^{\prime
}}\right) ^{\top }\partial _{yx}^{2}\sigma _{i}\left( t,\cdot \right)
Y_{t,i}^{\mathbf{u},u_{\ell }^{\prime \prime }}  \notag \\
&&+Y_{t,i}^{\mathbf{u},u_{h}^{\prime }}\partial _{xy}^{2}\sigma _{i}\left(
t,\cdot \right) \mathbf{Y}_{t}^{\mathbf{u},u_{\ell }^{\prime \prime}}.
\label{str}
\end{eqnarray}
Therefore,
\begin{eqnarray*}
\left[
\begin{array}{ccccc}
\mathfrak{f}_{1,t,1}^{\mathbf{u},u_{h}^{\prime },u_{\ell }^{\prime \prime }},
& \cdots , & \mathfrak{f}_{1,t,i}^{\mathbf{u},u_{h}^{\prime },u_{\ell
}^{\prime \prime }}, & \cdots , & \mathfrak{f}_{1,t,N}^{\mathbf{u}%
,u_{h}^{\prime },u_{\ell }^{\prime \prime }}%
\end{array}%
\right] _{N\times 1}^{\top } &=&\partial _{yy}^{2}\bar{b}\left( t,\mathbf{X}%
_{t},\mathbf{u}_{t}\right) \star \mathcal{Y}_{t}^{\mathbf{u},u_{h}^{\prime
},u_{\ell }^{\prime \prime }}+\Gamma _{t}^{\mathbf{u},u_{h}^{\prime
},u_{\ell }^{\prime \prime }}\left( \bar{b}\right) , \\
\left[
\begin{array}{ccccc}
0, & \cdots , & \mathfrak{g}_{1,t,j}^{\mathbf{u},u_{h}^{\prime },u_{\ell
}^{\prime \prime }}, & \cdots , & 0%
\end{array}%
\right] _{N\times 1}^{\top } &=&\partial _{yy}^{2}\bar{\sigma}_{j}\left( t,%
\mathbf{X}_{t},\mathbf{u}_{t}\right) \star \mathcal{Y}_{t}^{\mathbf{u}%
,u_{h}^{\prime },u_{\ell }^{\prime \prime }}+\bar{\Xi}_{t,j}^{\mathbf{u}%
,u_{h}^{\prime },u_{\ell }^{\prime \prime }}\left( \bar{\sigma}\right) ,
\end{eqnarray*}%
where
\begin{equation*}
\partial _{yy}^{2}\bar{b}\left( t,\mathbf{X}_{t},\mathbf{u}_{t}\right) \star
\mathcal{Y}_{t}^{\mathbf{u},u_{h}^{\prime },u_{\ell }^{\prime \prime }}=%
\left[
\begin{array}{c}
\text{tr}\left[ \partial _{yy}^{2}b_{1}\left( t,\cdot \right) \mathcal{Y}%
_{t}^{\mathbf{u},u_{h}^{\prime },u_{\ell }^{\prime \prime }}\right]  \\
\vdots  \\
\text{tr}\left[ \partial _{yy}^{2}b_{i}\left( t,\cdot \right) \mathcal{Y}%
_{t}^{\mathbf{u},u_{h}^{\prime },u_{\ell }^{\prime \prime }}\right]  \\
\vdots  \\
\text{tr}\left[ \partial _{yy}^{2}b_{N}\left( t,\cdot \right) \mathcal{Y}%
_{t}^{\mathbf{u},u_{h}^{\prime },u_{\ell }^{\prime \prime }}\right]
\end{array}%
\right]_{N\times 1},\quad \Gamma _{t}^{\mathbf{u},u_{h}^{\prime },u_{\ell }^{\prime
\prime }}\left( \bar{b}\right) =\left[
\begin{array}{c}
\Gamma _{t,1}^{\mathbf{u},u_{h}^{\prime },u_{\ell }^{\prime \prime }}\left(
\bar{b}\right)  \\
\vdots  \\
\Gamma _{t,i}^{\mathbf{u},u_{h}^{\prime },u_{\ell }^{\prime \prime }}\left(
\bar{b}\right)  \\
\vdots  \\
\Gamma _{t,N}^{\mathbf{u},u_{h}^{\prime },u_{\ell }^{\prime \prime }}\left(
\bar{b}\right)
\end{array}\right]_{N\times 1}
\end{equation*}
and
\begin{eqnarray*}
\partial _{yy}^{2}\bar{\sigma}_{j}\left( t,\mathbf{X}_{t},\mathbf{u}%
_{t}\right) \star \mathcal{Y}_{t}^{\mathbf{u},u_{h}^{\prime },u_{\ell
}^{\prime \prime }} &=&\left[
\begin{array}{ccccc}
0, & \ldots , & \text{tr}\left[ \partial _{yy}^{2}\sigma _{j}\left( t,\cdot
\right) \mathcal{Y}_{t}^{\mathbf{u},u_{h}^{\prime },u_{\ell }^{\prime \prime
}}\right] , & \cdots , & 0%
\end{array}%
\right] ^{\top }_{N\times 1}, \\
\quad \bar{\Xi}_{t,j}^{\mathbf{u},u_{h}^{\prime },u_{\ell }^{\prime \prime
}}\left( \bar{\sigma}\right)  &=&\left[
\begin{array}{ccccc}
0, & \cdots , & \Xi _{t,j}^{\mathbf{u},u_{h}^{\prime },u_{\ell }^{\prime
\prime }}\left( \bar{\sigma}\right) , & \cdots , & 0%
\end{array}\right] ^{\top }_{N\times 1}.
\end{eqnarray*}
Here
\begin{eqnarray}
\mathfrak{f}_{2,t,i}^{\mathbf{u},u_{h}^{\prime },u_{\ell }^{\prime \prime }}
&=&\,u_{t,h}^{\prime }\cdot \left( \partial _{ux}^{2}b_{i},\partial
_{uy}^{2}b_{i}\right) \left( t,\cdot \right) \cdot \left(
\begin{array}{c}
Y_{t,i}^{\mathbf{u},u_{\ell }^{\prime \prime }} \\
\mathbf{Y}_{t}^{\mathbf{u},u_{\ell }^{\prime \prime }}%
\end{array}%
\right) \delta _{h,i}  \notag \\
&&+u_{t,\ell }^{\prime \prime }\cdot \left( \partial _{ux}^{2}b_{i},\partial
_{uy}^{2}b_{i}\right) \left( t,\cdot \right) \cdot \left(
\begin{array}{c}
Y_{t,i}^{\mathbf{u},u_{h}^{\prime }} \\
\mathbf{Y}_{t}^{\mathbf{u},u_{h}^{\prime }}%
\end{array}%
\right) \delta _{\ell ,i}  \notag \\
&=&\,u_{t,h}^{\prime }\left( \partial _{ux}b_{i}\left( t,\cdot \right) \cdot
Y_{t,i}^{\mathbf{u},u_{\ell }^{\prime \prime }}+\partial
_{uy}^{2}b_{i}\left( t,\cdot \right) \cdot \mathbf{Y}_{t}^{\mathbf{u}%
,u_{\ell }^{\prime \prime }}\right) \delta _{h,i}  \notag \\
&&+u_{t,\ell }^{\prime \prime }\left( \partial _{xu}^{2}b_{i}\left( t,\cdot
\right) \cdot Y_{t,i}^{\mathbf{u},u_{h}^{\prime }}+\partial
_{yu}^{2}b_{i}\left( t,\cdot \right) \cdot \mathbf{Y}_{t}^{\mathbf{u}%
,u_{h}^{\prime }}\right) \delta _{\ell ,i},  \label{f2}
\end{eqnarray}
\begin{eqnarray}
\mathfrak{g}_{2,t,i}^{\mathbf{u},u_{h}^{\prime },u_{\ell }^{\prime \prime }}
&=&\,u_{t,h}^{\prime }\cdot \left( \partial _{ux}^{2}\sigma _{i},\partial
_{uy}^{2}\sigma _{i}\right) \left( t,\cdot \right) \cdot \left(
\begin{array}{c}
Y_{t,i}^{\mathbf{u},u_{\ell }^{\prime \prime }} \\
\mathbf{Y}_{t}^{\mathbf{u},u_{\ell }^{\prime \prime }}%
\end{array}%
\right) \delta _{h,i}  \notag \\
&&+u_{t,\ell }^{\prime \prime }\cdot \left( \partial _{ux}^{2}\sigma
_{i},\partial _{uy}^{2}\sigma _{i}\right) \left( t,\cdot \right) \cdot
\left(
\begin{array}{c}
Y_{t,i}^{\mathbf{u},u_{h}^{\prime }} \\
\mathbf{Y}_{t}^{\mathbf{u},u_{h}^{\prime }}%
\end{array}%
\right) \delta _{\ell ,i}  \notag \\
&=&\,u_{t,h}^{\prime }\left( \partial _{ux}^{2}\sigma _{i}\left( t,\cdot
\right) \cdot Y_{t,i}^{\mathbf{u},u_{\ell }^{\prime \prime }}+\partial
_{uy}^{2}\sigma _{i}\left( t,\cdot \right) \cdot \mathbf{Y}_{t}^{\mathbf{u}%
,u_{\ell }^{\prime \prime }}\right) \delta _{h,i}  \notag \\
&&+u_{t,\ell }^{\prime \prime }\left( \partial _{xu}^{2}\sigma _{i}\left(
t,\cdot \right) \cdot Y_{t,i}^{\mathbf{u},u_{h}^{\prime }}+\partial
_{yu}^{2}\sigma _{i}\left( t,\cdot \right) \cdot \mathbf{Y}_{t}^{\mathbf{u}%
,u_{h}^{\prime }}\right) \delta _{\ell ,i}.  \label{g2}
\end{eqnarray}
Meanwhile, assume $h<\ell $ without loss of generality,
\begin{equation}
\mathfrak{F}_{t}^{\mathbf{u},u_{h}^{\prime },u_{\ell }^{\prime \prime }}=%
\left[
\begin{array}{c}
\mathfrak{f}_{2,t,1}^{\mathbf{u},u_{h}^{\prime },u_{\ell }^{\prime \prime }}
\\
\vdots  \\
\mathfrak{f}_{2,t,h}^{\mathbf{u},u_{h}^{\prime },u_{\ell }^{\prime \prime }}
\\
\vdots  \\
\mathfrak{f}_{2,t,\ell }^{\mathbf{u},u_{h}^{\prime },u_{\ell }^{\prime
\prime }} \\
\vdots  \\
\mathfrak{f}_{2,t,N}^{\mathbf{u},u_{h}^{\prime },u_{\ell }^{\prime \prime }}%
\end{array}%
\right] _{N\times 1}=\left[
\begin{array}{c}
0 \\
\vdots  \\
u_{t,h}^{\prime }\left( \partial _{ux}^{2}b_{h}\left( t,\cdot \right) \cdot
Y_{t,h}^{\mathbf{u},u_{\ell }^{\prime \prime }}+\partial
_{uy}^{2}b_{h}\left( t,\cdot \right) \cdot \mathbf{Y}_{t}^{\mathbf{u}%
,u_{\ell }^{\prime \prime }}\right)  \\
\vdots  \\
u_{t,\ell }^{\prime \prime }\left( \partial _{ux}^{2}b_{\ell }\left( t,\cdot
\right) \cdot Y_{t,\ell }^{\mathbf{u},u_{h}^{\prime }}+\partial
_{uy}^{2}b_{\ell }\left( t,\cdot \right) \cdot \mathbf{Y}_{t}^{\mathbf{u}%
,u_{h}^{\prime }}\right)  \\
\vdots  \\
0%
\end{array}%
\right] _{N\times 1}  \label{ff}
\end{equation}%
and

\begin{equation}
\mathfrak{G}_{t,j}^{\mathbf{u},u_{h}^{\prime },u_{\ell }^{\prime \prime }}=%
\left[
\begin{array}{ccccc}
0, & \ldots , & \mathfrak{g}_{2,t,j}^{\mathbf{u},u_{h}^{\prime },u_{\ell
}^{\prime \prime }} & \cdots , & 0%
\end{array}%
\right] _{N\times 1}^{\top }.  \label{fb}
\end{equation}

\begin{remark}
In stochastic control theory, $%
\mathbf{Z}^{\mathbf{u},u_{h}^{\prime },u_{\ell }^{\prime \prime }}$
corresponds to the second-order variational equation. Moreover, when
the control domain is assumed to be compact and involves diffusion control, the second-order variational equation is
required through the use of the spike variation technique (see Yong and Zhou \cite{YZ}).
\end{remark}

\begin{remark}
Observe that $h \neq \ell$. Consequently, when analyzing the policy $\mathbf{%
u}^{\varepsilon} = \left( u_{\ell} + \varepsilon u_{\ell}^{\prime \prime},
u_{-\ell} \right)$, the $h$-th component of $\mathbf{u}^{\varepsilon}$
remains $u_h$. This results in the absence of $\partial_{uu}^2 b_i\left(t,
\cdot \right)$ and $\partial_{uu}^2 \sigma_i\left(t, \cdot \right)$ for all $%
i \in I_N$.
\end{remark}

Based on (\ref{v3}), we are able to present
\begin{equation}
\left\{\begin{aligned}
\mathrm{d}\mathbf{Z}_{t}^{\mathbf{u},u_{h}^{\prime },u_{\ell }^{\prime
\prime }} = & \Big [\mathbf{B}_{0,x,y}\left( t,\mathbf{X}_{t}^{\mathbf{u}},%
\mathbf{u}_{t}\right) \mathbf{Z}_{t}^{\mathbf{u},u_{h}^{\prime },u_{\ell
}^{\prime \prime }}+\mathfrak{F}_{t}^{\mathbf{u},u_{h}^{\prime },u_{\ell
}^{\prime \prime }} \\
&  +\partial _{yy}^{2}\bar{b}\left( t,\mathbf{X}_{t},\mathbf{u}_{t}\right)
\star \mathcal{Y}_{t}^{\mathbf{u},u_{h}^{\prime },u_{\ell }^{\prime \prime
}}+\Gamma _{t}^{\mathbf{u},u_{h}^{\prime },u_{\ell }^{\prime \prime }}\left(
\bar{b}\right) \Big ]\mathrm{d}t \\
&  +\sum_{j=1}^{N}\Big [\mathbf{\Pi }_{0,x,y}^{j}\left( t,\mathbf{X}_{t}^{%
\mathbf{u}},\mathbf{u}_{t}\right) \mathbf{Z}_{t}^{\mathbf{u},u_{h}^{\prime
},u_{\ell }^{\prime \prime }}+\mathfrak{G}_{t,j}^{\mathbf{u},u_{h}^{\prime
},u_{\ell }^{\prime \prime }} \\
&  +\partial _{yy}^{2}\bar{\sigma}_{j}\left( t,\mathbf{X}_{t},\mathbf{u}%
_{t}\right) \star \mathcal{Y}_{t}^{\mathbf{u},u_{h}^{\prime },u_{\ell
}^{\prime \prime }}+\bar{\Xi}_{t,j}^{\mathbf{u},u_{h}^{\prime },u_{\ell
}^{\prime \prime }}\left( \bar{\sigma}\right) \Big ]\mathrm{d}W_{t}^{j}, \\
\mathbf{Z}_{0}^{\mathbf{u},u_{h}^{\prime }u_{\ell }^{\prime \prime }} = &
\, 0, 
\end{aligned}\right.   \label{var2}
\end{equation}
where $\mathfrak{F}_{t}^{\mathbf{u},u_{h}^{\prime },u_{\ell }^{\prime \prime
}},\Gamma _{t}^{\mathbf{u},u_{h}^{\prime },u_{\ell }^{\prime \prime }}\left(
\bar{b}\right) ,\mathfrak{G}_{t,j}^{\mathbf{u},u_{h}^{\prime },u_{\ell
}^{\prime \prime }}$ and $\bar{\Xi}_{t,j}^{\mathbf{u},u_{h}^{\prime
},u_{\ell }^{\prime \prime }}\left( \bar{\sigma}\right) $ are defined in (%
\ref{ff}), (\ref{gama}), (\ref{fb}), (\ref{str}), respectively.
Consider the following $\mathbb{S}^{n}$-value second-order adjoint BSDE:
\begin{equation}
\left\{\begin{aligned}
-\mathrm{d}\mathcal{P}_{t,i}^{h,\ell } & = \digamma _{t,i}^{h,\ell }%
\mathrm{d}t-\sum_{j=1}^{N}\mathcal{Q}_{t,i,j}^{h,\ell }\mathrm{d}W_{t}^{j}, \\
\mathcal{P}_{T,i}^{h,\ell } & = \partial _{xx}g_{i}\left( \mathbf{X}_{T}^{%
\mathbf{u}}\right), \quad \quad \forall i\in I_{N}, 
\end{aligned}\right.  \label{adj2}
\end{equation}
where
\begin{eqnarray}
\digamma _{t,i}^{h,\ell } &=&\mathbf{B}_{0,x,y}\left( t,\mathbf{X}_{t}^{%
\mathbf{u}},\mathbf{u}_{t}\right) ^{\top }\mathcal{P}_{t,i}^{h,\ell }+%
\mathcal{P}_{t,i}^{h,\ell }\mathbf{B}_{0,x,y}\left( t,\mathbf{X}_{t}^{%
\mathbf{u}},\mathbf{u}_{t}\right)  \notag \\
&&+\sum_{i=1}^{N}\mathbf{\Pi }_{0,x,y}^{i}\left( t,\mathbf{X}_{t}^{\mathbf{u}%
},\mathbf{u}_{t}\right) ^{\top }\mathcal{P}_{t,i}^{h,\ell }\mathbf{\Pi }%
_{0,x,y}^{i}\left( t,\mathbf{X}_{t}^{\mathbf{u}},\mathbf{u}_{t}\right)
\notag \\
&&+\sum_{j=1}^{N}\left\{ \mathbf{\Pi }_{0,x,y}^{j}\left( t,\mathbf{X}_{t}^{%
\mathbf{u}},\mathbf{u}_{t}\right) ^{\top }\mathcal{Q}_{t,i,j}^{h,\ell }+%
\mathcal{Q}_{t,i,j}^{h,\ell }\mathbf{\Pi }_{0,x,y}^{j}\left( t,\mathbf{X}_{t}^{%
\mathbf{u}},\mathbf{u}_{t}\right) \right\}  \notag \\
&&+\partial _{yy}^{2}f_{i}+\partial _{yy}^{2}\left\langle P_{t,i},\bar{b}%
\left( t,\mathbf{X}_{t},\mathbf{u}_{t}\right) \right\rangle
+\sum_{j=1}^{N}\partial _{yy}^{2}\left\langle Q_{t,i,j},\bar{\sigma}_{j}\left(
t,\mathbf{X}_{t},\mathbf{u}_{t}\right) \right\rangle.  \label{c3}
\end{eqnarray}

Now we are able to express the second order derivative of $V_i$ as the following. 

\begin{proposition}
\label{V2} Assumptions \emph{(A1)} and \emph{(A2)} hold. Then
\begin{eqnarray}
&&\frac{\delta ^{2}V_{i}}{\delta u_{h}\delta u_{\ell }}\left( \mathbf{u}%
;u_{h}^{\prime },u_{\ell }^{\prime \prime }\right)   \notag \\
&=\!\!\! & \mathbb{E}\Bigg \{\int_{0}^{T}\text{tr}\Bigg [\Bigg (\mathcal{P}%
_{t,i}^{h,\ell }\mathbf{Y}_{t}^{\mathbf{u},u_{\ell }^{\prime \prime }}%
\mathbf{B}_{1,u}^{h}\left( t,\mathbf{X}_{t}^{\mathbf{u}},\mathbf{u}%
_{t}\right) ^{\top }+\mathcal{P}_{t,i}^{h,\ell }\mathbf{\Pi }%
_{0,x,y}^{h}\left( t,\mathbf{X}_{t}^{\mathbf{u}},\mathbf{u}_{t}\right)
\mathbf{Y}_{t}^{\mathbf{u},u_{\ell }^{\prime \prime }}\mathbf{\Pi }%
_{1,u}^{h}\left( t,\mathbf{X}_{t}^{\mathbf{u}},\mathbf{u}_{t}\right) ^{\top }
\notag \\
&&+\mathcal{Q}_{t,i,h}^{h,\ell }\cdot \mathbf{Y}_{t}^{\mathbf{u},u_{\ell
}^{\prime \prime }}\mathbf{\Pi }_{1,u}^{h}\left( t,\mathbf{X}_{t}^{\mathbf{u}%
},\mathbf{u}_{t}\right) ^{\top }+\partial _{u_{h}x}^{2}f_{i}\cdot \mathbf{Y}%
_{t}^{\mathbf{u},u_{\ell }^{\prime \prime }}\Bigg )u_{t,h}^{\prime }  \notag
\\
&&+\Bigg (\mathcal{P}_{t,i}^{h,\ell }\mathbf{B}_{1,u}^{\ell }\left( t,%
\mathbf{X}_{t}^{\mathbf{u}},\mathbf{u}_{t}\right) \left( \mathbf{Y}_{t}^{%
\mathbf{u},u_{h}^{\prime }}\right) ^{\top }+\mathcal{P}_{t,i}^{h,\ell }%
\mathbf{\Pi }_{1,u}^{\ell }\left( t,\mathbf{X}_{t}^{\mathbf{u}},\mathbf{u}%
_{t}\right) \left( \mathbf{Y}_{t}^{\mathbf{u},u_{h}^{\prime }}\right) ^{\top
}\mathbf{\Pi }_{0,x,y}^{\ell }\left( t,\mathbf{X}_{t}^{\mathbf{u}},\mathbf{u}%
_{t}\right) ^{\top }  \notag \\
&&+\mathcal{Q}_{t,i,\ell }^{h,\ell }\cdot \mathbf{\Pi }_{1,u}^{\ell }\left( t,%
\mathbf{X}_{t}^{\mathbf{u}},\mathbf{u}_{t}\right) \left( \mathbf{Y}_{t}^{%
\mathbf{u},u_{h}^{\prime }}\right) ^{\top }+\left( \mathbf{Y}_{t}^{\mathbf{u}%
,u_{h}^{\prime }}\right) ^{\top }\cdot \partial _{xu_{\ell }}^{2}f_{i}\Bigg )%
u_{t,\ell }^{\prime \prime }  \notag \\
&&+\partial _{u_{h}u_{\ell }}^{2}f_{i}\cdot u_{t,h}^{\prime }\cdot u_{t,\ell
}^{\prime \prime }+\left\langle P_{t,i},\mathfrak{F}_{t}^{\mathbf{u}%
,u_{h}^{\prime },u_{\ell }^{\prime \prime }}+\Gamma _{t}^{\mathbf{u}%
,u_{h}^{\prime },u_{\ell }^{\prime \prime }}\left( \bar{b}\right)
\right\rangle   \notag \\
&&+\sum_{j=1}^{N}\left\langle Q_{t,i,j},\mathfrak{G}_{t,j}^{\mathbf{u}%
,u_{h}^{\prime },u_{\ell }^{\prime \prime }}+\bar{\Xi}_{t,j}^{\mathbf{u}%
,u_{h}^{\prime },u_{\ell }^{\prime \prime }}\left( \bar{\sigma}\right)
\right\rangle \Bigg ]\mathrm{d}t\Bigg \},  \label{Vi}
\end{eqnarray}%
where $\mathbf{B}_{1,u}^{\theta },\mathbf{\Pi }_{0,x,y}^{\theta },\mathbf{%
\Pi }_{1,u}^{\theta },\theta =h,\ell $\textbf{\ }are defined in Theorem \ref%
{Thm:Y}, $\left( P_{t,i},\left\{ Q_{t,i,j}\right\} _{j=1,\ldots ,N}\right) $
and $\left( \mathcal{P}_{t,i}^{h,\ell },\left\{ \mathcal{Q}_{t,i,j}^{h,\ell
}\right\} _{j=1,\ldots ,N}\right) $ are solutions to BSDEs (\ref{adj1}) and (\ref%
{adj2}), respectively.
\end{proposition}
The above expression will be crucial for our subsequent analysis of $\alpha$,  especially for avoiding  the estimation of the second order sensitivity process $Z.$


\section{Estimation of $\alpha$}

\label{sect4}

  In this section, we derive an upper bound for $\alpha$ defined in (\ref{aph}) for the general open-loop framework (\ref{sde1})-(\ref{cost}), by estimating $X_{t,i}^{\mathbf{u}}$ and $\mathbf{Y}_{t}^{\mathbf{u},u_{\ell }^{\prime \prime }}$ via a BSDE formulation. 
As mentioned earlier, this BSDE approach enables us to avoid the estimation of $\mathbf{Z}^{\mathbf{u}, u_{h}^{\prime }, u_{\ell }^{\prime \prime }}$ (Proposition 7.3 of \cite{GLZ1}), and especially remove  the term $\left\Vert \cdot \right\Vert _{\mathcal{H}^{4}\left( \mathbb{R}\right) }^{2}$ with respect to the control process. 
We start with two technical lemmas.
\begin{lemma}
\label{l2}Suppose that Assumptions \emph{(A1)-(A2) }hold and that there
exists $p\geq 2$ such that $\xi _{i}\in L^{p}\left( \Omega ;\mathbb{R}%
\right) $ for all $i\in I_{N}$. For each $\mathbf{u}\in \mathcal{H}%
^{p}\left( \mathbb{R}\right) ^{N}$, the solution $\mathbf{X}^{\mathbf{u}}\in
\mathcal{H}^{p}\left( \mathbb{R}\right) ^{N}$ to (\ref{sde1}) satisfies for
all $i\in I_{N}$,
\begin{equation*}
\sup_{t\in \left[ 0,T\right] }\mathbb{E}\left[ \left\vert X_{t,i}^{\mathbf{u}%
}\right\vert ^{p}\right] \leq C_{X}^{i,p},
\end{equation*}%
with the constant $C_{X}^{i,p}$ defined by
\begin{equation*}
C_{X}^{i,p}:=\left[ I_{\xi _{i},b,\sigma ,p,u_{i}}^{0}+\frac{%
L_{y}^{b}+4\left( p-1\right) \left( L_{y}^{\sigma }\right) ^{2}}{N}%
\sum_{k=1}^{N}I_{\xi _{k},b,\sigma ,p,u_{k}}^{0}\cdot e^{I_{b,\sigma
,p}^{1}\cdot t}\right] \cdot e^{I_{b,\sigma ,p}^{2}\cdot T},
\end{equation*}%
where 
\begin{equation*}
I_{\xi _{i},b,\sigma ,p,u_{i}}^{0}=\mathbb{E}\left[ \left\vert \xi
_{i}\right\vert ^{p}\right] +\left( L^{b}+4\left( p-1\right) \left(
L^{\sigma }\right) ^{2}\right) T+\left( L^{b}+4\left( p-1\right) \left(
L^{\sigma }\right) ^{2}\right) \left\Vert u_{i}\right\Vert _{\mathcal{H}%
^{p}\left( \mathbb{R}\right) }^{p},
\end{equation*}%
\begin{equation*}
I_{b,\sigma ,p}^{1}=\left( 6\left( L^{\sigma }\right) ^{2}+2\left(
L_{y}^{\sigma }\right) ^{2}\right) p^{2}+\left( 3L^{b}+L_{y}^{b}+14\left(
L^{\sigma }\right) ^{2}-2\left( L_{y}^{\sigma }\right) ^{2}\right)
p-2L^{b}+8\left( L^{\sigma }\right) ^{2},
\end{equation*}%
and%
\begin{equation*}
I_{b,\sigma ,p}^{2}=L^{b}\left( 3p-2\right) +L_{y}^{b}\left( p-1\right)
+2\left( p-1\right) \left( 3p-4\right) \left( L^{\sigma }\right)
^{2}+2\left( p-1\right) \left( p-2\right) \left( L_{y}^{\sigma }\right) ^{2}.
\end{equation*}

\end{lemma}

\begin{lemma}
\label{l4}Suppose that Assumptions \emph{(A1)-(A2) }hold and that there
exists $p\geq 2$ such that $\xi _{i}\in L^{p}\left( \Omega ;\mathbb{R}%
\right) $ for all $i\in I_{N}$. For each $\mathbf{u}\in \mathcal{H}%
^{p}\left( \mathbb{R}\right) ^{N}$, $h\in I_{N},$ and $u_{h}^{\prime }\in
\mathcal{H}^{p}\left( \mathbb{R}\right) ,$ the solution $\mathbf{Y}^{\mathbf{%
u},u_{h}^{\prime }}\in \mathcal{H}^{p}\left( \mathbb{R}\right) ^{N}$ to (\ref%
{var1}) satisfies for all $i\in I_{N}$
\begin{eqnarray*}
\sup_{0\leq t\leq T}\mathbb{E}\left[ \left\vert Y_{t,i}^{\mathbf{u}%
,u_{h}^{\prime }}\right\vert ^{p}\right] 
\leq \left[ \left( \frac{L_{y}^{b}+\left( L_{y}^{\sigma }\right)
^{2}3\left( p-1\right) }{N}\right) Te^{\bar{I}_{b,\sigma ,p}^{3}\cdot
T}\left( L^{b}+3\left( p-1\right) L^{\sigma }\right) +\left( 3p-2\right)
\delta _{h,i}\right] e^{\bar{I}_{b,\sigma ,p}^{4}\cdot T}\left\Vert
u_{h}^{\prime }\right\Vert _{\mathcal{H}^{p}\left( \mathbb{R}\right) }^{p},
\end{eqnarray*}
with
\begin{eqnarray*}
& & \bar{I}_{b,\sigma ,p}^{3} = pL^{b}+L_{y}^{b}p+\frac{3}{2}\left( p-1\right)
p(\left( L^{\sigma }\right) ^{2}+\left( L_{y}^{\sigma }\right) ^{2})+\left(
p-1\right) \left( \frac{3}{2}p-2\right), \\
& & \bar{I}_{b,\sigma ,p}^{4} = \bar{I}_{b,\sigma ,p}^{3}-3(p-1)\left(
L_{y}^{\sigma }\right) ^{2}-L_{y}^{b}.
\end{eqnarray*}

\end{lemma}

\begin{theorem}
\label{the1} Given Assumptions \emph{(A1)-(A2)}, then for all $u\in \mathcal{H}^{2}\left( \mathbb{R}%
\right) ^{N}$ and $u_{i}^{\prime },u_{j}^{\prime \prime }\in \mathcal{H}^{2}(\mathbb{R}%
)$,
\begin{equation}
\left\vert \frac{\delta ^{2}V_{i}}{\delta u_{i}\delta u_{j}}\left( \mathbf{u}%
;u_{i}^{\prime },u_{j}^{\prime \prime }\right) -\frac{\delta ^{2}V_{j}}{%
\delta u_{j}\delta u_{i}}\left( \mathbf{u};u_{j}^{\prime \prime
},u_{i}^{\prime }\right) \right\vert \leq \tilde{C}^{i,j}\left\Vert
u_{i}^{\prime }\right\Vert _{\mathcal{H}^{2}\left( \mathbb{R}\right)
}\left\Vert u_{j}^{\prime \prime }\right\Vert _{\mathcal{H}^{2}\left(
\mathbb{R}\right) },  \label{estia}
\end{equation}
with
\begin{equation*}
\tilde{C}^{i,j}=\tilde{C}_{0}^{i,j}+\frac{1}{N}\tilde{C}_{1}^{i,j}+\frac{1}{N^{2}}%
\tilde{C}_{2}^{i,j},
\end{equation*}
where $L_{y}^{b,\sigma}=L_{y}^{b}+3\left(
L_{y}^{\sigma }\right) ^{2},$
\begin{eqnarray*}
\tilde{C}_{0}^{i,j} & \!\!\!=\!\!\! & \mathbb{E}\Big [\left\Vert \partial _{x_{i}x_{j}}^{2}\Delta^f_{i,j}\right\Vert _{L^{\infty }}+\left\Vert \partial
_{x_{i}u_{j}}^{2}\Delta^f_{i,j}\right\Vert _{L^{\infty }}+\left\Vert
\partial _{u_{i}x_{j}}^{2}\Delta^f_{i,j}\right\Vert _{L^{\infty
}}+\left\Vert \partial _{u_{i}u_{j}}^{2}\Delta^f_{i,j}\right\Vert
_{L^{\infty }}+\left\Vert \partial _{x_{i}x_{j}}^{2}\Delta^g_{i,j}\right\Vert _{L^{\infty }}\Big], \\
\tilde{C}_{1}^{i,j} & \!\!\!=\!\!\! & \mathbb{E}\Bigg \{L_{y}^{b,\sigma }\Bigg (\sum_{\ell \in
I_{N}\backslash \left\{ j\right\} }\left( \left\Vert \partial _{x_{i}x_{\ell
}}^{2}\Delta^f_{i,j}\right\Vert _{L^{\infty }}+\left\Vert \partial
_{u_{i}x_{\ell }}^{2}\Delta^f_{i,j}\right\Vert _{L^{\infty }}+\left\Vert
\partial _{x_{i}x_{\ell }}^{2}\Delta^g_{i,j}\right\Vert _{L^{\infty
}}\right)  \\
&&+\sum_{h\in I_{N}\backslash \left\{ i\right\} }\left( \left\Vert \partial
_{x_{h}x_{j}}^{2}\Delta^f_{i,j}\right\Vert _{L^{\infty }}+\left\Vert
\partial _{x_{h}u_{j}}^{2}\Delta^f_{i,j}\right\Vert _{L^{\infty
}}+\left\Vert \partial _{x_{h}x_{j}}^{2}\Delta^g_{i,j}\right\Vert
_{L^{\infty }}\right) \Bigg ) \\
&&+C\sqrt{\Lambda _{1}}L_{y}^{\sigma }\left( 1+L_{y}^{b,\sigma }+\left(
L_{y}^{b,\sigma }\right) ^{2}\right) +C\sqrt{\Lambda _{1}}\left( L^{\sigma
}\left( L_{y}^{b,\sigma }\right) ^{2}+2L_{y}^{\sigma }\right) \left(
1+L_{y}^{b,\sigma }+\left( L_{y}^{b,\sigma }\right) ^{2}\right)  \\
&&+2C\sqrt{\Lambda _{1}}L_{y}^{b,\sigma }\left( L^{\sigma }+L_{y}^{\sigma
}\right) +C\sqrt{\Lambda _{1}}L_{y}^{b}\left( 1+L_{y}^{b,\sigma }+\left(
L_{y}^{b,\sigma }\right) ^{2}\right)  \\
&& +C\sqrt{\Lambda _{1}}\left( L^{b}\left( L_{y}^{b,\sigma }\right)
^{2}+2L_{y}^{b}\right) \left( 1+L_{y}^{b,\sigma }+\left( L_{y}^{b,\sigma
}\right) ^{2}\right) +2C\sqrt{\Lambda _{1}}L_{y}^{b,\sigma }\left(
L^{b}+L_{y}^{b}\right) \Bigg \}, \\
\tilde{C}_{2}^{i,j} & \!\!\!=\!\!\! & \mathbb{E}\Bigg \{L_{y}^{b,\sigma }\Bigg [L_{y}^{b,\sigma
}\sum_{\substack{ \ell \in I_{N}\backslash \left\{ j\right\} , \\ h\in
I_{N}\backslash \left\{ i\right\} }}\left\Vert \partial _{x_{h}x_{\ell
}}^{2}\Delta^f_{i,j}\right\Vert _{L^{\infty }}+\sum_{\substack{ \ell \in
I_{N}\backslash \left\{ j\right\} , \\ h\in I_{N}\backslash \left\{
i\right\} }}\left\Vert \partial _{x_{h}x_{\ell }}^{2}\Delta^g_{i,j}\right\Vert _{L^{\infty }} \\
&&+C\sqrt{\Lambda _{1}}L_{y}^{\sigma }L_{y}^{b,\sigma }+C\sqrt{\Lambda _{1}}
L_{y}^{b}L_{y}^{b,\sigma }\Bigg ]\Bigg\},
\end{eqnarray*}
in which
\begin{eqnarray}
\Lambda_{1} & \!\!\!=\!\!\! & C_{1}\mathbb{E}\Bigg [\sum_{\ell \in I_{N}}\left\vert \left(
\partial _{x_{\ell }}\Delta _{i,j}^{g}\right) \left( 0\right) \right\vert
^{2}+\sum_{\ell ,k\in I_{N}}\left\Vert \partial _{x_{\ell }x_{k}}\Delta
_{i,j}^{g}\right\Vert _{L^{\infty }}^{2}  \notag \\
& & + 3\int_{0}^{T}\Bigg (\sum_{\ell }\left\vert \left( \partial _{x_{\ell
}}\Delta _{i,j}^{f}\right) \left( t,0,0\right) \right\vert ^{2}+\sum_{\ell
,k\in I_{N}}\left\Vert \partial _{x_{\ell }x_{k}}\Delta
_{i,j}^{f}\right\Vert _{L^{\infty }}^{2}  
+\sum_{\ell ,k\in I_{N}}\left\Vert \partial _{x_{\ell }u_{k}}\Delta
_{i,j}^{f}\right\Vert _{L^{\infty }}^{2}\Bigg )\mathrm{d}t\Bigg],
\label{lamda1}
\end{eqnarray}
where $C_{1}$ depends on the time $T$ and $\max \left\{ \left\Vert \mathbf{B}_{0,x,y}
\left( t,\mathbf{X}_{t}^{\mathbf{u}},\mathbf{u}_{t}\right) \right\Vert
_{2},\left\Vert \mathbf{\Pi }^{1}_{u}\left( t,\mathbf{X}_{t}^{\mathbf{u}},%
\mathbf{u}_{t}\right) \right\Vert _{2},\ldots ,\left\Vert \mathbf{\Pi }_{u}%
^{N}\left( t,\mathbf{X}_{t}^{\mathbf{u}},\mathbf{u}_{t}\right) \right\Vert
_{2}\right\} $. Consequently, if $\sup_{i\in I_{N},u_{i}\in \mathcal{A}_{i}}\left\Vert
u_{i}\right\Vert _{\mathcal{H}^{2}\left( \mathbb{R}\right) }<\infty $ and $%
0\in \mathcal{A}_{i},$ then $\mathcal{G}^{\mathsf{op}}$ is an $\alpha $%
-potential game with
\begin{equation}
\alpha \leq C\max_{i\in I_{N}}\sum_{j\in I_{N}\backslash \left\{ i\right\} }%
\tilde{C}^{i,j},  \label{alpha}
\end{equation}
where $C$ is a constant depending only on
the state coefficients and time horizon as specified in (\ref{potenfv}).
\end{theorem}


\begin{corollary}
Taking the same setting as in Theorem \ref{the1}, in the special case of $b_{i}\left( t,X,\mathbf{X},\mathbf{u}\right) =\bar{b}_{i}\left( t,X,%
\mathbf{X}\right) +u,$ $\sigma _{i}\left( t,X,\mathbf{X},\mathbf{u}\right)
=\sigma _{i}\left( t\right), $ as considered in \cite{GLZ1},  $%
L^{\sigma }=L_{y}^{\sigma }=0$, and the estimation of $\alpha$ in (\ref{alpha}) can be simplified such that
\begin{eqnarray*}
\tilde{C}^{i,j} 
&=&\mathbb{E}\Bigg \{\left\Vert \partial
_{x_{i}x_{j}}^{2}\Delta^f_{i,j}\right\Vert _{L^{\infty }}+\left\Vert
\partial _{x_{i}u_{j}}^{2}\Delta^f_{i,j}\right\Vert _{L^{\infty
}}+\left\Vert \partial _{u_{i}x_{j}}^{2}\Delta^f_{i,j}\right\Vert
_{L^{\infty }} +\left\Vert \partial _{u_{i}u_{j}}^{2}\Delta^f_{i,j}\right\Vert
_{L^{\infty }}+\left\Vert \partial _{x_{i}x_{j}}^{2}\Delta^g_{i,j}\right\Vert _{L^{\infty }} \\
&&+\frac{L_{y}^{b}}{N}\Bigg [\sum_{\ell \in I_{N}\backslash \left\{
j\right\} }\left( \left\Vert \partial _{x_{i}x_{\ell }}^{2}\Delta^f_{i,j}\right\Vert _{L^{\infty }}+\left\Vert \partial _{u_{i}x_{\ell
}}^{2}\Delta^f_{i,j}\right\Vert _{L^{\infty }}+\left\Vert \partial
_{x_{i}x_{\ell }}^{2}\Delta^g_{i,j}\right\Vert _{L^{\infty }}\right)  \\
&&+\sum_{h\in I_{N}\backslash \left\{ i\right\} }\left( \left\Vert \partial
_{x_{h}x_{j}}^{2}\Delta^f_{i,j}\right\Vert _{L^{\infty }}+\left\Vert
\partial _{x_{h}u_{j}}^{2}\Delta^f_{i,j}\right\Vert _{L^{\infty }}+\left\Vert \partial
_{x_{h}x_{j}}^{2}\Delta^g_{i,j}\right\Vert _{L^{\infty }}\right)  \Bigg ] \\
&&+\frac{\left( L_{y}^{b}\right) ^{2}}{N^{2}}\Bigg (\sum_{\ell \in
I_{N}\backslash \left\{ j\right\} ,h\in I_{N}\backslash \left\{ i\right\}
}\left\Vert \partial _{x_{h}x_{\ell }}^{2}\Delta^f_{i,j}\right\Vert
_{L^{\infty }}+\sum_{\ell \in I_{N}\backslash \left\{ j\right\} ,h\in
I_{N}\backslash \left\{ i\right\} }\left\Vert \partial _{x_{h}x_{\ell
}}^{2}\Delta^g_{i,j}\right\Vert _{L^{\infty }}\Bigg )\Bigg\} \\
&&+CL_{y}^{b}\sqrt{\Lambda _{1}}\Bigg \{\frac{1}{N}\bigg [\left(1+L_{y}^{b}+\left( L_{y}^{b}\right) ^{2}\right) +\left(L^{b}L_{y}^{b}+2\right) \left( 1+L_{y}^{b}+\left( L_{y}^{b}\right)^{2}\right)  +2\left( L^{b}+L_{y}^{b}\right) \bigg ]+\frac{\left( L_{y}^{b}\right)^{2}}{N^{2}}\Bigg\},
\end{eqnarray*}
where $\Lambda _{1}$ is defined in (\ref{lamda1}).
\end{corollary}




\begin{corollary}
\label{cor2} Given Assumptions \emph{(A1)-(A2)}, suppose that there exist $L,%
\tilde{L}\geq 0$, and $\beta >1/2$ such that
\begin{equation*}
\max \left\{ \sup_{i\in I_{N},u_{i}\in \mathcal{A}_{i}}\left\Vert
u_{i}\right\Vert _{\mathcal{H}^{2}\left( \mathbb{R}\right) },\max_{i\in
I_{N}}\mathbb{E}\left[ \left\vert \xi _{i}\right\vert ^{2}\right]
,\max_{i\in I_{N}}\left( L^{b_{i}}+L_{y}^{b_{i}}+L^{\sigma
_{i}}+L_{y}^{\sigma _{i}}\right) \right\} \leq L; 
\end{equation*}
for all $i,$ $j\in I_{N}$, 
satisfy for all $\left( t,x,u\right) \in \left[ 0,T%
\right] \times \mathbb{R}^{N}\times \mathbb{R}^{N},$ $\left\vert \left(
\partial _{u_{i}u_{j}}^{2}\Delta _{i,j}^{f}\right) \left( t,x,u\right)
\right\vert \leq \tilde{L}N^{-2\beta }$; and for all $h,$ $\ell \in I_{N}$, 
\begin{eqnarray*}
\begin{aligned}
\left\vert \left( \partial _{x_{h}}\Delta _{i,j}^{f}\right) \left(
t,0,0\right) \right\vert  & \leq \tilde{L}N^{-\beta }, \\
\text{ }\left\vert \left( \partial _{x_{h}x_{\ell }}^{2}\Delta
_{i,j}^{f}\right) \left( t,x,u\right) \right\vert +\left\vert \left(
\partial _{x_{h}u_{\ell }}^{2}\Delta _{i,j}^{f}\right) \left( t,x,u\right)
\right\vert  &\leq \tilde{L}\left(\delta_{h,\ell}N^{-\beta }+(1-\delta_{h,\ell})N^{-2\beta }\right), \\
\left\vert \left( \partial _{x_{h}}\Delta _{i,j}^{g}\right) \left( 0\right)
\right\vert  &\leq \tilde{L}N^{-\beta }, \\
\text{ }\left\vert \left( \partial _{x_{h}x_{\ell }}^{2}\Delta
_{i,j}^{g}\right) \left( x\right) \right\vert  &\leq \tilde{L}\left(\delta_{h,\ell}N^{-\beta }+(1-\delta_{h,\ell})N^{-2\beta}\right). 
\end{aligned} 
\end{eqnarray*}
Then $\mathcal{G}^{\mathsf{op}}$ is an $\alpha $-potential game 
\begin{eqnarray*} 
\begin{aligned}
\alpha \leq & \frac{\tilde{L}}{N^{2\beta }}\left( C+2L_{y}^{b,\sigma }+2\left(
L_{y}^{b,\sigma }\right) ^{2}\right) +\frac{4L_{y}^{b,\sigma }\tilde{L}}{N}%
\cdot \frac{1}{N^{\min \left\{ \beta ,2\beta -1\right\} }} \\
& +C\sqrt{C_{1}}\max \left\{ C^{1,b,\sigma },C^{2,b,\sigma }\right\} \frac{1%
}{N^{\frac{\beta +1}{2}}}, 
\end{aligned}
\end{eqnarray*}
where $C^{1,b,\sigma}$ and $C^{2,b,\sigma }$ are defined in (\ref{mc1})-(\ref{mc2}), and $C\geq 0$ is a constant independent of $N$ and $\beta $.
\end{corollary}

It is worth noting that if the cost function in (\ref{cost}) depends on the empirical measure of the joint state and control profiles, then Corollary  \ref{cor2} holds with $\beta = 1$ and a constant $\tilde{L}$ independent of $N$. Consequently, the $N$-player game becomes an $\alpha$-potential game with $\alpha = \mathcal{O}(1/N)$. 

For an $N$-player games with mean-field type interactions where the partial derivatives of the cost functions exhibit suitable decay as the number of players increases, we demonstrate that such a game constitutes an $\alpha$-potential game, with the associated constant decaying to zero as the number of players $N \to \infty$. More specifically, consider the following examples. 

\begin{example}
Suppose that there exists $L>0$ such that $\sup_{i\in I_{N},u_{i}\in
\mathcal{A}_{i}}\left\Vert u_{i}\right\Vert _{\mathcal{H}^{2}\left( \mathbb{R%
}\right) }<L$. Consider
\begin{eqnarray*} 
\begin{aligned}
b_{i}\left( t,x_i,x,u\right) &=\tilde{b}_{i}\left( t,x_{i},\frac{1}{N}%
\sum_{\ell =1}^{N}\delta _{x_{\ell }},u\right) ,\text{ }\sigma _{i}\left(
t,x_{i},x,u\right) =\tilde{\sigma}_{i}\left( t,x_{i},\frac{1}{N}\sum_{\ell
=1}^{N}\delta _{x_{\ell}},u\right), \\
f_{i}\left( t,x,u\right) &=f_{0}\left( t,x,u\right) +c_{i}\left(
u_{i}\right) +\tilde{f}_{i}\left( t,\frac{1}{N}\sum_{\ell =1}^{N}\delta
_{x_{\ell }}\right) , \\
g_{i}\left( x\right) &=g_{0}\left( x\right) +\tilde{g}_{i}\left( \frac{1}{N}%
\sum_{\ell =1}^{N}\delta _{x_{\ell }}\right), 
\end{aligned} 
\end{eqnarray*}
where $\tilde{b}_{i}:\left[ 0,T\right] \times \mathbb{R}\times \mathcal{P}_{2}\left( \mathbb{R}\right)\times \mathbb{R\rightarrow R},$ $\tilde{\sigma}_{i}:\left[ 0,T\right]
\times \mathbb{R}\times \mathcal{P}_{2}\left( \mathbb{R}\right)\times \mathbb{R\rightarrow R},$ $%
f_{0}:\left[ 0,T\right] \times \mathbb{R}^{N}\times \mathbb{R\rightarrow R},$
$c_{i}:\mathbb{R\rightarrow R},$ $\tilde{f}_{i}:\left[ 0,T\right] \times
\mathcal{P}_{2}\left( \mathbb{R}\right) \rightarrow \mathbb{R}$, $g_{0}:%
\mathbb{R}^{N}\mathbb{\rightarrow R}$ and $\tilde{g}_{i}:\mathcal{P}%
_{2}\left( \mathbb{R}\right) \rightarrow \mathbb{R}.$ According to Proposition 5.35 and 5.91 in \cite{cd18}, we have
\begin{equation*}
\left\vert \partial_{x_{h}}\tilde{b}_{i}\left(t,x_i,\mu,u\right) \right\vert
+\left\vert \partial _{x_{h}}\tilde{\sigma}_{i}\left(t,x_i,\mu,u\right)
\right\vert \leq \frac{C}{N},
\end{equation*}
\begin{equation*}
\left\vert \partial _{x_{h}}\Delta _{i,j}^{f}\left( t,0,0\right) \right\vert
+\left\vert \partial _{x_{h}}\Delta _{i,j}^{g}\left( 0\right) \right\vert
\leq \frac{C}{N},
\end{equation*}
and
\begin{equation*}
\left\vert \partial _{x_{\ell }x_{h}}^{2}\Delta _{i,j}^{f}\left(
t,x,u\right) \right\vert +\left\vert \partial _{x_{\ell }x_{h}}^{2}\Delta
_{i,j}^{g}\left(x\right) \right\vert \leq C\left( \frac{1}{N}\delta_{h,\ell}+\frac{1}{N^{2}}(1-\delta_{h,\ell})\right),
\end{equation*}
for some positive constant $C$ independent of $N$. Therefore this game is an $\alpha $-potential game with $\alpha \leq C/N$, for a positive
constant $C$ independent of $N.$

\end{example}

\begin{example}
\label{exa1} Consider a system with $N$ weakly-coupled
negligible agents, where the state $X_{t,i}$ satisfies the following linear controlled  stochastic system:
\begin{equation}
\left\{\begin{aligned}
\mathrm{d}X_{t,i} = & \left( A_{t,i}X_{t,i}+\bar{A}_{t,i}X_{t}^{\left(
N\right) }+B_{t,i}u_{t,i}+b_{t,i}\right) \mathrm{d}t \\
& +\left( C_{t,i}X_{t,i}+\bar{C}_{t,i}X_{t}^{\left( N\right)}+D_{t,i}u_{t,i}+\sigma _{t,i}\right) \mathrm{d}W_{t}^{i}, \\
X_{0,i} = & \, \xi_{i},
\end{aligned}\right.   \label{lar}
\end{equation}
where $X_{t}^{(N)}=\frac{1}{N}\sum_{i=1}^{N}X_{t,i}$, $A,\bar{A},b,B,C,\bar{C},D$, and $\sigma$ are deterministic matrix-valued
functions with appropriate dimensions, and the cost functional
is defined as
\begin{equation}
V_{i}\left( \mathbf{u}\right) =\frac{1}{2}\mathbb{E}\left[
\int_{0}^{T}\left( \hat{Q}_{t,i}\left( X_{t,i}-X_{t}^{(N)}\right)
^{2}+R_{t,i}u_{t,i}^{2}\right) \mathrm{d}t+G_{i}\left(
X_{T,i}-X_{T}^{(N)}\right) ^{2}\right].  \label{larcost}
\end{equation}
Here for all $t\in \left[ 0,T\right] ,$ $0\leq \hat{Q}_{i,t}\in \mathbb{S}%
^{n},$ $0<R_{i,t}\in \mathbb{S}^{n},$ $0\leq G_{i}\in \mathbb{S}^{n}.$
Suppose that all the coefficients are bounded by a positive constant $C.$
Noting that $L^{b}=C,$ $L_{y}^{b}=\frac{C}{N},$ $L^{\sigma }=C,$ $%
L_{y}^{\sigma }=\frac{C}{N},$ by Theorem \ref{the1}, we get, after some
simple calculations, 
\begin{eqnarray*}
\begin{aligned}
\tilde{C}^{i,j} = &\,C\Bigg \{\left\Vert \partial _{x_{i}x_{j}}^{2}\Delta
^f_{i,j}\right\Vert _{L^{\infty }}+\frac{1}{N^{2}}\Bigg [\sum_{\ell \in
I_{N}\backslash \left\{ j\right\} }\left\Vert \partial _{x_{i}x_{\ell
}}^{2}\Delta^f_{i,j}\right\Vert _{L^{\infty }}+\sum_{h\in I_{N}\backslash
\left\{ i\right\} }\left\Vert \partial _{x_{h}x_{j}}^{2}\Delta^f_{i,j}\right\Vert _{L^{\infty }}\Bigg ] \\
&+\frac{1}{N^{3}}\Bigg (\sum_{\ell \in I_{N}\backslash \left\{ j\right\}
,h\in I_{N}\backslash \left\{ i\right\} }\left\Vert \partial _{x_{h}x_{\ell
}}^{2}\Delta^f_{i,j}\right\Vert _{L^{\infty }}\Bigg )+\sqrt{\Lambda _{1}}%
\Bigg (\frac{1}{N^{2}}+\frac{1}{N^{5}}\Bigg )\Bigg \} \\
=&\,C\Bigg [\left( \hat{Q}_{t,i}-\hat{Q}_{t,j}\right) \frac{1}{N}\left( 1-%
\frac{1}{N}\right) -\frac{1}{N^{2}}\left( N-1\right) \left( \hat{Q}_{t,i}%
\frac{1}{N}\left( 1-\frac{1}{N}\right) +\hat{Q}_{t,j}\frac{1}{N^{2}}\right)
\\
&+\frac{1}{N^{2}}\left( N-1\right) \left( \hat{Q}_{t,i}\frac{1}{N^{2}}+\hat{%
Q}_{t,j}\frac{1}{N}\left( 1-\frac{1}{N}\right) \right) +\frac{1}{N^{3}}%
\left( N^{2}-2\right) \left( \hat{Q}_{t,i}-\hat{Q}_{t,j}\right) \frac{1}{%
N^{2}} \\
&+\sqrt{\Lambda _{1}}\frac{1}{N^{2}}+\left( G_{i}-G_{j}\right) \frac{1}{N}%
\left( 1-\frac{1}{N}\right) -\frac{1}{N^{2}}\left( N-1\right) \left( G_{i}%
\frac{1}{N}\left( 1-\frac{1}{N}\right) +G_{j}\frac{1}{N^{2}}\right)  \\
&+\frac{1}{N^{2}}\left( N-1\right) \left( G_{i}\frac{1}{N^{2}}+G_{j}\frac{1%
}{N}\left( 1-\frac{1}{N}\right) \right) +\frac{1}{N^{3}}\left(
N^{2}-2\right) \left( G_{i}-G_{j}\right) \frac{1}{N^{2}}+\sqrt{\Lambda _{1}}%
\frac{1}{N^{2}}\Bigg ] \\
\leq &\,C\Bigg [\frac{1}{N}\left( \left\vert \hat{Q}_{t,i}-\hat{Q}%
_{t,j}\right\vert +\left\vert G_{i}-G_{j}\right\vert \right) +\sqrt{\Lambda
_{1}}\frac{1}{N^{2}}\Bigg], 
\end{aligned}
\end{eqnarray*}
where $\Lambda _{1}$ is given in (\ref{lamda1}), and $\alpha$ given in (\ref{alpha}).

In the special case of $G_{i}=G_{j}$ and $\hat{Q}_{t,i}=\hat{Q}_{t,j}$, even when the
state dynamics are heterogenous, this is a potential game with $\alpha =0$.
For general heterogenous cases, as $N\rightarrow \infty $, $\alpha
\rightarrow 0$ as long as the rest of the terms are bounded.
\end{example}

\begin{example}  [Games with common noise]

Now we analyze a class of stochastic differential games with common noise (cf. \cite{cd18}). 
 Take a finite time horizon $T > 0$, consider an $N$-player
stochastic differential game, where the private state process $X_i$ of
player $i$ is determined by the solution of the following SDE: 

\begin{equation}
\left\{
\begin{array}{rcl}
\mathrm{d}X_{t,i} & = & \mathsf{b}_{t,i}u_{t,i}\mathrm{d}t+\sigma _{t,i}%
\mathrm{d}W_{t,i}+\mathrm{d}W_{t}^{0}, \\
X_{0,i} & = & \xi _{i},%
\end{array}%
\right.   \label{comnoise}
\end{equation}%
Here, $W^{0}$, the \emph{common noise}, is a Brownian motion defined on $\left( \Omega ^{0},\mathcal{F}^{0},\mathbb{P}^{0}\right) $, and $W_{t,1},\ldots ,W_{t,N}$ defined on $\left( \Omega ,\mathcal{F},\mathbb{P}%
\right) $ are independent Brownian motions independent of $W^{0}$. 
We denote by $\mathbb{F}=(\mathbb{F}_{t})_{t\geq 0}$ the $\mathbb{P}$-completion of the natural
filtration generated by $W$ and $ W^{0}$, and without loss of generality we assume that $\mathcal{F}%
=\mathcal{F}_{\infty }$. 
For each $i\in I_{N}$, assume that $u_{i}\in \mathcal{A}_{i}$, where $\mathcal{A}_{i}$
denotes a convex subset of $\mathcal{H}^{2}(\mathbb{R})$ representing
the set of admissible open-loop controls of player $i$.  
Here, $L^{b}=0,$ $L_{y}^{b}=0,$ $L^{\sigma }=0,$ $L_{y}^{\sigma }=0$.

\begin{corollary}
    
 \label{comthe}
    Assume that $\mathsf{b}_{t,i}$ and $\sigma_{t,i}$ are bounded $\mathcal{F}_{t}$-adapted processes, and  the cost functional
is defined as in (\ref{larcost}). Then this game is an $\alpha $-potential game with $\alpha$ given in (\ref{alpha}) and 
\begin{equation*}
\tilde{C}^{i,j}\leq C\Big [\frac{1}{N}\left( \left\vert
\hat{Q}_{i,t}-\hat{Q}_{j,t}\right\vert +\left\vert G_{i}-G_{j}\right\vert \right) +\sqrt{%
\Lambda _{1}}\frac{1}{N^{2}}\Big ],
\end{equation*}
where $\Lambda _{1}$ is given in (\ref{lamda1}) after adjustment of
coefficients.
\end{corollary}
\end{example}

The proof can be seen in Appendix \ref{APP}.

It is worth noting that the framework studied in \cite{GLZ1} cannot cover this example due to the presence of common noise. In fact, by considering the limiting problem (denoted by $X$) and assuming that each player adopts the same strategy, one can
represent the player distribution by a conditional law of a single representative player given a common noise. Therefore, the term $X_{t}^{\left( N\right) }$ in (\ref{larcost}) converges to $\mathbb{E}\left[ \left. X_{t}\right\vert W^{0}\right] $, which denotes the conditional expectation of $X_{t}$ given the $\sigma $-algebra $\mathcal{F}^{0}$ generated by $W^{0}$ augmented by $\mathbb{P}^{0}$-null sets.

If the control process is adapted to the filtration generated by $W^{0}$, it can be addressed through the dynamic programming principle (see \cite{PW2017}).
For the general case of $\alpha$-potential stochastic differential games with common noise, it can be analyzed via similar technique, with details left for 
future investigation.


\section{Proof of Main Results}\label{sect6} 

\paragraph{Proof of Theorem \ref{Thm:Y}.}

By Lemma \ref{bs1}, there exists a unique adapted solution $\left(P_{i},\left\{ Q_{j}\right\} _{j=1,\ldots ,N}\right) $ to BSDE (\ref{adj1}).
By It\^{o}'s formula to $\left\langle\mathbf{Y}^{\mathbf{u},u_{h}^{\prime}},P_{i}\right\rangle $ on $\left[ 0,T\right] $, we obtain 

\begin{eqnarray*}
\begin{aligned}
\mathbb{E}\left[ \left\langle \mathbf{Y}_{T}^{\mathbf{u},u_{h}^{\prime
}},P_{T,i}\right\rangle \right]
=& \, \mathbb{E}\left[ \left\langle \mathbf{Y}_{T}^{\mathbf{u},u_{h}^{\prime
}},\left( \partial _{x}g_{i}\right) \left( \mathbf{X}_{T}^{\mathbf{u}%
}\right) \right\rangle \right] \\
=& \, \mathbb{E}\left[ \int_{0}^{T}\left\langle P_{t,i},\mathbf{B}%
_{0,x,y}\left( t,\mathbf{X}_{t}^{\mathbf{u}},\mathbf{u}_{t}\right) \mathbf{Y}%
_{t}^{\mathbf{u},u_{h}^{\prime }}+\mathbf{B}_{1,u}^{h}\left( t,\mathbf{X}%
_{t}^{\mathbf{u}},\mathbf{u}_{t}\right) u_{t,h}^{\prime }\right\rangle
\mathrm{d}t\right] \\
&+\mathbb{E}\Bigg [\int_{0}^{T}\Big (-\left\langle \mathbf{Y}_{t}^{\mathbf{u%
},u_{h}^{\prime }},\mathbf{B}_{0,x,y}\left( t,\mathbf{X}_{t}^{\mathbf{u}},%
\mathbf{u}_{t}\right) ^{\top }P_{t,i}\right\rangle -\left\langle \mathbf{Y}%
_{t}^{\mathbf{u},u_{h}^{\prime }},\sum_{j=1}^{N}\mathbf{\Pi }%
_{0,x,y}^{j}\left( t,\mathbf{X}_{t}^{\mathbf{u}},\mathbf{u}_{t}\right)
^{\top }Q_{t,j}\right\rangle \\
&-\left\langle \mathbf{Y}_{t}^{\mathbf{u},u_{h}^{\prime }},\left( \partial
_{x}f_{i}\right) \left( t,\mathbf{X}_{t}^{\mathbf{u}},\mathbf{u}_{t}\right)
\right\rangle \Big )\mathrm{d}t\Bigg ] \\
&+\mathbb{E}\Bigg \{\int_{0}^{T}\Big [\sum_{j=1}^{N}\left( \mathbf{\Pi }%
_{0,x,y}^{j}\left( t,\mathbf{X}_{t}^{\mathbf{u}},\mathbf{u}_{t}\right)
\mathbf{Y}_{t}^{\mathbf{u},u_{h}^{\prime }}+\mathbf{\Pi }_{1,u}^{j}\left( t,%
\mathbf{X}_{t}^{\mathbf{u}},\mathbf{u}_{t}\right) \delta
_{i,h}u_{t,h}^{\prime }\right) Q_{t,j}\Big ]\mathrm{d}t\Bigg \} \\
=& \, \mathbb{E}\Bigg [\int_{0}^{T}\Big (\left\langle P_{t,i},\mathbf{B}%
_{1,u}^{h}\left( t,\mathbf{X}_{t}^{\mathbf{u}},\mathbf{u}_{t}\right)
u_{t,h}^{\prime }\right\rangle -\left\langle \mathbf{Y}_{t}^{\mathbf{u}%
,u_{h}^{\prime }},\left( \partial _{x}f_{i}\right) \left( t,\mathbf{X}_{t}^{%
\mathbf{u}},\mathbf{u}_{t}\right) \right\rangle \\
&+\left\langle \mathbf{\Pi }_{1,u}^{h}\left( t,\mathbf{X}_{t}^{\mathbf{u}},%
\mathbf{u}_{t}\right) ,Q_{t,h}u_{t,h}^{\prime }\right\rangle \Big )\mathrm{d}%
t\Bigg],
\end{aligned}
\end{eqnarray*}
from which we get the desired result. \hfill $\Box$

\paragraph{Proof of Proposition \ref{V2}.}

Recall (3.8) in \cite{GLZ1}, we have 
\begin{eqnarray}
\frac{\delta ^{2}V_{i}}{\delta u_{h}\delta u_{\ell }}\left( \mathbf{u}%
;u_{h}^{\prime },u_{\ell }^{\prime \prime }\right)  &=&\mathbb{E}\Bigg [%
\int_{0}^{T}\Bigg (\left(
\begin{array}{c}
\mathbf{Y}_{t}^{\mathbf{u},u_{h}^{\prime }} \\
u_{t,h}^{\prime }%
\end{array}%
\right) ^{\top }\left(
\begin{array}{cc}
\partial _{xx}^{2}f_{i} & \partial _{xu_{\ell }}^{2}f_{i} \\
\partial _{u_{h}x}^{2}f_{i} & \partial _{u_{h}u_{\ell }}^{2}f_{i}%
\end{array}%
\right) \left(
\begin{array}{c}
\mathbf{Y}_{t}^{\mathbf{u},u_{\ell }^{\prime \prime }} \\
u_{t,\ell }^{\prime \prime }%
\end{array}%
\right)   \notag \\
&&+\left( \partial _{x}f_{i}\right) ^{\top }\left( t,\mathbf{X}_{t}^{\mathbf{%
u}},\mathbf{u}_{t}\right) \mathbf{Z}_{t}^{\mathbf{u},u_{h}^{\prime },u_{\ell
}^{\prime \prime }}\Bigg )\mathrm{d}t\Bigg ]  \notag \\
&&+\mathbb{E}\left[ \left( \mathbf{Y}_{T}^{\mathbf{u},u_{h}^{\prime
}}\right) ^{\top }\partial _{xx}^{2}g_{i}\left( \mathbf{X}_{T}^{\mathbf{u}%
}\right) \mathbf{Y}_{T}^{\mathbf{u},u_{\ell }^{\prime \prime }}+\left(
\partial _{x}g_{i}\right) ^{\top }\left( \mathbf{X}_{T}^{\mathbf{u}}\right)
\mathbf{Z}_{T}^{\mathbf{u},u_{h}^{\prime },u_{\ell }^{\prime \prime }}\right]
\notag \\
&=&\mathbb{E}\Bigg [\int_{0}^{T}\Bigg (\text{tr}\left[ \partial
_{xx}^{2}f_{i}\cdot \mathbf{Y}_{t}^{\mathbf{u},u_{\ell }^{\prime \prime
}}\cdot \left( \mathbf{Y}_{t}^{\mathbf{u},u_{h}^{\prime }}\right) ^{\top }%
\right] +\partial _{u_{h}x}^{2}f_{i}\cdot \mathbf{Y}_{t}^{\mathbf{u},u_{\ell
}^{\prime \prime }}\cdot u_{t,h}^{\prime }  \notag \\
&&+\left( \mathbf{Y}_{t}^{\mathbf{u},u_{h}^{\prime }}\right) ^{\top }\cdot
\partial _{xu_{\ell }}^{2}f_{i}\cdot u_{t,\ell }^{\prime \prime }+\partial
_{u_{h}u_{\ell }}^{2}f_{i}\cdot u_{t,h}^{\prime }\cdot u_{t,\ell }^{\prime
\prime }  \notag \\
&&+\left( \partial _{x}f_{i}\right) ^{\top }\left( t,\mathbf{X}_{t}^{\mathbf{%
u}},\mathbf{u}_{t}\right) \mathbf{Z}_{t}^{\mathbf{u},u_{h}^{\prime },u_{\ell
}^{\prime \prime }}\Bigg )\mathrm{d}t\Bigg ]  \notag \\
&&+\mathbb{E}\left[ \left( \mathbf{Y}_{T}^{\mathbf{u},u_{h}^{\prime
}}\right) ^{\top }\partial _{xx}^{2}g_{i}\left( \mathbf{X}_{T}^{\mathbf{u}%
}\right) \mathbf{Y}_{T}^{\mathbf{u},u_{\ell }^{\prime \prime }}+\left(
\partial _{x}g_{i}\right) ^{\top }\left( \mathbf{X}_{T}^{\mathbf{u}}\right)
\mathbf{Z}_{T}^{\mathbf{u},u_{h}^{\prime },u_{\ell }^{\prime \prime }}\right]
.  \label{aerfa}
\end{eqnarray}
Applying It\^{o}'s formula to $\left\langle \mathbf{Z}_{t}^{\mathbf{u}%
,u_{h}^{\prime },u_{\ell }^{\prime \prime }},P_{t,i}\right\rangle $ on $%
[0,T] $, we have
\begin{eqnarray*}
\begin{aligned}
\mathbb{E}\left[ \left\langle \mathbf{Z}_{T}^{\mathbf{u},u_{h}^{\prime
},u_{\ell }^{\prime \prime }},\partial _{x}g_{i}\left( \mathbf{X}_{T}^{%
\mathbf{u}}\right) \right\rangle \right] 
=& \, \mathbb{E}\Bigg [-\int_{0}^{T}\Big <\mathbf{Z}_{t}^{\mathbf{u}%
,u_{h}^{\prime },u_{\ell }^{\prime \prime }},\left( \partial
_{x}f_{i}\right) \left( t,\mathbf{X}_{t}^{\mathbf{u}},\mathbf{u}_{t}\right) %
\Big >\mathrm{d}t \\
&+\int_{0}^{T}\left\langle P_{t,i},\mathfrak{F}_{t}^{\mathbf{u}%
,u_{h}^{\prime },u_{\ell }^{\prime \prime }}+\partial _{yy}^{2}\bar{b}\left(
t,\mathbf{X}_{t},\mathbf{u}_{t}\right) \star \mathcal{Y}_{t}^{\mathbf{u}%
,u_{h}^{\prime },u_{\ell }^{\prime \prime }}+\Gamma _{t}^{\mathbf{u}%
,u_{h}^{\prime },u_{\ell }^{\prime \prime }}\left( \bar{b}\right)
\right\rangle \mathrm{d}t \\
&+\int_{0}^{T}\sum_{j=1}^{N}\left\langle Q_{t,j},\mathfrak{G}_{t,j}^{%
\mathbf{u},u_{h}^{\prime },u_{\ell }^{\prime \prime }}+\partial _{yy}^{2}%
\bar{\sigma}_{j}\left( t,\mathbf{X}_{t},\mathbf{u}_{t}\right) \star \mathcal{%
Y}_{t}^{\mathbf{u},u_{h}^{\prime },u_{\ell }^{\prime \prime }}+\bar{\Xi}%
_{t,j}^{\mathbf{u},u_{h}^{\prime },u_{\ell }^{\prime \prime }}\left( \bar{%
\sigma}\right) \right\rangle \mathrm{d}t\Bigg]. 
\end{aligned}
\end{eqnarray*}
%
Applying It\^{o}'s formula to (\ref{fy}) yields
\begin{eqnarray}
\mathrm{d}\left( \mathbf{Y}_{T}^{\mathbf{u},u_{\ell }^{\prime \prime
}}\cdot \left( \mathbf{Y}_{T}^{\mathbf{u},u_{h}^{\prime }}\right) ^{\top
}\right) 
& \!\!\!=\!\!\! & \mathbf{Y}_{t}^{\mathbf{u},u_{\ell }^{\prime \prime }}\Bigg [\left(
\left( \mathbf{Y}_{t}^{\mathbf{u},u_{h}^{\prime }}\right) ^{\top }\mathbf{B}%
_{0,x,y}\left( t,\mathbf{X}_{t}^{\mathbf{u}},\mathbf{u}_{t}\right) ^{\top }+%
\mathbf{B}_{1,u}^{h}\left( t,\mathbf{X}_{t}^{\mathbf{u}},\mathbf{u}%
_{t}\right) ^{\top }u_{t,h}^{\prime }\right) \mathrm{d}t  \notag \\
&&+\sum_{i=1}^{N}\left( \left( \mathbf{Y}_{t}^{\mathbf{u},u_{h}^{\prime
}}\right) ^{\top }\mathbf{\Pi }_{0,x,y}^{i}\left( t,\mathbf{X}_{t}^{\mathbf{u%
}},\mathbf{u}_{t}\right) ^{\top }+\mathbf{\Pi }_{1,u}^{i}\left( t,\mathbf{X}%
_{t}^{\mathbf{u}},\mathbf{u}_{t}\right) ^{\top }\delta _{i,h}u_{t,h}^{\prime
}\right) \mathrm{d}W_{t}^{i}\Bigg ]  \notag \\
&&+\Bigg [\left( \mathbf{B}_{0,x,y}\left( t,\mathbf{X}_{t}^{\mathbf{u}},%
\mathbf{u}_{t}\right) \mathbf{Y}_{t}^{\mathbf{u},u_{\ell }^{\prime \prime }}+%
\mathbf{B}_{1,u}^{\ell }\left( t,\mathbf{X}_{t}^{\mathbf{u}},\mathbf{u}%
_{t}\right) u_{t,\ell }^{\prime \prime }\right) \mathrm{d}t  \notag \\
&&+\sum_{i=1}^{N}\left( \mathbf{\Pi }_{0,x,y}^{i}\left( t,\mathbf{X}_{t}^{%
\mathbf{u}},\mathbf{u}_{t}\right) \mathbf{Y}_{t}^{\mathbf{u},u_{\ell
}^{\prime \prime }}+\mathbf{\Pi }_{1,u}^{i}\left( t,\mathbf{X}_{t}^{\mathbf{u%
}},\mathbf{u}_{t}\right) \delta _{i,\ell }u_{t,\ell }^{\prime \prime
}\right) \mathrm{d}W_{t}^{i}\Bigg ]\cdot \left( \mathbf{Y}_{T}^{\mathbf{u}%
,u_{h}^{\prime }}\right) ^{\top }  \notag \\
&&+\sum_{i=1}^{N}\Bigg [\left( \mathbf{\Pi }_{0,x,y}^{i}\left( t,\mathbf{X}%
_{t}^{\mathbf{u}},\mathbf{u}_{t}\right) \mathbf{Y}_{t}^{\mathbf{u},u_{\ell
}^{\prime \prime }}+\mathbf{\Pi }_{1,u}^{i}\left( t,\mathbf{X}_{t}^{\mathbf{u%
}},\mathbf{u}_{t}\right) \delta _{i,\ell }u_{t,\ell }^{\prime \prime
}\right)   \notag \\
&&\cdot \left( \left( \mathbf{Y}_{t}^{\mathbf{u},u_{h}^{\prime }}\right)
^{\top }\mathbf{\Pi }_{0,x,y}^{i}\left( t,\mathbf{X}_{t}^{\mathbf{u}},%
\mathbf{u}_{t}\right) ^{\top }+\mathbf{\Pi }_{1,u}^{i}\left( t,\mathbf{X}%
_{t}^{\mathbf{u}},\mathbf{u}_{t}\right) ^{\top }\delta _{i,h}u_{t,h}^{\prime
}\right) \Bigg ]\mathrm{d}t,  \label{y2}
\end{eqnarray}
Now applying Lemma \ref{l1} to (\ref{y2}) and (\ref{adj2}), together with
\begin{eqnarray*}
\begin{aligned}
\Phi _{t} =& \, \mathcal{Y}_{t}^{\mathbf{u},u_{h}^{\prime },u_{\ell }^{\prime
\prime }}\mathbf{B}_{0,x,y}\left( t,\mathbf{X}_{t}^{\mathbf{u}},\mathbf{u}%
_{t}\right) ^{\top }+\mathbf{B}_{0,x,y}\left( t,\mathbf{X}_{t}^{\mathbf{u}},%
\mathbf{u}_{t}\right) \mathcal{Y}_{t}^{\mathbf{u},u_{h}^{\prime },u_{\ell
}^{\prime \prime }} \\
&+\mathbf{Y}_{t}^{\mathbf{u},u_{\ell }^{\prime \prime }}\mathbf{B}%
_{1,u}^{h}\left( t,\mathbf{X}_{t}^{\mathbf{u}},\mathbf{u}_{t}\right) ^{\top
}u_{t,h}^{\prime }+\mathbf{B}_{1,u}^{\ell }\left( t,\mathbf{X}_{t}^{\mathbf{u%
}},\mathbf{u}_{t}\right) \left( \mathbf{Y}_{T}^{\mathbf{u},u_{h}^{\prime
}}\right) ^{\top }u_{t,\ell }^{\prime \prime } \\
&+\sum_{i=1}^{N}\mathbf{\Pi }_{0,x,y}^{i}\left( t,\mathbf{X}_{t}^{\mathbf{u}%
},\mathbf{u}_{t}\right) \mathcal{Y}_{t}^{\mathbf{u},u_{h}^{\prime },u_{\ell
}^{\prime \prime }}\mathbf{\Pi }_{0,x,y}^{i}\left( t,\mathbf{X}_{t}^{\mathbf{%
u}},\mathbf{u}_{t}\right) ^{\top } \\
&+\mathbf{\Pi }_{0,x,y}^{h}\left( t,\mathbf{X}_{t}^{\mathbf{u}},\mathbf{u}%
_{t}\right) \mathbf{Y}_{t}^{\mathbf{u},u_{\ell }^{\prime \prime }}\mathbf{%
\Pi }_{1,u}^{h}\left( t,\mathbf{X}_{t}^{\mathbf{u}},\mathbf{u}_{t}\right)
^{\top }u_{t,h}^{\prime } \\
&+\mathbf{\Pi }_{1,u}^{\ell }\left( t,\mathbf{X}_{t}^{\mathbf{u}},\mathbf{u}%
_{t}\right) \left( \mathbf{Y}_{t}^{\mathbf{u},u_{h}^{\prime }}\right) ^{\top
}\mathbf{\Pi }_{0,x,y}^{\ell }\left( t,\mathbf{X}_{t}^{\mathbf{u}},\mathbf{u}%
_{t}\right) ^{\top }u_{t,\ell }^{\prime \prime }, \\
\Psi _{t}^{j} = & \, \mathcal{Y}_{t}^{\mathbf{u},u_{h}^{\prime },u_{\ell
}^{\prime \prime }}\mathbf{\Pi }_{0,x,y}^{j}\left( t,\mathbf{X}_{t}^{\mathbf{%
u}},\mathbf{u}_{t}\right) ^{\top }+\mathbf{\Pi }_{0,x,y}^{j}\left( t,\mathbf{%
X}_{t}^{\mathbf{u}},\mathbf{u}_{t}\right) \mathcal{Y}_{t}^{\mathbf{u}%
,u_{h}^{\prime },u_{\ell }^{\prime \prime }} \\
&+\mathbf{Y}_{t}^{\mathbf{u},u_{\ell }^{\prime \prime }}\mathbf{\Pi }%
_{1,u}^{j}\left( t,\mathbf{X}_{t}^{\mathbf{u}},\mathbf{u}_{t}\right) ^{\top
}\delta _{j,h}u_{t,h}^{\prime }+\mathbf{\Pi }_{1,u}^{j}\left( t,\mathbf{X}%
_{t}^{\mathbf{u}},\mathbf{u}_{t}\right) \left( \mathbf{Y}_{T}^{\mathbf{u}%
,u_{h}^{\prime }}\right) ^{\top }\delta _{j,\ell }u_{t,\ell }^{\prime \prime
}, \\
\Theta_{t} = & -\digamma _{t,i}^{h,\ell }, \  \ 
\mathbf{Q}_{t}^{j} =  \, \mathcal{Q}_{t,i,j}^{h,\ell},  
\end{aligned}
\end{eqnarray*} 
yields the following (noting tr$[A_{n\times m}B_{m\times n}]=$tr$[B_{m\times n}A_{n\times m}],$ $m,n\in N$ and $Y_{0}=0$): 
\begin{equation*}
\mathbb{E}\left\{ \text{tr}\left[ \mathcal{P}_{T,i}^{h,\ell }\mathcal{Y}%
_{T}^{\mathbf{u},u_{h}^{\prime },u_{\ell }^{\prime \prime }}\right] \right\}
=\mathbb{E}\Bigg [\int_{0}^{T}\left\{ \text{tr}\left[ -\digamma
_{t,i}^{h,\ell }\mathcal{Y}_{t}^{\mathbf{u},u_{h}^{\prime },u_{\ell
}^{\prime \prime }}+\mathcal{P}_{t,i}^{h,\ell }\Phi _{t}+\sum_{j=1}^{N}%
\mathcal{Q}_{t,i,j}^{h,\ell }\Psi _{t}^{j}\right] \right\} \mathrm{d}t, 
\end{equation*}
with
\begin{eqnarray}
\begin{aligned}
\mathcal{P}_{t,i}^{h,\ell }\Phi _{t} = & \, \mathcal{P}_{t,i}^{h,\ell }\mathcal{Y%
}_{t}^{\mathbf{u},u_{h}^{\prime },u_{\ell }^{\prime \prime }}\mathbf{B}%
_{0,x,y}\left( t,\mathbf{X}_{t}^{\mathbf{u}},\mathbf{u}_{t}\right) ^{\top }+%
\mathcal{P}_{t,i}^{h,\ell }\mathbf{Y}_{t}^{\mathbf{u},u_{\ell }^{\prime
\prime }}\mathbf{B}_{1,u}^{h}\left( t,\mathbf{X}_{t}^{\mathbf{u}},\mathbf{u}%
_{t}\right) ^{\top }u_{t,h}^{\prime }  \notag \\
&+\mathcal{P}_{t,i}^{h,\ell }\mathbf{B}_{0,x,y}\left( t,\mathbf{X}_{t}^{%
\mathbf{u}},\mathbf{u}_{t}\right) \mathcal{Y}_{t}^{\mathbf{u},u_{h}^{\prime
},u_{\ell }^{\prime \prime }}+\mathcal{P}_{t,i}^{h,\ell }\mathbf{B}%
_{1,u}^{\ell }\left( t,\mathbf{X}_{t}^{\mathbf{u}},\mathbf{u}_{t}\right)
u_{t,\ell }^{\prime \prime }\left( \mathbf{Y}_{T}^{\mathbf{u},u_{h}^{\prime
}}\right) ^{\top }  \notag \\
&+\mathcal{P}_{t,i}^{h,\ell }\sum_{i=1}^{N}\Bigg [\left( \mathbf{\Pi }%
_{0,x,y}^{i}\left( t,\mathbf{X}_{t}^{\mathbf{u}},\mathbf{u}_{t}\right)
\mathbf{Y}_{t}^{\mathbf{u},u_{\ell }^{\prime \prime }}+\mathbf{\Pi }%
_{1,u}^{i}\left( t,\mathbf{X}_{t}^{\mathbf{u}},\mathbf{u}_{t}\right) \delta
_{i,\ell }u_{t,\ell }^{\prime \prime }\right) \cdot   \notag \\
&\left( \left( \mathbf{Y}_{t}^{\mathbf{u},u_{h}^{\prime }}\right) ^{\top }%
\mathbf{\Pi }_{0,x,y}^{i}\left( t,\mathbf{X}_{t}^{\mathbf{u}},\mathbf{u}%
_{t}\right) ^{\top }+\mathbf{\Pi }_{1,u}^{i}\left( t,\mathbf{X}_{t}^{\mathbf{%
u}},\mathbf{u}_{t}\right) ^{\top }\delta _{i,h}u_{t,h}^{\prime }\right) %
\Bigg ]  \notag \\
= & \, \mathcal{P}_{t,i}^{h,\ell }\mathcal{Y}_{t}^{\mathbf{u},u_{h}^{\prime
},u_{\ell }^{\prime \prime }}\mathbf{B}_{0,x,y}\left( t,\mathbf{X}_{t}^{%
\mathbf{u}},\mathbf{u}_{t}\right) ^{\top }+\mathcal{P}_{t,i}^{h,\ell }%
\mathbf{B}_{0,x,y}\left( t,\mathbf{X}_{t}^{\mathbf{u}},\mathbf{u}_{t}\right)
\mathcal{Y}_{t}^{\mathbf{u},u_{h}^{\prime },u_{\ell }^{\prime \prime }}
\notag \\
&+\mathcal{P}_{t,i}^{h,\ell }\mathbf{Y}_{t}^{\mathbf{u},u_{\ell }^{\prime
\prime }}\mathbf{B}_{1,u}^{h}\left( t,\mathbf{X}_{t}^{\mathbf{u}},\mathbf{u}%
_{t}\right) ^{\top }u_{t,h}^{\prime }+\mathcal{P}_{t,i}^{h,\ell }\mathbf{B}%
_{1,u}^{\ell }\left( t,\mathbf{X}_{t}^{\mathbf{u}},\mathbf{u}_{t}\right)
\left( \mathbf{Y}_{T}^{\mathbf{u},u_{h}^{\prime }}\right) ^{\top }u_{t,\ell
}^{\prime \prime }  \notag \\
&+\mathcal{P}_{t,i}^{h,\ell }\Bigg [\sum_{i=1}^{N}\mathbf{\Pi }%
_{0,x,y}^{i}\left( t,\mathbf{X}_{t}^{\mathbf{u}},\mathbf{u}_{t}\right)
\mathcal{Y}_{t}^{\mathbf{u},u_{h}^{\prime },u_{\ell }^{\prime \prime }}%
\mathbf{\Pi }_{0,x,y}^{i}\left( t,\mathbf{X}_{t}^{\mathbf{u}},\mathbf{u}%
_{t}\right) ^{\top }  \notag \\
&+\mathbf{\Pi }_{0,x,y}^{h}\left( t,\mathbf{X}_{t}^{\mathbf{u}},\mathbf{u}%
_{t}\right) \mathbf{Y}_{t}^{\mathbf{u},u_{\ell }^{\prime \prime }}\mathbf{%
\Pi }_{1,u}^{h}\left( t,\mathbf{X}_{t}^{\mathbf{u}},\mathbf{u}_{t}\right)
^{\top }u_{t,h}^{\prime }  \notag \\
&+\mathbf{\Pi }_{1,u}^{\ell }\left( t,\mathbf{X}_{t}^{\mathbf{u}},\mathbf{u}%
_{t}\right) \left( \mathbf{Y}_{t}^{\mathbf{u},u_{h}^{\prime }}\right) ^{\top
}\mathbf{\Pi }_{0,x,y}^{\ell }\left( t,\mathbf{X}_{t}^{\mathbf{u}},\mathbf{u}%
_{t}\right) ^{\top }u_{t,\ell }^{\prime \prime }\Bigg ]  \label{c1}
\end{aligned}
\end{eqnarray}
and
\begin{eqnarray}
\sum_{j=1}^{N}\mathcal{Q}_{t,i,j}^{h,\ell }\Psi _{t}^{j} & \!\!\!=\!\!\! & \sum_{j=1}^{N}%
\mathcal{Q}_{t,i,j}^{h,\ell }\cdot \Bigg [\mathbf{Y}_{t}^{\mathbf{u},u_{\ell
}^{\prime \prime }}\cdot \left( \left( \mathbf{Y}_{t}^{\mathbf{u}%
,u_{h}^{\prime }}\right) ^{\top }\mathbf{\Pi }_{0,x,y}^{j}\left( t,\mathbf{X}%
_{t}^{\mathbf{u}},\mathbf{u}_{t}\right) ^{\top }+\mathbf{\Pi }%
_{1,u}^{j}\left( t,\mathbf{X}_{t}^{\mathbf{u}},\mathbf{u}_{t}\right) ^{\top
}\delta _{j,h}u_{t,h}^{\prime }\right)   \notag \\
&&+\left( \mathbf{\Pi }_{0,x,y}^{j}\left( t,\mathbf{X}_{t}^{\mathbf{u}},%
\mathbf{u}_{t}\right) \mathbf{Y}_{t}^{\mathbf{u},u_{\ell }^{\prime \prime }}+%
\mathbf{\Pi }_{1,u}^{j}\left( t,\mathbf{X}_{t}^{\mathbf{u}},\mathbf{u}%
_{t}\right) \delta _{j,\ell }u_{t,\ell }^{\prime \prime }\right) \cdot
\left( \mathbf{Y}_{t}^{\mathbf{u},u_{h}^{\prime }}\right) ^{\top }\Bigg ]
\notag \\
& \!\!\!=\!\!\! & \sum_{j=1}^{N}\mathcal{Q}_{t,i,j}^{h,\ell }\cdot \Bigg [\mathcal{Y}_{t}^{%
\mathbf{u},u_{h}^{\prime },u_{\ell }^{\prime \prime }}\mathbf{\Pi }%
_{0,x,y}^{j}\left( t,\mathbf{X}_{t}^{\mathbf{u}},\mathbf{u}_{t}\right)
^{\top }+\mathbf{\Pi }_{0,x,y}^{j}\left( t,\mathbf{X}_{t}^{\mathbf{u}},%
\mathbf{u}_{t}\right) \mathcal{Y}_{t}^{\mathbf{u},u_{h}^{\prime },u_{\ell
}^{\prime \prime }}  \notag \\
&&+\mathbf{Y}_{t}^{\mathbf{u},u_{\ell }^{\prime \prime }}\mathbf{\Pi }%
_{1,u}^{j}\left( t,\mathbf{X}_{t}^{\mathbf{u}},\mathbf{u}_{t}\right) ^{\top
}\delta _{j,h}u_{t,h}^{\prime }+\mathbf{\Pi }_{1,u}^{j}\left( t,\mathbf{X}%
_{t}^{\mathbf{u}},\mathbf{u}_{t}\right) \left( \mathbf{Y}_{t}^{\mathbf{u}%
,u_{h}^{\prime }}\right) ^{\top }\delta _{j,\ell }u_{t,\ell }^{\prime \prime
}\Bigg ]  \notag \\
& \!\!\!=\!\!\! & \sum_{j=1}^{N}\mathcal{Q}_{t,i,j}^{h,\ell }\cdot \left( \mathcal{Y}_{t}^{%
\mathbf{u},u_{h}^{\prime },u_{\ell }^{\prime \prime }}\mathbf{\Pi }%
_{0,x,y}^{j}\left( t,\mathbf{X}_{t}^{\mathbf{u}},\mathbf{u}_{t}\right)
^{\top }+\mathbf{\Pi }_{0,x,y}^{j}\left( t,\mathbf{X}_{t}^{\mathbf{u}},%
\mathbf{u}_{t}\right) \mathcal{Y}_{t}^{\mathbf{u},u_{h}^{\prime },u_{\ell
}^{\prime \prime }}\right)   \notag \\
&&+\mathcal{Q}_{t,i,h}^{h,\ell }\cdot \mathbf{Y}_{t}^{\mathbf{u},u_{\ell
}^{\prime \prime }}\mathbf{\Pi }_{1,u}^{h}\left( t,\mathbf{X}_{t}^{\mathbf{u}%
},\mathbf{u}_{t}\right) ^{\top }u_{t,h}^{\prime }+\mathcal{Q}_{t,i,\ell
}^{h,\ell }\cdot \mathbf{\Pi }_{1,u}^{\ell }\left( t,\mathbf{X}_{t}^{\mathbf{%
u}},\mathbf{u}_{t}\right) \left( \mathbf{Y}_{t}^{\mathbf{u},u_{h}^{\prime
}}\right) ^{\top }u_{t,\ell }^{\prime \prime }.  \label{c2}
\end{eqnarray}
Note that equation (\ref{adj2}) is also a BSDE with matrix-valued unknowns.
As with (\ref{adj1}), under Assumptions (A1)-(A2), there exists a unique
adapted solution $\left( \mathcal{P}_{t,i}^{h,\ell },\left\{ \mathcal{Q}%
_{t,i,j}^{h,\ell }\right\} _{j=1,\ldots ,N}\right) _{0\leq t\leq T}$ to (\ref%
{adj2}).
Consequently, combining (\ref{c1}) and (\ref{c3}), we have
\begin{eqnarray*}
\mathbb{E}\left\{ \text{tr}\left[ \mathcal{P}_{T,i}^{h,\ell }\mathcal{Y}%
_{T}^{\mathbf{u},u_{h}^{\prime },u_{\ell }^{\prime \prime }}\right] \right\}
& \!\!\!=\!\!\! & \mathbb{E}\left[ \int_{0}^{T}\text{tr}\left[ -\digamma _{t,i}^{h,\ell }%
\mathcal{Y}_{t}^{\mathbf{u},u_{h}^{\prime },u_{\ell }^{\prime \prime }}+%
\mathcal{P}_{t,i}^{h,\ell }\Phi _{t}+\sum_{j=1}^{N}\mathcal{Q}_{t,i,j}^{h,\ell
}\Psi _{t}^{j}\right] \mathrm{d}t\right] \\
& \!\!\!=\!\!\! & \mathbb{E}\Bigg [\int_{0}^{T}\text{tr}\Bigg \{\mathcal{P}_{t,i}^{h,\ell }%
\Big [\mathbf{Y}_{t}^{\mathbf{u},u_{\ell }^{\prime \prime }}\mathbf{B}%
_{1,u}^{h}\left( t,\mathbf{X}_{t}^{\mathbf{u}},\mathbf{u}_{t}\right) ^{\top
}u_{t,h}^{\prime }+\mathbf{B}_{1,u}^{\ell }\left( t,\mathbf{X}_{t}^{\mathbf{u%
}},\mathbf{u}_{t}\right) \left( \mathbf{Y}_{T}^{\mathbf{u},u_{h}^{\prime
}}\right) ^{\top }u_{t,\ell }^{\prime \prime } \\
&&+\mathbf{\Pi }_{0,x,y}^{h}\left( t,\mathbf{X}_{t}^{\mathbf{u}},\mathbf{u}%
_{t}\right) \mathbf{Y}_{t}^{\mathbf{u},u_{\ell }^{\prime \prime }}\mathbf{%
\Pi }_{1,u}^{h}\left( t,\mathbf{X}_{t}^{\mathbf{u}},\mathbf{u}_{t}\right)
^{\top }u_{t,h}^{\prime } \\
&&+\mathbf{\Pi }_{1,u}^{\ell }\left( t,\mathbf{X}_{t}^{\mathbf{u}},\mathbf{u}%
_{t}\right) \left( \mathbf{Y}_{t}^{\mathbf{u},u_{h}^{\prime }}\right) ^{\top
}\mathbf{\Pi }_{0,x,y}^{\ell }\left( t,\mathbf{X}_{t}^{\mathbf{u}},\mathbf{u}%
_{t}\right) ^{\top }u_{t,\ell }^{\prime \prime }\Big ]-\partial
_{yy}^{2}f_{i}\mathcal{Y}_{t}^{\mathbf{u},u_{h}^{\prime },u_{\ell }^{\prime
\prime }} \\
&&-\left\langle P_{t,i},\partial _{yy}^{2}\bar{b}\left( t,\mathbf{X}_{t},%
\mathbf{u}_{t}\right) \star \mathcal{Y}_{t}^{\mathbf{u},u_{h}^{\prime
},u_{\ell }^{\prime \prime }}\right\rangle -\sum_{j=1}^{N}\left\langle
Q_{t,i,j},\partial _{yy}^{2}\bar{\sigma}_{j}\left( t,\mathbf{X}_{t},\mathbf{u}%
_{t}\right) \star \mathcal{Y}_{t}^{\mathbf{u},u_{h}^{\prime },u_{\ell
}^{\prime \prime }}\right\rangle \\
&&+\mathcal{Q}_{t,i,h}^{h,\ell }\cdot \mathbf{Y}_{t}^{\mathbf{u},u_{\ell
}^{\prime \prime }}\mathbf{\Pi }_{1,u}^{h}\left( t,\mathbf{X}_{t}^{\mathbf{u}%
},\mathbf{u}_{t}\right) ^{\top }u_{t,h}^{\prime }+\mathcal{Q}_{t,i,\ell
}^{h,\ell }\cdot \mathbf{\Pi }_{1,u}^{\ell }\left( t,\mathbf{X}_{t}^{\mathbf{%
u}},\mathbf{u}_{t}\right) \left( \mathbf{Y}_{t}^{\mathbf{u},u_{h}^{\prime
}}\right) ^{\top }u_{t,\ell }^{\prime \prime }\Bigg \}\mathrm{d}t.
\end{eqnarray*}
From (\ref{aerfa}) for $\frac{\delta ^{2}V_{i}}{\delta u_{h}\delta u_{\ell}}$, one can get
\begin{eqnarray}
\frac{\delta ^{2}V_{i}}{\delta u_{h}\delta u_{\ell }}\left( \mathbf{u}%
;u_{h}^{\prime },u_{\ell }^{\prime \prime }\right) 
& \!\!\!=\!\!\! & \mathbb{E}\left[ \left( \mathbf{Y}_{T}^{\mathbf{u},u_{h}^{\prime
}}\right) ^{\top }\partial _{xx}^{2}g_{i}\left( \mathbf{X}_{T}^{\mathbf{u}%
}\right) \mathbf{Y}_{T}^{\mathbf{u},u_{\ell }^{\prime \prime }}+\left(
\partial _{x}g_{i}\right) ^{\top }\left( \mathbf{X}_{t}^{\mathbf{u}}\right)
\mathbf{Z}_{T}^{\mathbf{u},u_{h}^{\prime },u_{\ell }^{\prime \prime }}\right]
\notag \\
&&+\mathbb{E}\Bigg \{\int_{0}^{T}\Big [\text{tr}\left[ \partial
_{xx}^{2}f_{i}\cdot \mathcal{Y}_{t}^{\mathbf{u},u_{h}^{\prime },u_{\ell
}^{\prime \prime }}\right] +\left( \partial _{x}f_{i}\right) ^{\top }\left(
t,\mathbf{X}_{t}^{\mathbf{u}},\mathbf{u}_{t}\right) \mathbf{Z}_{t}^{\mathbf{u%
},u_{h}^{\prime },u_{\ell }^{\prime \prime }}  \notag \\
&&+\partial _{u_{h}x}^{2}f_{i}\cdot \mathbf{Y}_{t}^{\mathbf{u},u_{\ell
}^{\prime \prime }}\cdot u_{t,h}^{\prime }+\left( \mathbf{Y}_{t}^{\mathbf{u}%
,u_{h}^{\prime }}\right) ^{\top }\cdot \partial _{xu_{\ell }}^{2}f_{i}\cdot
u_{t,\ell }^{\prime \prime }  \notag \\
&&+\partial _{u_{h}u_{\ell }}^{2}f_{i}\cdot u_{t,h}^{\prime }\cdot u_{t,\ell
}^{\prime \prime }\Big ]\mathrm{d}t\Bigg\}.  \label{vsecl}
\end{eqnarray}
Meanwhile,
\begin{eqnarray*}
&&\mathbb{E}\left\{ \text{tr}\left[ \partial _{xx}^{2}g_{i}\left( \mathbf{X}%
_{T}^{\mathbf{u}}\right) \mathcal{Y}_{T}^{\mathbf{u},u_{h}^{\prime },u_{\ell
}^{\prime \prime }}\right] \right\} +\mathbb{E}\left[ \left\langle \mathbf{Z}%
_{T}^{\mathbf{u},u_{h}^{\prime },u_{\ell }^{\prime \prime }},\partial
_{x}g_{i}\left( \mathbf{X}_{T}^{\mathbf{u}}\right) \right\rangle \right]  \\
&&+\mathbb{E}\left[ \int_{0}^{T}\left( \left( \partial _{x}f_{i}\right)
^{\top }\left( t,\mathbf{X}_{t}^{\mathbf{u}},\mathbf{u}_{t}\right) \mathbf{Z}%
_{t}^{\mathbf{u},u_{h}^{\prime },u_{\ell }^{\prime \prime }}+\text{tr}\left[
\partial _{xx}^{2}f_{i}\cdot \mathcal{Y}_{t}^{\mathbf{u},u_{h}^{\prime
},u_{\ell }^{\prime \prime }}\right] \right) \mathrm{d}t\right]  \\
& =\!\!\! &\mathbb{E}\Bigg [\int_{0}^{T}\text{tr}\Bigg \{\mathcal{P}_{t,i}^{h,\ell }%
\Big [\mathbf{Y}_{t}^{\mathbf{u},u_{\ell }^{\prime \prime }}\mathbf{B}%
_{1,u}^{h}\left( t,\mathbf{X}_{t}^{\mathbf{u}},\mathbf{u}_{t}\right) ^{\top
}u_{t,h}^{\prime }+\mathbf{B}_{1,u}^{\ell }\left( t,\mathbf{X}_{t}^{\mathbf{u%
}},\mathbf{u}_{t}\right) \left( \mathbf{Y}_{T}^{\mathbf{u},u_{h}^{\prime
}}\right) ^{\top }u_{t,\ell }^{\prime \prime } \\
&&+\mathbf{\Pi }_{0,x,y}^{h}\left( t,\mathbf{X}_{t}^{\mathbf{u}},\mathbf{u}%
_{t}\right) \mathbf{Y}_{t}^{\mathbf{u},u_{\ell }^{\prime \prime }}\mathbf{%
\Pi }_{1,u}^{h}\left( t,\mathbf{X}_{t}^{\mathbf{u}},\mathbf{u}_{t}\right)
^{\top }u_{t,h}^{\prime } \\
&&+\mathbf{\Pi }_{1,u}^{\ell }\left( t,\mathbf{X}_{t}^{\mathbf{u}},\mathbf{u}%
_{t}\right) \left( \mathbf{Y}_{t}^{\mathbf{u},u_{h}^{\prime }}\right) ^{\top
}\mathbf{\Pi }_{0,x,y}^{\ell }\left( t,\mathbf{X}_{t}^{\mathbf{u}},\mathbf{u}%
_{t}\right) ^{\top }u_{t,\ell }^{\prime \prime }\Big ] \\
&&+\mathcal{Q}_{t,i,h}^{h,\ell }\cdot \mathbf{Y}_{t}^{\mathbf{u},u_{\ell
}^{\prime \prime }}\mathbf{\Pi }_{1,u}^{h}\left( t,\mathbf{X}_{t}^{\mathbf{u}%
},\mathbf{u}_{t}\right) ^{\top }u_{t,h}^{\prime }+\mathcal{Q}_{t,i,\ell
}^{h,\ell }\cdot \mathbf{\Pi }_{1,u}^{\ell }\left( t,\mathbf{X}_{t}^{\mathbf{%
u}},\mathbf{u}_{t}\right) \left( \mathbf{Y}_{t}^{\mathbf{u},u_{h}^{\prime
}}\right) ^{\top }u_{t,\ell }^{\prime \prime } \\
&&+\left\langle P_{t,i},\mathfrak{F}_{t}^{\mathbf{u},u_{h}^{\prime },u_{\ell
}^{\prime \prime }}+\Gamma _{t}^{\mathbf{u},u_{h}^{\prime },u_{\ell
}^{\prime \prime }}\left( \bar{b}\right) \right\rangle
+\sum_{j=1}^{N}\left\langle Q_{t,i,j},\mathfrak{G}_{t,j}^{\mathbf{u}%
,u_{h}^{\prime },u_{\ell }^{\prime \prime }}+\bar{\Xi}_{t,j}^{\mathbf{u}%
,u_{h}^{\prime },u_{\ell }^{\prime \prime }}\left( \bar{\sigma}\right)
\right\rangle \Bigg \}\mathrm{d}t.
\end{eqnarray*}%
Hence,
\begin{eqnarray*}
\frac{\delta ^{2}V_{i}}{\delta u_{h}\delta u_{\ell }}\left( \mathbf{u}%
;u_{h}^{\prime },u_{\ell }^{\prime \prime }\right)  
&=&\mathbb{E}\Bigg \{\int_{0}^{T}\text{tr}\Bigg [\mathcal{P}_{t,i}^{h,\ell }%
\Bigg (\mathbf{Y}_{t}^{\mathbf{u},u_{\ell }^{\prime \prime }}\mathbf{B}%
_{1,u}^{h}\left( t,\mathbf{X}_{t}^{\mathbf{u}},\mathbf{u}_{t}\right) ^{\top
}u_{t,h}^{\prime } \\
&&+\mathbf{B}_{1,u}^{\ell }\left( t,\mathbf{X}_{t}^{\mathbf{u}},\mathbf{u}%
_{t}\right) \left( \mathbf{Y}_{T}^{\mathbf{u},u_{h}^{\prime }}\right) ^{\top
}u_{t,\ell }^{\prime \prime }+\mathbf{\Pi }_{0,x,y}^{h}\left( t,\mathbf{X}%
_{t}^{\mathbf{u}},\mathbf{u}_{t}\right) \mathbf{Y}_{t}^{\mathbf{u},u_{\ell
}^{\prime \prime }}\mathbf{\Pi }_{1,u}^{h}\left( t,\mathbf{X}_{t}^{\mathbf{u}%
},\mathbf{u}_{t}\right) ^{\top }u_{t,h}^{\prime } \\
&&+\mathbf{\Pi }_{1,u}^{\ell }\left( t,\mathbf{X}_{t}^{\mathbf{u}},\mathbf{u}%
_{t}\right) \left( \mathbf{Y}_{t}^{\mathbf{u},u_{h}^{\prime }}\right) ^{\top
}\mathbf{\Pi }_{0,x,y}^{\ell }\left( t,\mathbf{X}_{t}^{\mathbf{u}},\mathbf{u}%
_{t}\right) ^{\top }u_{t,\ell }^{\prime \prime }\Bigg ) \\
&&+\mathcal{Q}_{t,i,h}^{h,\ell }\cdot \mathbf{Y}_{t}^{\mathbf{u},u_{\ell
}^{\prime \prime }}\mathbf{\Pi }_{1,u}^{h}\left( t,\mathbf{X}_{t}^{\mathbf{u}%
},\mathbf{u}_{t}\right) ^{\top }u_{t,h}^{\prime }+\mathcal{Q}_{t,i,\ell
}^{h,\ell }\cdot \mathbf{\Pi }_{1,u}^{\ell }\left( t,\mathbf{X}_{t}^{\mathbf{%
u}},\mathbf{u}_{t}\right) \left( \mathbf{Y}_{t}^{\mathbf{u},u_{h}^{\prime
}}\right) ^{\top }u_{t,\ell }^{\prime \prime } \\
&&+\partial _{u_{h}x}^{2}f_{i}\cdot \mathbf{Y}_{t}^{\mathbf{u},u_{\ell
}^{\prime \prime }}\cdot u_{t,h}^{\prime }+\left( \mathbf{Y}_{t}^{\mathbf{u}%
,u_{h}^{\prime }}\right) ^{\top }\cdot \partial _{xu_{\ell }}^{2}f_{i}\cdot
u_{t,\ell }^{\prime \prime }+\partial _{u_{h}u_{\ell }}^{2}f_{i}\cdot
u_{t,h}^{\prime }\cdot u_{t,\ell }^{\prime \prime }\Bigg ]\mathrm{d}t\Bigg \}%
,
\end{eqnarray*}
from which we get the desired result. \hfill $\Box $

\paragraph{Proof of Theorem \protect\ref{the1}.}

The proof of Theorem \ref{the1} relies on the precise estimates of the sensitivity process $Y^{\mathbf{u}, u_{i}^{\prime}}$, which quantify the dependence of each player's state process on the variations of other players' controls,  the coupling strengths $L_{y}^{b}, L_{y}^{\sigma}$ in the drift and diffusion coefficients, as well as the total number of players.

For notational simplicity we omit the dependence on $u$ in the superscript
of all processes, i.e., $X=X^{\mathbf{u}}$, $Y^{i}=Y^{\mathbf{u}%
,u_{i}^{\prime }},$ $Z^{i,j}=Z^{\mathbf{u},u_{i}^{\prime },u_{j}^{\prime
\prime }}$. We denote by $C>0$ a generic constant depending only on the
upper bounds of $T$, $\max_{i\in I_{N}}\mathbb{E}\left[ \left\vert \xi
_{i}\right\vert ^{2}\right] $, $\max_{i\in I_{N}}\left\Vert u_{i}\right\Vert
_{\mathcal{H}^{2}\left( \mathbb{R}\right) }$. 
Now, note that $Z^{\mathbf{u},u_{i}^{\prime },u_{j}^{\prime \prime }}=Z^{\mathbf{u%
},u_{j}^{\prime \prime },u_{i}^{\prime }}$, 
\begin{eqnarray}
&&\frac{\delta ^{2}V_{i}}{\delta u_{i}\delta u_{j}}\left( \mathbf{u}%
;u_{i}^{\prime },u_{j}^{\prime \prime }\right) -\frac{\delta ^{2}V_{j}}{%
\delta u_{j}\delta u_{i}}\left( \mathbf{u};u_{j}^{\prime \prime
},u_{i}^{\prime }\right)   \notag \\
=\!\!\! &&\mathbb{E}\Bigg [\int_{0}^{T}\Bigg \{\left(
\begin{array}{c}
\mathbf{Y}_{t}^{i} \\
u_{t,i}^{\prime }%
\end{array}%
\right) ^{\top }\left(
\begin{array}{cc}
\partial _{xx}^{2}\Delta _{i,j}^{f} & \partial _{xu_{j}}^{2}\Delta _{i,j}^{f}
\\
\partial _{u_{i}x}^{2}\Delta _{i,j}^{f} & \partial _{u_{i}u_{j}}^{2}\Delta
_{i,j}^{f}%
\end{array}%
\right) \left( t,\cdot \right) \left(
\begin{array}{c}
\mathbf{Y}_{t}^{j} \\
u_{t,j}^{\prime \prime }%
\end{array}%
\right) +\left( \mathbf{Z}_{t}^{i,j}\right) ^{\top }\left( \partial
_{x}\Delta _{i,j}^{f}\right) \left( t,\cdot \right) \Bigg \}\mathrm{d}t\Bigg
]  \notag \\
&&+\mathbb{E}\left[ \left( \mathbf{Y}_{T}^{i}\right) ^{\top }\left( \partial
_{xx}^{2}\Delta _{i,j}^{g}\right) \left( \mathbf{X}_{T}\right) \mathbf{Y}%
_{T}^{j}+\left( \mathbf{Z}_{T}^{i,j}\right) ^{\top }\left( \partial
_{x}\Delta _{i,j}^{g}\right) \left( \mathbf{X}_{T}\right) \right] .
\label{vz}
\end{eqnarray}
%
Consider the following BSDEs
\begin{equation}
\left\{\begin{aligned}
-\mathrm{d}P_{t}^{i,j} = & \, \Big [\mathbf{B}_{0,x,y}\left( t,\mathbf{X}_{t},\mathbf{u}_{t}\right)^{\top}P_{t}^{i,j}+\sum_{k=1}^{N}\mathbf{\Pi}_{0,x,y}^{k}\left(t,\mathbf{X}_{t},\mathbf{u}_{t}\right) ^{\top
}Q_{t,k}^{i,j} \\
& +\left( \partial _{x}\Delta _{i,j}^{f}\right) \left(t,\mathbf{X}_{t},\mathbf{u}_{t}\right) \Big]\mathrm{d}t-\sum_{k=1}^{N}Q_{t,k}^{i,j}\mathrm{d}W_{t}^{k}, \\
P_{T}^{i,j} = & \, \partial _{x}\Delta _{i,j}^{g}\left( \mathbf{X}_{T}\right),\text{ }\forall i,j\in I_{N}\text{ with }i\neq j%
\end{aligned}\right.   \label{adjn1}
\end{equation}
and
\begin{equation}
\left\{\begin{aligned}
-\mathrm{d}\mathcal{P}_{t}^{i,j} = & \Big [\mathbf{B}_{0,x,y}\left( t,%
\mathbf{X}_{t},\mathbf{u}_{t}\right) ^{\top }\mathcal{P}_{t}^{i,j}+\mathcal{P%
}_{t}^{i,j}\mathbf{B}_{0,x,y}\left( t,\mathbf{X}_{t},\mathbf{u}_{t}\right)
\\
&  +\sum_{k=1}^{N}\mathbf{\Pi }_{0,x,y}^{k}\left( t,\mathbf{X}_{t},\mathbf{%
u}_{t}\right) ^{\top }\mathcal{P}_{t}^{i,j}\mathbf{\Pi }_{0,x,y}^{k}\left( t,%
\mathbf{X}_{t},\mathbf{u}_{t}\right)  \\
& +\sum_{k=1}^{N}\left( \mathbf{\Pi }_{0,x,y}^{k}\left( t,\mathbf{X}_{t},%
\mathbf{u}_{t}\right) ^{\top }\mathcal{Q}_{t,k}^{i,j}+\mathcal{Q}_{t,k}^{i,j}%
\mathbf{\Pi }_{0,x,y}^{k}\left( t,\mathbf{X}_{t},\mathbf{u}_{t}\right)
\right)  \\
&  +\partial _{yy}^{2}\mathbb{H}_{i,j}\left( t,\mathbf{X}_{t},\mathbf{u}%
_{t},P_{t}^{i,j},Q_{t}^{i,j}\right) \Big ]\mathrm{d}t-\sum_{k=1}^{N}\mathcal{%
Q}_{t,k}^{i,j}\mathrm{d}W_{t}^{k}, \\
\mathcal{P}_{T}^{i,j} = & \, \partial _{xx}\Delta _{i,j}^{g}\left( \mathbf{X}_{T}\right), \text{ }\forall i,j\in I_{N}\text{ with }i\neq j.%
\end{aligned}\right.   \label{adjn2}
\end{equation}%
By (\ref{Vi}), we have
\begin{eqnarray}
&&\frac{\delta ^{2}V_{i}}{\delta u_{i}\delta u_{j}}\left( \mathbf{u}%
;u_{i}^{\prime },u_{j}^{\prime \prime }\right) -\frac{\delta ^{2}V_{j}}{%
\delta u_{j}\delta u_{i}}\left( \mathbf{u};u_{j}^{\prime \prime
},u_{i}^{\prime }\right)   \notag \\
& =\!\!\! & \mathbb{E}\Bigg \{\int_{0}^{T}\text{tr}\Bigg [\Big (\mathcal{P}%
_{t}^{i,j}\mathbf{Y}_{t}^{j}\mathbf{B}_{1,u}^{i}\left( t,\mathbf{X}_{t},%
\mathbf{u}_{t}\right) ^{\top }+\mathcal{P}_{t}^{i,j}\mathbf{\Pi }%
_{0,x,y}^{i}\left( t,\mathbf{X}_{t},\mathbf{u}_{t}\right) \mathbf{Y}_{t}^{j}%
\mathbf{\Pi }_{1,u}^{i}\left( t,\mathbf{X}_{t},\mathbf{u}_{t}\right) ^{\top }
\notag \\
&&+\mathcal{Q}_{t,i}^{i,j}\cdot \mathbf{Y}_{t}^{j}\mathbf{\Pi }%
_{1,u}^{i}\left( t,\mathbf{X}_{t},\mathbf{u}_{t}\right) ^{\top }+\partial
_{u_{i}x}^{2}\Delta _{i,j}^{f}\cdot \mathbf{Y}_{t}^{j}\Big )u_{t,i}^{\prime }
\notag \\
&&+\Big (\mathcal{P}_{t}^{i,j}\mathbf{B}_{1,u}^{j}\left( t,\mathbf{X}_{t},%
\mathbf{u}_{t}\right) \left( \mathbf{Y}_{t}^{i}\right) ^{\top }+\mathcal{P}%
_{t}^{i,j}\mathbf{\Pi }_{1,u}^{j}\left( t,\mathbf{X}_{t},\mathbf{u}%
_{t}\right) \left( \mathbf{Y}_{t}^{i}\right) ^{\top }\mathbf{\Pi }%
_{0,x,y}^{j}\left( t,\mathbf{X}_{t},\mathbf{u}_{t}\right) ^{\top }  \notag \\
&&+\mathcal{Q}_{t,j}^{i,j}\cdot \mathbf{\Pi }_{1,u}^{j}\left( t,\mathbf{X}%
_{t},\mathbf{u}_{t}\right) \left( \mathbf{Y}_{t}^{i}\right) ^{\top }+\left(
\mathbf{Y}_{t}^{i}\right) ^{\top }\cdot \partial _{xu_{j}}^{2}\Delta
_{i,j}^{f}\Big )u_{t,j}^{\prime \prime }  \notag \\
&&+\partial _{u_{i}u_{j}}^{2}\Delta _{i,j}^{f}\cdot u_{t,i}^{\prime }\cdot
u_{t,j}^{\prime \prime }+\left\langle P_{t}^{i,j},\mathfrak{F}_{t}^{i,j}+\Gamma _{t}^{i,j}\left( \bar{b}\right) \right\rangle
+\sum_{k=1}^{N}\left\langle Q_{t,k}^{i,j},\mathfrak{G}_{t,k}^{i,j}+\bar{\Xi}_{t,k}^{i,j}\left( \bar{\sigma}\right) \right\rangle \Bigg ]\mathrm{d}t\Bigg\},  \label{mineq1}
\end{eqnarray}%
where $\mathbf{B}_{1,u}^{i},\mathbf{\Pi }_{0,x,y}^{i},\mathbf{\Pi }%
_{1,u}^{i},\mathfrak{F}_{t}^{i,j},\Gamma _{t}^{i,j}\left( \bar{b}\right) ,%
\mathfrak{G}_{t,k}^{i,j},\bar{\Xi}_{t,k}^{i,j}\left( \bar{\sigma}\right) $
are defined in (\ref{bi1u}), (\ref{pii0y}), (\ref{pij1u}), (\ref{ff}), (\ref{gama}), (\ref{fb}),  (\ref{str}), respectively.

In order to obtain the estimation of (\ref{mineq1}), observe that 
\begin{eqnarray}
\mathbb{E}\left[ \left\langle \mathbf{Z}_{T}^{i,j},\partial_{x}\Delta_{i,j}^{g}\left( \mathbf{X}_{T}\right) \right\rangle \right]  
& \!\!\!=\!\!\! & \mathbb{E}\Bigg [-\int_{0}^{T}\Big <\mathbf{Z}_{t}^{i,j},\left( \partial_{x}\Delta _{i,j}^{f}\right) \left( t,\cdot \right) \Big >\mathrm{d}t  \notag
\\
&&+\int_{0}^{T}\left\langle P_{t}^{i,j},\mathfrak{F}_{t}^{i,j}+\partial
_{yy}^{2}\bar{b}\left( t,\mathbf{X}_{t},\mathbf{u}_{t}\right) \star \mathcal{%
Y}_{t}^{i,j}+\Gamma _{t}^{i,j}\left( \bar{b}\right) \right\rangle \mathrm{d}t
\notag \\
&&+\int_{0}^{T}\sum_{j=1}^{N}\left\langle Q_{t}^{i,j},\mathfrak{G}%
_{t,j}^{i,j}+\partial _{yy}^{2}\bar{\sigma}_{j}\left( t,\mathbf{X}_{t},%
\mathbf{u}_{t}\right) \star \mathcal{Y}_{t}^{i,j}+\bar{\Xi}%
_{t,j}^{i,j}\left( \bar{\sigma}\right) \right\rangle \mathrm{d}t\Bigg ].
\label{zij}
\end{eqnarray}
Substituting (\ref{zij}) into (\ref{vz}), 
\begin{eqnarray}
& & \frac{\delta ^{2}V_{i}}{\delta u_{i}\delta u_{j}}\left( \mathbf{u};u_{i}^{\prime },u_{j}^{\prime \prime }\right) -\frac{\delta ^{2}V_{j}}{\delta u_{j}\delta u_{i}}\left( \mathbf{u};u_{j}^{\prime \prime},u_{i}^{\prime }\right) \\ 
& =\!\!\! & \mathbb{E}\Bigg [\left( \mathbf{Y}_{T}^{i}\right) ^{\top }\left( \partial
_{xx}^{2}\Delta _{i,j}^{g}\right) \left( \mathbf{X}_{T}\right) \mathbf{Y}%
_{T}^{j}  \notag \\
&&+\int_{0}^{T}\left(
\begin{array}{c}
\mathbf{Y}_{t}^{i} \\
u_{t,i}^{\prime }%
\end{array}%
\right) ^{\top }\left(
\begin{array}{cc}
\partial _{xx}^{2}\Delta _{i,j}^{f} & \partial _{xu_{j}}^{2}\Delta _{i,j}^{f}
\\
\partial _{u_{i}x}^{2}\Delta _{i,j}^{f} & \partial _{u_{i}u_{j}}^{2}\Delta
_{i,j}^{f}%
\end{array}%
\right) \left( t,\cdot \right) \left(
\begin{array}{c}
\mathbf{Y}_{t}^{j} \\
u_{t,j}^{\prime \prime }%
\end{array}%
\right) \mathrm{d}t\Bigg ]  \notag \\
&&+\mathbb{E}\left[ \left( \mathbf{Z}_{T}^{i,j}\right) ^{\top }\left(
\partial _{x}\Delta _{i,j}^{g}\right) \left( \mathbf{X}_{T}\right)
+\int_{0}^{T}\left( \mathbf{Z}_{t}^{i,j}\right) ^{\top }\left( \partial
_{x}\Delta _{i,j}^{f}\right) \left( t,\cdot \right) \mathrm{d}t\right]
\notag \\
& =\!\!\! & \mathbb{E}\Bigg [\left( \mathbf{Y}_{T}^{i}\right) ^{\top }\left( \partial
_{xx}^{2}\Delta _{i,j}^{g}\right) \left( \mathbf{X}_{T}\right) \mathbf{Y}%
_{T}^{j}  \notag \\
&&+\int_{0}^{T}\left(
\begin{array}{c}
\mathbf{Y}_{t}^{i} \\
u_{t,i}^{\prime }%
\end{array}%
\right) ^{\top }\left(
\begin{array}{cc}
\partial _{xx}^{2}\Delta _{i,j}^{f} & \partial _{xu_{j}}^{2}\Delta _{i,j}^{f}
\\
\partial _{u_{i}x}^{2}\Delta _{i,j}^{f} & \partial _{u_{i}u_{j}}^{2}\Delta_{i,j}^{f}%
\end{array}%
\right) \left( t,\cdot \right) \left(
\begin{array}{c}
\mathbf{Y}_{t}^{j} \\
u_{t,j}^{\prime \prime }%
\end{array}%
\right) \mathrm{d}t\Bigg ]  \notag \\
&&+\mathbb{E}\Bigg [\int_{0}^{T}\left\langle P_{t}^{i,j},\mathfrak{F}%
_{t}^{i,j}+\partial _{yy}^{2}\bar{b}\left( t,\mathbf{X}_{t},\mathbf{u}%
_{t}\right) \star \mathcal{Y}_{t}^{i,j}+\Gamma _{t}^{i,j}\left( \bar{b}%
\right) \right\rangle \mathrm{d}t  \notag \\
&&+\int_{0}^{T}\sum_{k=1}^{N}\left\langle Q_{t,k}^{i,j},\mathfrak{G}%
_{t,k}^{i,j}+\partial _{yy}^{2}\bar{\sigma}_{k}\left( t,\mathbf{X}_{t},%
\mathbf{u}_{t}\right) \star \mathcal{Y}_{t}^{i,j}+\bar{\Xi}%
_{t,k}^{i,j}\left( \bar{\sigma}\right) \right\rangle \mathrm{d}t\Bigg ].
\label{newv1}
\end{eqnarray}
We can estimate the first and second terms from the end of (\ref{newv1}) as
in \cite{GLZ1}:
\begin{eqnarray}
&&\left\vert \mathbb{E}\Bigg [\int_{0}^{T}\left(
\begin{array}{c}
\mathbf{Y}_{t}^{i} \\
u_{t,i}^{\prime }%
\end{array}%
\right) ^{\top }\left(
\begin{array}{cc}
\partial _{xx}^{2}\Delta _{i,j}^{f} & \partial _{xu_{j}}^{2}\Delta _{i,j}^{f}
\\
\partial _{u_{i}x}^{2}\Delta _{i,j}^{f} & \partial _{u_{i}u_{j}}^{2}\Delta
_{i,j}^{f}%
\end{array}%
\right) \left( t,\cdot \right) \left(
\begin{array}{c}
\mathbf{Y}_{t}^{j} \\
u_{t,j}^{\prime \prime }%
\end{array}%
\right) \mathrm{d}t\Bigg ]\right\vert   \notag \\
& \leq\!\!\! & C\Bigg \{\left\Vert \partial _{x_{i}x_{j}}^{2}\Delta^f_{i,j}\right\Vert _{L^{\infty }}+\left\Vert \partial
_{x_{i}u_{j}}^{2}\Delta^f_{i,j}\right\Vert _{L^{\infty }}+\left\Vert
\partial _{u_{i}x_{j}}^{2}\Delta^f_{i,j}\right\Vert _{L^{\infty
}}+\left\Vert \partial _{u_{i}u_{j}}^{2}\Delta^f_{i,j}\right\Vert
_{L^{\infty }}  \notag \\
&&+\frac{L_{y}^{b,\sigma }}{N}\Bigg (\sum_{\ell \in I_{N}\backslash \left\{
j\right\} }\left( \left\Vert \partial _{x_{i}x_{\ell }}^{2}\Delta^f_{i,j}\right\Vert _{L^{\infty }}+\left\Vert \partial _{u_{i}x_{\ell
}}^{2}\Delta^f_{i,j}\right\Vert _{L^{\infty }}\right)   \notag \\
&&+\sum_{h\in I_{N}\backslash \left\{ i\right\} }\left( \left\Vert \partial
_{x_{h}x_{j}}^{2}\Delta^f_{i,j}\right\Vert _{L^{\infty }}+\left\Vert
\partial _{x_{h}u_{j}}^{2}\Delta^f_{i,j}\right\Vert _{L^{\infty }}\right) %
\Bigg )  \notag \\
&&+\frac{\left( L_{y}^{b,\sigma }\right) ^{2}}{N^{2}}\left( \sum_{\ell \in
I_{N}\backslash \left\{ j\right\} ,h\in I_{N}\backslash \left\{ i\right\}
}\left\Vert \partial _{x_{h}x_{\ell }}^{2}\Delta^f_{i,j}\right\Vert
_{L^{\infty }}\right) \Bigg \}\left\Vert u_{i}^{\prime }\right\Vert _{%
\mathcal{H}^{2}\left( \mathbb{R}\right) }\left\Vert u_{j}^{\prime \prime
}\right\Vert _{\mathcal{H}^{2}\left( \mathbb{R}\right) }  \label{m1}
\end{eqnarray}
and
\begin{eqnarray}
&&\left\vert \mathbb{E}\left[ \left( \mathbf{Y}_{T}^{i}\right) ^{\top
}\left( \partial _{xx}^{2}\Delta _{i,j}^{g}\right) \left( \mathbf{X}%
_{T}\right) \mathbf{Y}_{T}^{j}\right] \right\vert   \notag \\
& \leq\!\!\! & C\Bigg [\left\Vert \partial _{x_{i}x_{j}}^{2}\Delta^g_{i,j}\right\Vert _{L^{\infty }}+\frac{L_{y}^{b,\sigma }}{N}\Bigg (%
\sum_{\ell \in I_{N}\backslash \left\{ j\right\} }\left\Vert \partial
_{x_{i}x_{\ell }}^{2}\Delta^g_{i,j}\right\Vert _{L^{\infty }}+\sum_{h\in
I_{N}\backslash \left\{ i\right\} }\left\Vert \partial
_{x_{h}x_{j}}^{2}\Delta^g_{i,j}\right\Vert _{L^{\infty }}\Bigg )  \notag \\
&&+\frac{\left( L_{y}^{b,\sigma }\right) ^{2}}{N^{2}}\sum_{\ell \in
I_{N}\backslash \left\{ j\right\} ,h\in I_{N}\backslash \left\{ i\right\}
}\left\Vert \partial _{x_{h}x_{\ell }}^{2}\Delta^g_{i,j}\right\Vert
_{L^{\infty }}\Bigg ]\left\Vert u_{i}^{\prime }\right\Vert _{\mathcal{H}%
^{2}\left( \mathbb{R}\right) }\left\Vert u_{j}^{\prime \prime }\right\Vert _{%
\mathcal{H}^{2}\left( \mathbb{R}\right) }.  \label{m2}
\end{eqnarray}
Next we focus on the remaining two terms in (\ref{newv1}). We will
derive the estimation of solutions to BSDE (\ref{adjn1}) as follows. From
Lemma \ref{bs2}, by letting
\begin{eqnarray*}
\xi & \!\!\!=\!\!\! & \partial _{x}\Delta _{i,j}^{g}\left( \mathbf{X}_{T}^{\mathbf{u}%
}\right) ,\text{ }\forall i,j\in I_{N}\text{ with }i<j, \\
\mathsf{g}_{t}\left( P_{t}^{i,j},\left\{ Q_{t,k}^{i,j}\right\} _{k=1,\ldots
,N}\right) & \!\!\!=\!\!\! & \mathbf{B}\left( t,\mathbf{X}_{t}^{\mathbf{u}},\mathbf{u}%
_{t}\right) ^{\top }P_{t}^{i,j}+\sum_{k=1}^{N}\mathbf{\Pi }^{k}\left( t,%
\mathbf{X}_{t}^{\mathbf{u}},\mathbf{u}_{t}\right) ^{\top }Q_{t,k}^{i,j} \\
&&+\left( \partial _{x}\Delta _{i,j}^{f}\right) \left( t,\mathbf{X}_{t}^{%
\mathbf{u}},\mathbf{u}_{t}\right), 
\end{eqnarray*}
we derive, 
\begin{eqnarray*}
\mathbb{E}\left[ \sup_{0\leq t\leq T}\left\vert P_{t}^{i,j}\right\vert
^{2}+\int_{0}^{T}\left\vert Q_{t}^{i,j}\right\vert ^{2}\mathrm{d}s\right] 
& \!\!\!\leq\!\!\! & C\mathbb{E}\Bigg [\left\vert \partial _{x}\Delta _{i,j}^{g}\left(
\mathbf{X}_{T}^{\mathbf{u}}\right) \right\vert ^{2}+\left(
\int_{0}^{T}\left\vert \left( \partial _{x}\Delta _{i,j}^{f}\right) \left(t,\mathbf{X}_{t}^{\mathbf{u}},\mathbf{u}_{t}\right) \right\vert \mathrm{d}t\right) ^{2}\Bigg] \\
& \!\!\!\leq\!\!\! & C_{1}\mathbb{E}\Bigg [\left\vert \partial _{x}\Delta _{i,j}^{g}\left(
\mathbf{X}_{T}^{\mathbf{u}}\right) \right\vert ^{2}+T\int_{0}^{T}\left\vert
\left( \partial _{x}\Delta _{i,j}^{f}\right) \left( t,\mathbf{X}_{t}^{\mathbf{u}},\mathbf{u}_{t}\right) \right\vert ^{2}\mathrm{d}t\Bigg],
\end{eqnarray*}
where $C_{1}$ depends on the time $T$, $\max \left\{ \left\Vert \mathbf{B}%
_{0,x,y}\left( t,\mathbf{X}_{t}^{\mathbf{u}},\mathbf{u}_{t}\right)
\right\Vert _{2},\left\Vert \mathbf{\Pi }_{0,x,y}^{1}\left( t,\mathbf{X}%
_{t}^{\mathbf{u}},\mathbf{u}_{t}\right) \right\Vert _{2},\ldots ,\left\Vert
\mathbf{\Pi }_{0,x,y}^{N}\left( t,\mathbf{X}_{t}^{\mathbf{u}},\mathbf{u}%
_{t}\right) \right\Vert _{2}\right\} $ with
\begin{eqnarray*}
\left\Vert \mathbf{B}_{0,x,y}\left( t,\mathbf{X}_{t}^{\mathbf{u}},\mathbf{u}%
_{t}\right) \right\Vert _{2} &\leq &\left\Vert \mathbf{\bar{B}}_{x}\left( t,%
\mathbf{X}_{t}^{\mathbf{u}},\mathbf{u}_{t}\right) +\mathbf{B}_{0,y}\left( t,%
\mathbf{X}_{t}^{\mathbf{u}},\mathbf{u}_{t}\right) \right\Vert _{2} \\
&\leq &\left\Vert \mathbf{\bar{B}}_{x}\left( t,\mathbf{X}_{t}^{\mathbf{u}},%
\mathbf{u}_{t}\right) \right\Vert _{2}+\left\Vert \mathbf{B}_{0,y}\left( t,%
\mathbf{X}_{t}^{\mathbf{u}},\mathbf{u}_{t}\right) \right\Vert _{2} \\
&\leq &L^{b}+L_{y}^{b}\left( \frac{2}{N}-\frac{1}{N^{2}}\right),
\end{eqnarray*}
and
\begin{eqnarray*}
\left\Vert \mathbf{\Pi }_{0,x,y}^{j}\left( t,\mathbf{X}_{t}^{\mathbf{u}},%
\mathbf{u}_{t}\right) \right\Vert _{2} &\leq &\left\Vert \mathbf{\bar{\Pi}}%
_{x}^{j}\left( t,\mathbf{X}_{t}^{\mathbf{u}},\mathbf{u}_{t}\right) +\mathbf{%
\Pi }_{0,y}^{j}\left( t,\mathbf{X}_{t}^{\mathbf{u}},\mathbf{u}_{t}\right)
\right\Vert _{2} \\
&\leq &\left\Vert \mathbf{\bar{\Pi}}_{x}^{j}\left( t,\mathbf{X}_{t}^{\mathbf{%
u}},\mathbf{u}_{t}\right) \right\Vert _{2}+\left\Vert \mathbf{\Pi }%
_{0,y}^{j}\left( t,\mathbf{X}_{t}^{\mathbf{u}},\mathbf{u}_{t}\right)
\right\Vert _{2} \\
&\leq &L^{\sigma }+\frac{L_{y}^{\sigma }}{\sqrt{N}}.
\end{eqnarray*}%
The fundamental theorem of calculus implies that for all $\left(
t,x,u\right) =\left( t,\left( x_{\ell }\right) _{\ell =1}^{N},\left( u_{\ell
}\right) _{\ell =1}^{N}\right) \in \left[ 0,T\right] \times \mathbb{R}%
^{n}\times \mathbb{R}^{n}$ and $i,j\in I_{N}$,
\begin{eqnarray*}
\left\vert \left( \partial _{x_{\ell }}\Delta _{i,j}^{f}\right) \left(
t,x,u\right) \right\vert ^{2} &\leq &3\Bigg [\left\vert \left( \partial
_{x_{\ell }}\Delta _{i,j}^{f}\right) \left( t,0,0\right) \right\vert^{2} 
+\sum_{k=1}^{N}\left( \left\Vert \partial _{x_{_{\ell }}x_{k}}\Delta
_{i,j}^{f}\right\Vert _{L^{\infty }}^{2}\left\vert x_{k}\right\vert
^{2}+\left\Vert \partial _{x_{\ell }u_{k}}\Delta _{i,j}^{f}\right\Vert
_{L^{\infty }}^{2}\left\vert u_{k}\right\vert ^{2}\right) \Bigg ], \\
\left\vert \partial _{x_{\ell }}\Delta _{i,j}^{g}\left( x\right) \right\vert
^{2} &\leq &\left\vert \left( \partial _{x_{\ell }}\Delta _{i,j}^{g}\right)
\left( 0\right) \right\vert ^{2}+\sum_{k=1}^{N}\left\Vert \partial
_{x_{_{\ell }}x_{k}}\Delta _{i,j}^{g}\right\Vert _{L^{\infty
}}^{2}\left\vert x_{k}\right\vert ^{2}.
\end{eqnarray*}
Thus, from Lemma \ref{bs2}
\begin{eqnarray*}
&&\mathbb{E}\left[ \sup_{0\leq t\leq T}\left\vert P_{t}^{i,j}\right\vert
^{2}+\int_{0}^{T}\left\vert Q_{t}^{i,j}\right\vert ^{2}\mathrm{d}s\right] \\
& \leq\!\!\! & C_{1}\mathbb{E}\Bigg [\left\vert \partial _{x}\Delta _{i,j}^{g}\left(
\mathbf{X}_{T}^{\mathbf{u}}\right) \right\vert ^{2}+\int_{0}^{T}\left\vert
\left( \partial _{x}\Delta _{i,j}^{f}\right) \left( t,\mathbf{X}_{t}^{%
\mathbf{u}},\mathbf{u}_{t}\right) \right\vert ^{2}\mathrm{d}t\Bigg ] \\
& \leq\!\!\! & C_{1}\mathbb{E}\Bigg [\sum_{\ell \in I_{N}}\left\vert \left( \partial
_{x_{\ell }}\Delta _{i,j}^{g}\right) \left( 0\right) \right\vert
^{2}+\sum_{\ell ,k\in I_{N}}\left\Vert \partial _{x_{\ell }x_{k}}\Delta
_{i,j}^{g}\right\Vert _{L^{\infty }}^{2} \\
&&+3\int_{0}^{T}\left( \sum_{\ell }\left\vert \left( \partial _{x_{\ell
}}\Delta _{i,j}^{f}\right) \left( t,0,0\right) \right\vert ^{2}+\sum_{\ell
,k\in I_{N}}\left\Vert \partial _{x_{\ell }x_{k}}\Delta
_{i,j}^{f}\right\Vert _{L^{\infty }}^{2}+\sum_{\ell ,k\in I_{N}}\left\Vert
\partial _{x_{\ell }u_{k}}\Delta _{i,j}^{f}\right\Vert _{L^{\infty
}}^{2}\right) \mathrm{d}t\Bigg].
\end{eqnarray*}
Now consider the term
\begin{equation*}
\mathbb{E}\left[ \int_{0}^{T}\left\langle P_{t}^{i,j},\mathfrak{F}%
_{t}^{i,j}+\partial _{yy}^{2}\bar{b}\left( t,\mathbf{X}_{t},\mathbf{u}%
_{t}\right) \star \mathcal{Y}_{t}^{i,j}+\Gamma _{t}^{i,j}\left( \bar{b}%
\right) \right\rangle \mathrm{d}t\right] .
\end{equation*}
Recall
\begin{eqnarray*}
\partial _{yy}^{2}\bar{b}\left( t,\mathbf{X}_{t},\mathbf{u}_{t}\right) \star\mathcal{Y}_{t}^{i,j} 
&=&\left[
\begin{array}{c}
\text{tr}\left[ \partial _{yy}^{2}b_{1}\left( t,\cdot \right) \mathcal{Y}%
_{t}^{i,j}\right] \\
\vdots \\
\text{tr}\left[ \partial _{yy}^{2}b_{i}\left( t,\cdot \right) \mathcal{Y}%
_{t}^{i,j}\right] \\
\vdots \\
\text{tr}\left[ \partial _{yy}^{2}b_{N}\left( t,\cdot \right) \mathcal{Y}%
_{t}^{i,j}\right]%
\end{array}%
\right]_{N\times 1} =\left[
\begin{array}{c}
\text{tr}\left[ \partial _{yy}^{2}b_{1}\left( t,\cdot \right) \mathbf{Y}%
_{t}^{j}\left( \mathbf{Y}_{t}^{i}\right) ^{\top }\right] \\
\vdots \\
\text{tr}\left[ \partial _{yy}^{2}b_{i}\left( t,\cdot \right) \mathbf{Y}%
_{t}^{j}\left( \mathbf{Y}_{t}^{i}\right) ^{\top }\right] \\
\vdots \\
\text{tr}\left[ \partial _{yy}^{2}b_{N}\left( t,\cdot \right) \mathbf{Y}%
_{t}^{j}\left( \mathbf{Y}_{t}^{i}\right) ^{\top }\right]%
\end{array}\right]_{N\times 1} \\
&=&\left[
\begin{array}{c}
\left( \mathbf{Y}_{t}^{i}\right) ^{\top }\partial _{yy}^{2}b_{1}\left(
t,\cdot \right) \mathbf{Y}_{t}^{j} \\
\vdots \\
\left( \mathbf{Y}_{t}^{i}\right) ^{\top }\partial _{yy}^{2}b_{i}\left(
t,\cdot \right) \mathbf{Y}_{t}^{j} \\
\vdots \\
\left( \mathbf{Y}_{t}^{i}\right)^{\top }\partial _{yy}^{2}b_{N}\left(
t,\cdot \right) \mathbf{Y}_{t}^{j}
\end{array}\right]_{N\times 1}. 
\end{eqnarray*}
Then
\begin{eqnarray*}
\mathbb{E}\Bigg [\int_{0}^{T}\left\langle P_{t}^{i,j},\partial _{yy}^{2}%
\bar{b}\left( t,\mathbf{X}_{t},\mathbf{u}_{t}\right) \star \mathcal{Y}%
_{t}^{i,j}\right\rangle \mathrm{d}t\Bigg] 
& \!\!\!=\!\!\! & \mathbb{E}\Bigg [\int_{0}^{T}\sum_{k=1}^{N}P_{t,k}^{i,j}\left( \mathbf{Y}%
_{t}^{i}\right)^{\top }\partial _{yy}^{2}b_{k}\left( t,\cdot \right)
\mathbf{Y}_{t}^{j}\mathrm{d}t\Bigg ] \\
& \!\!\!=\!\!\! & \sum_{k=1}^{N}\mathbb{E}\Bigg [\int_{0}^{T}P_{t,k}^{i,j}\left( \mathbf{Y}%
_{t}^{i}\right) ^{\top }\partial _{yy}^{2}b_{k}\left( t,\cdot \right)
\mathbf{Y}_{t}^{j}\mathrm{d}t\Bigg].
\end{eqnarray*}
For simplicity, we first estimate the term $P_{t,k}^{i,j}\left( \mathbf{Y}%
_{t}^{i}\right) ^{\top }\partial _{yy}^{2}b_{k}\left( t,\cdot \right)
\mathbf{Y}_{t}^{j}$. Note 
\begin{equation*}
\mathbb{E}\left[ \int_{0}^{T}P_{t,k}^{i,j}\left( \mathbf{Y}_{t}^{i}\right)
^{\top }\partial _{yy}^{2}b_{k}\left( t,\cdot \right) \mathbf{Y}_{t}^{j}%
\mathrm{d}t\right] =\mathbb{E}\left[ \int_{0}^{T}P_{t,k}^{i,j}\sum_{h,\ell =1}^{N}\left(
\partial _{y_{h}y_{\ell }}^{2}b_{k}\left( t,\cdot \right) \right)
Y_{h,t}^{i}Y_{\ell ,t}^{j}\mathrm{d}t\right] .
\end{equation*}
For $i<j$,
\begin{eqnarray*}
I_{N}\times I_{N} & \!\!\!=\!\!\! & \left\{ \left( i,j\right) \right\} \cup \left\{ \left(
i,i\right) \right\} \cup \left\{ \left( j,j\right) \right\} \cup \left\{
\left( i,\ell \right) |\ell \in I_{N}\backslash \left\{ i,j\right\} \right\}
\\
&&\cup \left\{ \left( h,j\right) |h\in I_{N}\backslash \left\{ i,j\right\}
\right\} \cup \left\{ \left( h,\ell \right) |h\in I_{N}\backslash \left\{
i\right\} ,\ell \in I_{N}\backslash \left\{ j\right\} \right\}. 
\end{eqnarray*}
Then
\begin{eqnarray}
&&\mathbb{E}\Bigg [\int_{0}^{T}P_{t,k}^{i,j}\left( \mathbf{Y}_{t}^{i}\right)
^{\top }\partial _{yy}^{2}b_{k}\left( t,\cdot \right) \mathbf{Y}_{t}^{j}%
\Bigg ]  \notag \\
& =\!\!\! & \mathbb{E}\Bigg \{\int_{0}^{T}P_{t,k}^{i,j}\Bigg [\underset{I_{1}}{%
\underbrace{\left( \partial _{y_{i}y_{j}}^{2}b_{k}\left( t,\cdot \right)
\right) Y_{t,i}^{i}Y_{t,j}^{j}}}+\underset{I_{2}}{\underbrace{\left(
\partial _{y_{i}y_{i}}^{2}b_{k}\left( t,\cdot \right) \right)
Y_{t,i}^{i}Y_{t,i}^{j}}}  \notag \\
&&+\underset{I_{3}}{\underbrace{\left( \partial _{y_{j}y_{j}}^{2}b_{k}\left(
t,\cdot \right) \right) Y_{t,j}^{i}Y_{t,j}^{j}}}+\underset{I_{4}}{%
\underbrace{\sum_{\ell \in I_{N}\backslash \left\{ i,j\right\} }^{N}\left(
\partial _{y_{i}y_{\ell }}^{2}b_{k}\left( t,\cdot \right) \right)
Y_{t,i}^{i}Y_{t,\ell }^{j}}}  \notag \\
&&+\underset{I_{5}}{\underbrace{\sum_{h\in I_{N}\backslash \left\{
i,j\right\} }^{N}\left( \partial _{y_{h}y_{j}}^{2}b_{k}\left( t,\cdot
\right) \right) Y_{t,h}^{i}Y_{t,j}^{j}}}+\underset{I_{6}}{\underbrace{%
\sum_{h\in I_{N}\backslash \left\{ i\right\} }^{N}\sum_{\ell \in
I_{N}\backslash \left\{ j\right\} }^{N}\left( \partial _{y_{h}y_{\ell}}^{2}b_{k}\left(t,\cdot \right) \right) Y_{t,h}^{i}Y_{t,\ell }^{j}}}\Bigg]\Bigg \}\mathrm{d}t.  \label{est1}
\end{eqnarray}
For $I_{1}$,
\begin{eqnarray}
\mathbb{E}\left[ \int_{0}^{T}P_{t,k}^{i,j}\left( \partial
_{y_{i}y_{j}}^{2}b_{k}\left( t,\cdot \right) \right) Y_{t,i}^{i}Y_{t,j}^{j}%
\mathrm{d}t\right]  
& \!\!\!\leq\!\!\! & \frac{\sqrt{\Lambda _{1}}L_{y}^{b}}{N^{2}}\mathbb{E}\Bigg[%
\int_{0}^{T}Y_{t,i}^{i}Y_{t,j}^{j}\mathrm{d}t\Bigg ]  \notag \\
& \!\!\!\leq\!\!\! & \frac{L_{y}^{b}\sqrt{\Lambda _{1}}C}{N^{2}}\left\Vert u_{i}^{\prime
}\right\Vert _{\mathcal{H}^{2}\left( \mathbb{R}\right) }\left\Vert
u_{j}^{\prime \prime }\right\Vert _{\mathcal{H}^{2}\left( \mathbb{R}\right)}.  \label{est2}
\end{eqnarray}
For $I_{2}$,
\begin{eqnarray}
\mathbb{E}\left[ \int_{0}^{T}P_{t,k}^{i,j}\left( \partial
_{y_{i}y_{i}}^{2}b_{k}\left( t,\cdot \right) \right) Y_{t,i}^{i}Y_{t,i}^{j}%
\mathrm{d}t\right]  &\leq &\frac{L_{y}^{b}}{N}\mathbb{E}\left[
\int_{0}^{T}P_{t,k}^{i,j}Y_{t,i}^{i}Y_{t,i}^{j}\mathrm{d}t\right]   \notag \\
&\leq &\frac{L_{y}^{b}\sqrt{\Lambda _{1}}}{N}\mathbb{E}\Bigg [%
\int_{0}^{T}Y_{t,i}^{i}Y_{t,i}^{j}\mathrm{d}t\Bigg ]  \notag \\
&\leq &\frac{L_{y}^{b}\sqrt{\Lambda _{1}}}{N}\cdot \frac{CL_{y}^{b,\sigma }}{%
N}\left\Vert u_{i}^{\prime }\right\Vert _{\mathcal{H}^{2}\left( \mathbb{R}%
\right) }\left\Vert u_{j}^{\prime \prime }\right\Vert _{\mathcal{H}%
^{2}\left( \mathbb{R}\right) },  \label{est3}
\end{eqnarray}
For $I_{3}$,
\begin{equation*}
\mathbb{E}\left[ \int_{0}^{T}P_{t,k}^{i,j}\left( \partial
_{y_{j}y_{j}}^{2}b_{k}\left( t,\cdot \right) \right) Y_{t,j}^{i}Y_{t,j}^{j}%
\mathrm{d}t\right] \leq \frac{L_{y}^{b}\sqrt{\Lambda _{1}}}{N}\cdot \frac{%
CL_{y}^{b,\sigma }}{N}\left\Vert u_{i}^{\prime }\right\Vert _{\mathcal{H}%
^{2}\left( \mathbb{R}\right) }\left\Vert u_{j}^{\prime \prime }\right\Vert _{%
\mathcal{H}^{2}\left( \mathbb{R}\right) }.
\end{equation*}
For $I_{4}$,
\begin{eqnarray}
&&\mathbb{E}\left[ \int_{0}^{T}P_{t,k}^{i,j}\sum_{\ell \in I_{N}\backslash
\left\{ i,j\right\} }^{N}\left( \partial _{y_{i}y_{\ell }}^{2}b_{k}\left(
t,\cdot \right) \right) Y_{t,i}^{i}Y_{t,\ell }^{j}\mathrm{d}t\right]  \leq \frac{L_{y}^{b}\sqrt{\Lambda _{1}}}{2N^{2}}\mathbb{E}\left[
\int_{0}^{T}\sum_{\ell \in I_{N}\backslash \left\{ i,j\right\}
}^{N}Y_{t,i}^{i}Y_{t,\ell }^{j}\mathrm{d}t\right]   \notag \\
&\leq &\frac{L_{y}^{b}\sqrt{\Lambda _{1}}}{2N^{2}}\cdot \frac{%
CNL_{y}^{b,\sigma }}{N}\left\Vert u_{i}^{\prime }\right\Vert _{\mathcal{H}%
^{2}\left( \mathbb{R}\right) }\left\Vert u_{j}^{\prime \prime }\right\Vert _{%
\mathcal{H}^{2}\left( \mathbb{R}\right) }.  \label{est5}
\end{eqnarray}%
For $I_{5}$,
\begin{eqnarray}
&&\mathbb{E}\left[ \int_{0}^{T}P_{t,k}^{i,j}\sum_{h\in I_{N}\backslash
\left\{ i,j\right\} }^{N}\left( \partial _{y_{h}y_{j}}^{2}b_{k}\left(
t,\cdot \right) \right) Y_{t,h}^{i}Y_{t,j}^{j}\mathrm{d}t\right]  
\leq \frac{L_{y}^{b}\sqrt{\Lambda _{1}}}{2N^{2}}\mathbb{\cdot }\frac{%
NCL_{y}^{b,\sigma }}{N}\left\Vert u_{i}^{\prime }\right\Vert _{\mathcal{H}%
^{2}\left( \mathbb{R}\right) }\left\Vert u_{j}^{\prime \prime }\right\Vert _{%
\mathcal{H}^{2}\left( \mathbb{R}\right) }.  \label{est55}
\end{eqnarray}%
As for the last term in (\ref{est1}), we have
\begin{eqnarray}
&&\mathbb{E}\left[ \int_{0}^{T}P_{t,k}^{i,j}\sum_{h\in I_{N}\backslash
\left\{ i\right\} }^{N}\sum_{\ell \in I_{N}\backslash \left\{ j\right\}
}^{N}\left( \partial _{y_{h}y_{\ell }}^{2}b_{k}\left( t,\cdot \right)
\right) Y_{t,h}^{i}Y_{t,\ell }^{j}\mathrm{d}t\right]   \notag \\
& \leq\!\!\! &\mathbb{E}\Bigg [\int_{0}^{T}P_{t,k}^{i,j}\Bigg (\sum_{\ell \in
I_{N}\backslash \left\{ j\right\} }^{N}\left( \partial _{y_{j}y_{\ell
}}^{2}b_{k}\left( t,\cdot \right) \right) Y_{t,j}^{i}Y_{t,\ell }^{j}  \notag
\\
&&+\sum_{h\in I_{N}\backslash \left\{ i,j\right\} }^{N}\Bigg (\left(
\partial _{y_{h}y_{h}}^{2}b_{k}\left( t,\cdot \right) \right)
Y_{t,h}^{i}Y_{t,h}^{j}+\sum_{\ell \in I_{N}\backslash \left\{ j,h\right\}
}^{N}\left( \partial _{y_{h}y_{h}}^{2}b_{k}\left( t,\cdot \right) \right)
Y_{t,h}^{i}Y_{t,\ell }^{j}\Bigg )\Bigg )\mathrm{d}t\Bigg ]  \notag \\
& \leq\!\!\! & \frac{L_{y}^{b}\sqrt{\Lambda _{1}}}{N^{2}}\mathbb{E}\left[
\int_{0}^{T}\sum_{\ell \in I_{N}\backslash \left\{ j\right\}
}^{N}Y_{t,j}^{i}Y_{t,\ell }^{j}\mathrm{d}t\right] +\frac{L_{y}^{b}\sqrt{%
\Lambda _{1}}}{N}\mathbb{E}\left[ \int_{0}^{T}\sum_{h\in I_{N}\backslash
\left\{ i,j\right\} }^{N}Y_{t,h}^{i}Y_{t,h}^{j}\mathrm{d}t\right]   \notag \\
&&+\frac{L_{y}^{b}\sqrt{\Lambda _{1}}}{N^{2}}\mathbb{E}\Bigg [%
\int_{0}^{T}\sum_{h\in I_{N}\backslash \left\{ i,j\right\} }^{N}\sum_{\ell
\in I_{N}\backslash \left\{ j,h\right\} }^{N}Y_{t,h}^{i}Y_{t,\ell }^{j}%
\mathrm{d}t\Bigg]  \notag \\
& \leq\!\!\! & \Bigg[\frac{L_{y}^{b}}{N^{2}}\mathbb{\cdot }\frac{\sqrt{\Lambda _{1}}%
NC\left( L_{y}^{b,\sigma }\right) ^{2}}{N^{2}}+\frac{L_{y}^{b}}{N}\mathbb{%
\cdot }NC\frac{\sqrt{\Lambda _{1}}\left( L_{y}^{b,\sigma }\right) ^{2}}{N^{2}%
}  \notag \\
&&+\frac{L_{y}^{b}}{N^{2}}\mathbb{\cdot }N^{2}\frac{\sqrt{\Lambda _{1}}%
\left( L_{y}^{b,\sigma }\right) ^{2}}{N^{2}}\Bigg ]\left\Vert u_{i}^{\prime
}\right\Vert _{\mathcal{H}^{2}\left( \mathbb{R}\right) }\left\Vert
u_{j}^{\prime \prime }\right\Vert _{\mathcal{H}^{2}\left( \mathbb{R}\right)
}.  \label{est8}
\end{eqnarray}%
Combining (\ref{est2})-(\ref{est8}), we have
\begin{eqnarray*}
&&\mathbb{E}\Bigg [\int_{0}^{T}P_{t,k}^{i,j}\left( \mathbf{Y}_{t}^{i}\right)
^{\top }\partial _{yy}^{2}b_{k}\left( t,\cdot \right) \mathbf{Y}_{t}^{j}%
\Bigg ] \leq L_{y}^{b}C\sqrt{\Lambda _{1}}\Bigg [\frac{1}{N^{2}}+\frac{\left(
L_{y}^{b}+3\left( L_{y}^{\sigma }\right) ^{2}\right) }{N^{2}} \\
&&+\frac{\left( L_{y}^{b}+3\left( L_{y}^{\sigma }\right) ^{2}\right) ^{2}}{%
N^{3}}+\frac{\left( L_{y}^{b}+3\left( L_{y}^{\sigma }\right) ^{2}\right) ^{2}%
}{N^{2}}+\Bigg ]\left\Vert u_{i}^{\prime }\right\Vert _{\mathcal{H}%
^{2}\left( \mathbb{R}\right) }\left\Vert u_{j}^{\prime \prime }\right\Vert _{%
\mathcal{H}^{2}\left( \mathbb{R}\right)}.
\end{eqnarray*}
Therefore, 
\begin{eqnarray}
&&\mathbb{E}\Bigg [\int_{0}^{T}\left\langle P_{t}^{i,j},\partial _{yy}^{2}%
\bar{b}\left( t,\mathbf{X}_{t},\mathbf{u}_{t}\right) \star \mathcal{Y}%
_{t}^{i,j}\right\rangle \mathrm{d}t\Bigg ]  \notag \\
& \leq\!\!\! & L_{y}^{b}C\sqrt{\Lambda _{1}}\Bigg [\frac{1}{N}+\frac{\left(
L_{y}^{b}+3\left( L_{y}^{\sigma }\right) ^{2}\right) }{N}+\frac{\left(
L_{y}^{b}+3\left( L_{y}^{\sigma }\right) ^{2}\right) ^{2}}{N^{2}}  \notag \\
&&+\frac{\left( L_{y}^{b}+3\left( L_{y}^{\sigma }\right) ^{2}\right) ^{2}}{N}%
\Bigg ]\left\Vert u_{i}^{\prime }\right\Vert _{\mathcal{H}^{2}\left( \mathbb{%
R}\right) }\left\Vert u_{j}^{\prime \prime }\right\Vert _{\mathcal{H}%
^{2}\left( \mathbb{R}\right) }.  \label{p11}
\end{eqnarray}
Now
\begin{eqnarray*}
\mathbb{E}\left[ \int_{0}^{T}\left\langle P_{t}^{i,j},\mathfrak{F}_{t}^{i,j}+\Gamma _{t}^{i,j}\left( \bar{b}\right) \right\rangle \mathrm{d}t\right] 
& \!\!\!=\!\!\! & \mathbb{E}\Bigg \{\int_{0}^{T}\Big [\sum_{\ell =1}^{N}\underset{%
J_{1}}{\underbrace{P_{t,\ell }^{i,j}\Big (Y_{t,\ell }^{i}\partial
_{xx}^{2}b_{\ell }\left( t,\cdot \right) Y_{t,\ell }^{j}}} \\
&&+\underset{J_{2}}{\underbrace{\left( \mathbf{Y}_{t}^{i}\right) ^{\top
}\partial _{yx}^{2}b_{\ell }\left( t,\cdot \right) Y_{t,\ell }^{j}}}+%
\underset{J_{3}}{\underbrace{Y_{t,\ell }^{i}\partial _{xy}^{2}b_{\ell
}\left( t,\cdot \right) \mathbf{Y}_{t}^{j}}} \\
&&+\underset{J_{4}}{\underbrace{u_{t,\ell }^{\prime }\left( \partial
_{ux}^{2}b_{\ell }\left( t,\cdot \right) Y_{t,\ell }^{j}+\partial
_{uy}^{2}b_{\ell }\left( t,\cdot \right) \cdot \mathbf{Y}_{t}^{j}\right)
\delta _{\ell,i}}} \\
&&+\underset{J_{5}}{\underbrace{u_{t,\ell }^{\prime \prime }\left( \partial
_{ux}^{2}b_{\ell }\left( t,\cdot \right) Y_{t,\ell }^{i}+\partial
_{uy}^{2}b_{\ell }\left( t,\cdot \right) \cdot \mathbf{Y}_{t}^{i}\right)
\delta _{\ell,j}}}\Big )\Big ]\mathrm{d}t\Bigg \}.
\end{eqnarray*}
For $J_{1}$, 
\begin{eqnarray}
&&\mathbb{E}\Bigg [\int_{0}^{T}\sum_{\ell =1}^{N}P_{t,\ell }^{i,j}Y_{t,\ell}^{i}\partial _{xx}^{2}b_{\ell }\left( t,\cdot \right) Y_{t,\ell }^{j}%
\mathrm{d}t\Bigg ]  \notag \\
& \leq\!\!\! & L^{b}\sqrt{\Lambda _{1}}\mathbb{E}\Bigg [\int_{0}^{T}\Big (P_{t,\ell
}^{i,j}Y_{t,i}^{i}Y_{t,i}^{j}+P_{t,\ell
}^{i,j}Y_{t,j}^{i}Y_{t,j}^{j}+\sum_{\ell =I_{N}\backslash \left\{
i,j\right\} }^{N}P_{t,\ell }^{i,j}Y_{t,\ell }^{i}Y_{t,\ell }^{j}\Big )%
\mathrm{d}t\Bigg ]  \notag \\
& \leq\!\!\! & L^{b}\sqrt{\Lambda _{1}}\Bigg [\frac{C\left( L_{y}^{b}+3\left(
L_{y}^{\sigma }\right) ^{2}\right) }{N}+N\frac{C\left( L_{y}^{b}+3\left(
L_{y}^{\sigma }\right) ^{2}\right) ^{2}}{N^{2}}\Bigg ]\left\Vert
u_{i}^{\prime }\right\Vert _{\mathcal{H}^{2}\left( \mathbb{R}\right)
}\left\Vert u_{j}^{\prime \prime }\right\Vert _{\mathcal{H}^{2}\left(
\mathbb{R}\right) }  \notag \\
& \leq\!\!\! & L^{b}\sqrt{\Lambda _{1}}\Bigg [\frac{C\left( L_{y}^{b}+3\left(
L_{y}^{\sigma }\right) ^{2}\right) ^{2}}{N}\Bigg ]\left\Vert u_{i}^{\prime
}\right\Vert _{\mathcal{H}^{2}\left( \mathbb{R}\right) }\left\Vert
u_{j}^{\prime \prime }\right\Vert _{\mathcal{H}^{2}\left( \mathbb{R}\right)
}.  \label{j1}
\end{eqnarray}
For $J_{2}$, 
\begin{eqnarray}
&&\mathbb{E}\Bigg [\int_{0}^{T}\sum_{\ell =1}^{N}P_{t,\ell }^{i,j}\left(
\mathbf{Y}_{t}^{i}\right) ^{\top }\partial _{yx}^{2}b_{\ell }\left( t,\cdot
\right) Y_{t,\ell }^{j}\mathrm{d}t\Bigg ]  \notag \\
& \leq\!\!\! & \frac{L_{y}^{b}}{N}\mathbb{E}\Bigg [\int_{0}^{T}\sum_{\ell
=1}^{N}P_{t,\ell }^{i,j}\sum_{k=1}^{N}Y_{t,k}^{i}Y_{t,\ell }^{j}\mathrm{d}t%
\Bigg ]  \notag \\
& \leq\!\!\! & \frac{L_{y}^{b}}{N}\mathbb{E}\Bigg [\int_{0}^{T}\Big (%
P_{t,j}^{i,j}Y_{t,i}^{i}Y_{t,j}^{j}+\sum_{k\in I_{N}\backslash \left\{
j\right\} }^{N}P_{t,k}^{i,j}Y_{t,i}^{i}Y_{t,k}^{j}+\sum_{k\in
I_{N}\backslash \left\{ i\right\} }^{N}P_{t,j}^{i,j}Y_{t,k}^{i}Y_{t,j}^{j}
\notag \\
&&+\sum_{k\in I_{N}\backslash \left\{ j\right\} ,h\in I_{N}\backslash
\left\{ i\right\} }^{N}P_{t,k}^{i,j}Y_{t,h}^{i}Y_{t,k}^{j}\Big )\mathrm{d}t%
\Bigg ]  \notag \\
& \leq\!\!\! & \frac{L_{y}^{b}\sqrt{\Lambda _{1}}C}{N}\Bigg (1+\frac{2\left(
N-1\right) \left( L_{y}^{b}+3\left( L_{y}^{\sigma }\right) ^{2}\right) }{N}
\notag \\
&&+\frac{\left( N^{2}-2N+1\right) \left( L_{y}^{b}+3\left( L_{y}^{\sigma
}\right) ^{2}\right) ^{2}}{N^{2}}\Bigg )\left\Vert u_{i}^{\prime
}\right\Vert _{\mathcal{H}^{2}\left( \mathbb{R}\right) }\left\Vert
u_{j}^{\prime \prime }\right\Vert _{\mathcal{H}^{2}\left( \mathbb{R}\right) }
\notag \\
& \leq\!\!\! & \frac{L_{y}^{b}\sqrt{\Lambda _{1}}C}{N}\Bigg (1+\left(
L_{y}^{b}+3\left( L_{y}^{\sigma }\right) ^{2}\right) +\left(
L_{y}^{b}+3\left( L_{y}^{\sigma }\right) ^{2}\right) ^{2}\Bigg )\left\Vert
u_{i}^{\prime }\right\Vert _{\mathcal{H}^{2}\left( \mathbb{R}\right)
}\left\Vert u_{j}^{\prime \prime }\right\Vert _{\mathcal{H}^{2}\left(
\mathbb{R}\right)}.  \label{j2}
\end{eqnarray}
For $J_{3}$, 
\begin{eqnarray}
&&\mathbb{E}\Bigg [\int_{0}^{T}\sum_{\ell =1}^{N}P_{t,\ell }^{i,j}Y_{t,\ell
}^{i}\partial _{xy}^{2}b_{\ell }\left( t,\cdot \right) \mathbf{Y}_{t}^{j}%
\mathrm{d}t\Bigg]  \notag \\
&\leq\!\!\! &\frac{L_{y}^{b}\sqrt{\Lambda _{1}}C}{N}\Bigg (1+\left(
L_{y}^{b}+3\left( L_{y}^{\sigma }\right) ^{2}\right) +\left(
L_{y}^{b}+3\left( L_{y}^{\sigma }\right) ^{2}\right) ^{2}\Bigg )\left\Vert
u_{i}^{\prime }\right\Vert _{\mathcal{H}^{2}\left( \mathbb{R}\right)
}\left\Vert u_{j}^{\prime \prime }\right\Vert _{\mathcal{H}^{2}\left(
\mathbb{R}\right)}.  \label{j3}
\end{eqnarray}
For $J_{4}$, 
\begin{eqnarray}
&&\mathbb{E}\Bigg \{\int_{0}^{T}\sum_{\ell =1}^{N}P_{t,\ell }^{i,j}u_{\ell
,t}^{\prime }\left( \partial _{ux}^{2}b_{\ell }\left( t,\cdot \right)
Y_{t,\ell }^{j}+\partial _{uy}^{2}b_{\ell }\left( t,\cdot \right) \cdot
\mathbf{Y}_{t}^{j}\right) \delta _{\ell =i}\mathrm{d}t\Bigg \}  \notag \\
& =\!\!\! & \mathbb{E}\Bigg \{\int_{0}^{T}P_{t,i}^{i,j}\left( \partial
_{ux}^{2}b_{i}\left( t,\cdot \right) Y_{t,i}^{j}+\partial
_{uy}^{2}b_{i}\left( t,\cdot \right) \cdot \mathbf{Y}_{t}^{j}\right)
u_{t,i}^{\prime }\mathrm{d}t\Bigg \}  \notag \\
& \leq\!\!\! &L ^{b}\mathbb{E}\left[ \int_{0}^{T}P_{t,i}^{i,j}Y_{t,i}^{j}\mathrm{d}t%
\right] +\mathbb{E}\left[ \int_{0}^{T}P_{t,i}^{i,j}\partial
_{uy}^{2}b_{i}\left( t,\cdot \right) \cdot \mathbf{Y}_{t}^{j}u_{t,i}^{\prime
}\mathrm{d}t\right]   \notag \\
& \leq\!\!\! & L^{b}\sqrt{\Lambda _{1}}\frac{\left( L_{y}^{b}+3\left( L_{y}^{\sigma
}\right) ^{2}\right) }{N}  \notag \\
&&+\mathbb{E}\left[ \int_{0}^{T}P_{t,i}^{i,j}\Big (\partial
_{uy_{j}}^{2}b_{i}\left( t,\cdot \right) Y_{t,j}^{j}+\sum_{k\in
I_{N}\backslash \left\{ j\right\} }\partial _{uy_{k}}^{2}b_{i}\left( t,\cdot
\right) Y_{t,k}^{j}\Big )\mathrm{d}t\right]   \notag \\
&\leq &C\sqrt{\Lambda _{1}}\left( L_{y}^{b}+3\left( L_{y}^{\sigma }\right)
^{2}\right) \left( \frac{L^{b}+L_{y}^{b}}{N}\right) \left\Vert u_{i}^{\prime
}\right\Vert _{\mathcal{H}^{2}\left( \mathbb{R}\right) }\left\Vert
u_{j}^{\prime \prime }\right\Vert _{\mathcal{H}^{2}\left( \mathbb{R}\right)
}.  \label{j4}
\end{eqnarray}
For $J_{5}$, 
\begin{eqnarray}
&&\mathbb{E}\Bigg \{\int_{0}^{T}\sum_{\ell =1}^{N}P_{t,\ell }^{i,j}u_{t,\ell
}^{\prime \prime }\left( \partial _{ux}^{2}b_{\ell }\left( t,\cdot \right)
Y_{t,\ell }^{i}+\partial _{uy}^{2}b_{\ell }\left( t,\cdot \right) \cdot
\mathbf{Y}_{t}^{i}\right) \delta _{\ell,j}\mathrm{d}t\Bigg \}  \notag \\
& \leq\!\!\! & C\sqrt{\Lambda _{1}}\left( L_{y}^{b}+3\left( L_{y}^{\sigma }\right)
^{2}\right) \left( \frac{L^{b}+L_{y}^{b}}{N}\right) \left\Vert u_{i}^{\prime
}\right\Vert _{\mathcal{H}^{2}\left( \mathbb{R}\right) }\left\Vert
u_{j}^{\prime \prime }\right\Vert _{\mathcal{H}^{2}\left( \mathbb{R}\right)
}.  \label{j5}
\end{eqnarray}
From (\ref{j1})-(\ref{j5}), we can get
\begin{eqnarray}
\mathbb{E}\left[ \int_{0}^{T}\left\langle P_{t}^{i,j},\mathfrak{F}_{t}^{i,j}+\Gamma _{t}^{i,j}\left( \bar{b}\right) \right\rangle \mathrm{d}t\right]  
& \!\!\!\leq\!\!\! & C\sqrt{\Lambda _{1}}\Bigg \{\frac{L^{b}\left(L_{y}^{b}+3\left( L_{y}^{\sigma }\right) ^{2}\right) ^{2}}{N}  \notag \\
&&+\frac{2L_{y}^{b}}{N}\Big [1+\left( L_{y}^{b}+3\left( L_{y}^{\sigma}\right) ^{2}\right) +\left( L_{y}^{b}+3\left( L_{y}^{\sigma }\right)^{2}\right) ^{2}\Big ]  \notag \\
&&+2\left(L_{y}^{b}+3\left(L_{y}^{\sigma }\right)^{2}\right)\left(\frac{L^{b}+L_{y}^{b}}{N}\right)\Bigg\}\left\Vert u_{i}^{\prime}\right\Vert_{\mathcal{H}^{2}\left(\mathbb{R}\right)}\left\Vert u_{j}^{\prime \prime}\right\Vert _{\mathcal{H}^{2}\left( \mathbb{R}\right) }.  \label{p22}
\end{eqnarray}
Therefore, from (\ref{p11})-(\ref{p22}), it yields
\begin{eqnarray}
&&\mathbb{E}\left[ \int_{0}^{T}\left\langle P_{t}^{i,j},\mathfrak{F}%
_{t}^{i,j}+\partial _{yy}^{2}\bar{b}\left( t,\mathbf{X}_{t},\mathbf{u}%
_{t}\right) \star \mathcal{Y}_{t}^{i,j}+\Gamma _{t}^{i,j}\left( \bar{b}%
\right) \right\rangle \mathrm{d}t\right]   \notag \\
&\leq\!\!\! &C\sqrt{\Lambda _{1}}\Bigg \{\frac{L_{y}^{b}}{N}+\frac{\left(
L_{y}^{b}+3\left( L_{y}^{\sigma }\right) ^{2}\right) L_{y}^{b}}{N}  \notag \\
&&+\frac{L_{y}^{b}\left( L_{y}^{b}+3\left( L_{y}^{\sigma }\right)
^{2}\right) }{N}+\frac{L_{y}^{b}\left( L_{y}^{b}+3\left( L_{y}^{\sigma
}\right) ^{2}\right) ^{2}}{N^{2}}  \notag \\
&&+\frac{L_{y}^{b}\left( L_{y}^{b}+3\left( L_{y}^{\sigma }\right)
^{2}\right) ^{2}}{N}+\frac{L^{b}\left( L_{y}^{b}+3\left( L_{y}^{\sigma
}\right) ^{2}\right) ^{2}}{N}  \notag \\
&&+\frac{2L_{y}^{b}}{N}\Big [1+\left( L_{y}^{b}+3\left( L_{y}^{\sigma
}\right) ^{2}\right) +\left( L_{y}^{b}+3\left( L_{y}^{\sigma }\right)
^{2}\right) ^{2}\Big ]  \notag \\
&&+2\left( L_{y}^{b}+3\left( L_{y}^{\sigma }\right) ^{2}\right) \left( \frac{%
L^{b}+L_{y}^{b}}{N}\right) \Bigg \}\left\Vert u_{i}^{\prime }\right\Vert _{%
\mathcal{H}^{2}\left( \mathbb{R}\right) }\left\Vert u_{j}^{\prime \prime
}\right\Vert _{\mathcal{H}^{2}\left( \mathbb{R}\right) }.  \label{m3}
\end{eqnarray}
Next, set
\begin{eqnarray*}
\left( \mathfrak{G}_{t,k}^{i,j}+\bar{\Xi}_{t,k}^{i,j}\left( \bar{\sigma}%
\right) \right) _{k} & \!\!\!=\!\!\! & Y_{t,k}^{\mathbf{u},u_{i}^{\prime }}\partial
_{xx}^{2}\sigma _{k}\left( t,\cdot \right) Y_{t,k}^{\mathbf{u},u_{j}^{\prime
\prime }}+\left( \mathbf{Y}_{t}^{\mathbf{u},u_{i}^{\prime }}\right) ^{\top
}\partial _{yx}^{2}\sigma _{k}\left( t,\cdot \right) Y_{t,k}^{\mathbf{u}%
,u_{j}^{\prime \prime }} \\
&&+Y_{t,k}^{\mathbf{u},u_{i}^{\prime }}\partial _{xy}^{2}\sigma _{k}\left(
t,\cdot \right) \mathbf{Y}_{t}^{\mathbf{u},u_{j}^{\prime \prime }} \\
&&+u_{t,i}^{\prime }\left( \partial _{ux}^{2}\sigma _{k}\left( t,\cdot
\right) \cdot Y_{t,k}^{\mathbf{u},u_{j}^{\prime \prime }}+\partial
_{uy}^{2}\sigma _{k}\left( t,\cdot \right) \cdot \mathbf{Y}_{t}^{\mathbf{u}%
,u_{j}^{\prime \prime }}\right) \delta _{i,k} \\
&&+u_{t,j}^{\prime \prime }\left( \partial _{xu}^{2}\sigma _{k}\left(
t,\cdot \right) \cdot Y_{t,k}^{\mathbf{u},u_{i}^{\prime }}+\partial
_{yu}^{2}\sigma _{k}\left( t,\cdot \right) \cdot \mathbf{Y}_{t}^{\mathbf{u}%
,u_{i}^{\prime }}\right) \delta _{j,k}.
\end{eqnarray*}%
Then
\begin{equation*}
\left\langle Q_{t,k}^{i,j},\mathfrak{G}_{t,k}^{i,j}+\bar{\Xi}%
_{t,k}^{i,j}\left( \bar{\sigma}\right) \right\rangle
=\sum_{k=1}^{N}Q_{t,k,k}^{i,j}\left( \mathfrak{G}_{t,k}^{i,j}+\bar{\Xi}%
_{t,k}^{i,j}\left( \bar{\sigma}\right) \right) _{k}.
\end{equation*}%
Repeating similar analysis, we have
\begin{eqnarray}
&&\mathbb{E}\Bigg [\int_{0}^{T}\sum_{k=1}^{N}\left\langle Q_{t,k}^{i,j},%
\mathfrak{G}_{t,k}^{i,j}+\bar{\Xi}_{t,k}^{i,j}\left( \bar{\sigma}\right)
\right\rangle \mathrm{d}t\Bigg ]  \notag \\
&=\!\!\!&\sum_{k=1}^{N}\mathbb{E}\Bigg [\int_{0}^{T}\left\langle Q_{t,k}^{i,j},%
\mathfrak{G}_{t,k}^{i,j}+\bar{\Xi}_{t,k}^{i,j}\left( \bar{\sigma}\right)
\right\rangle \mathrm{d}t\Bigg ]  \notag \\
&\leq\!\!\! &\sum_{k=1}^{N}\left[ \mathbb{E}\left[ \int_{0}^{T}\left(
Q_{t,k,k}^{i,j}\right) ^{2}\mathrm{d}t\right] ^{\frac{1}{2}}\cdot \mathbb{E}%
\left[ \int_{0}^{T}\left( \mathfrak{G}_{t,k}^{i,j}+\bar{\Xi}%
_{t,k}^{i,j}\left( \bar{\sigma}_{t}\right) \right) _{k}^{2}\mathrm{d}t\right]
^{\frac{1}{2}}\right]   \notag \\
&=\!\!\!&\sqrt{\Lambda _{1}}\sum_{k=1}^{N}\left[ \mathbb{E}\left[
\int_{0}^{T}\left( \mathfrak{G}_{t,k}^{i,j}+\bar{\Xi}_{t,k}^{i,j}\left( \bar{%
\sigma}_{t}\right) \right) _{k}^{2}\mathrm{d}t\right] ^{\frac{1}{2}}\right]
\notag \\
&\leq\!\!\! &C\sqrt{\Lambda _{1}}\Bigg \{\frac{L^{\sigma }\left( L_{y}^{b}+3\left(
L_{y}^{\sigma }\right) ^{2}\right) ^{2}}{N}+\frac{2L_{y}^{\sigma }}{N}\Big [%
1+\left( L_{y}^{b}+3\left( L_{y}^{\sigma }\right) ^{2}\right) +\left(
L_{y}^{b}+3\left( L_{y}^{\sigma }\right) ^{2}\right) ^{2}\Big ]  \notag \\
&&+2\left( L_{y}^{b}+3\left( L_{y}^{\sigma }\right) ^{2}\right) \left( \frac{%
L^{\sigma }+L_{y}^{\sigma }}{N}\right) \Bigg ]\Bigg \}\left\Vert
u_{i}^{\prime }\right\Vert _{\mathcal{H}^{2}\left( \mathbb{R}\right)
}\left\Vert u_{j}^{\prime \prime }\right\Vert _{\mathcal{H}^{2}\left(
\mathbb{R}\right)},  \label{m44}
\end{eqnarray}
where the third line follows by H\"{o}lder inequality and the fifth line
follows by noting $\sqrt{a_{1}+\cdots +a_{N}}\leq \sqrt{a_{1}}+\cdots +\sqrt{%
a_{N}}$ for any $a_{1},\cdots ,a_{N}\geq 0.$

Finally, 
\begin{eqnarray}
&&\mathbb{E}\left[ \int_{0}^{T}\sum_{k=1}^{N}\left\langle
Q_{t,k}^{i,j},\partial _{yy}^{2}\bar{\sigma}_{t,k}\left( t,\mathbf{X}_{t},%
\mathbf{u}_{t}\right) \star \mathcal{Y}_{t}^{i,j}\right\rangle \mathrm{d}t%
\right]   \notag \\
&\!\!\!=\!\!\!&\sum_{k=1}^{N}\mathbb{E}\left[ \int_{0}^{T}Q_{t,k,k}^{i,j}\left( \mathbf{Y%
}_{t}^{i}\right) ^{\top }\partial _{yy}^{2}\bar{\sigma}_{t,k}\left( t,%
\mathbf{X}_{t},\mathbf{u}_{t}\right) \mathbf{Y}_{t}^{j}\mathrm{d}t\right]
\notag \\
&\!\!\!\leq\!\!\! &\sum_{k=1}^{N}\left[ \mathbb{E}\left[ \int_{0}^{T}\left(
Q_{t,k,k}^{i,j}\right) ^{2}\mathrm{d}t\right] ^{\frac{1}{2}}\cdot \mathbb{E}%
\left[ \int_{0}^{T}\left( \left( \mathbf{Y}_{t}^{i}\right) ^{\top }\partial
_{yy}^{2}\bar{\sigma}_{t,i}\left( X_{t,i},\mathbf{X}_{t}\right) \mathbf{Y}%
_{t}^{j}\right) _{k}^{2}\mathrm{d}t\right] ^{\frac{1}{2}}\right]   \notag \\
&\!\!\!\leq\!\!\! &L_{y}^{\sigma }C\sqrt{\Lambda _{1}}\Bigg [\frac{1\!+\!\left(
L_{y}^{b}\!+\!3\left( L_{y}^{\sigma }\right) ^{2}\right) }{N}\!+\!\frac{\left(
L_{y}^{b}\!+\!3\left( L_{y}^{\sigma }\right) ^{2}\right) ^{2}}{N}\!+\!\frac{\left(
L_{y}^{b}\!+\!3\left( L_{y}^{\sigma }\right) ^{2}\right) ^{2}}{N^{2}}\Bigg ]%
\left\Vert u_{i}^{\prime }\right\Vert _{\mathcal{H}^{2}\left( \mathbb{R}%
\right) }\left\Vert u_{j}^{\prime \prime }\right\Vert _{\mathcal{H}%
^{2}\left( \mathbb{R}\right)}.  \label{m5}
\end{eqnarray}
Thus combining (\ref{m1})-(\ref{m5}), we have the desired results. \hfill $\Box$

\paragraph{Proof of Corollary \ref{cor2}.}
If $i=j,$ this case is trivial. So we consider the case $i\neq j$ and deal with the estimate in the following way: 

\begin{eqnarray*}
\begin{aligned}
\tilde{C}^{i,j} = & \, C\Bigg \{\left\Vert \partial _{x_{i}x_{j}}^{2}\Delta^f_{i,j}\right\Vert _{L^{\infty }}+\left\Vert \partial
_{x_{i}u_{j}}^{2}\Delta^f_{i,j}\right\Vert _{L^{\infty }}+\left\Vert
\partial _{u_{i}x_{j}}^{2}\Delta^f_{i,j}\right\Vert _{L^{\infty
}}+\left\Vert \partial _{u_{i}u_{j}}^{2}\Delta^f_{i,j}\right\Vert
_{L^{\infty }}+\left\Vert \partial _{x_{i}x_{j}}^{2}\Delta^g_{i,j}\right\Vert _{L^{\infty }} \\
&+\frac{L_{y}^{b,\sigma }}{N}\Bigg [\sum_{\ell \in I_{N}\backslash \left\{
j\right\} }\left( \left\Vert \partial _{x_{i}x_{\ell }}^{2}\Delta^f_{i,j}\right\Vert _{L^{\infty }}+\left\Vert \partial _{u_{i}x_{\ell
}}^{2}\Delta^f_{i,j}\right\Vert _{L^{\infty }}\right)  \\
&+\sum_{h\in I_{N}\backslash \left\{ i\right\} }\left( \left\Vert \partial
_{x_{h}x_{j}}^{2}\Delta^f_{i,j}\right\Vert _{L^{\infty }}+\left\Vert
\partial _{x_{h}u_{j}}^{2}\Delta^f_{i,j}\right\Vert _{L^{\infty }}\right) %
\Bigg ]+\frac{\left( L_{y}^{b,\sigma }\right) ^{2}}{N^{2}}\Bigg (\sum_{\ell
\in I_{N}\backslash \left\{ j\right\} ,h\in I_{N}\backslash \left\{
i\right\} }\left\Vert \partial _{x_{h}x_{\ell }}^{2}\Delta^f_{i,j}\right\Vert _{L^{\infty }}\Bigg ) \\
&+\frac{L_{y}^{b,\sigma }}{N}\Bigg (\sum_{\ell \in I_{N}\backslash \left\{
j\right\} }\left\Vert \partial _{x_{i}x_{\ell }}^{2}\Delta^g_{i,j}\right\Vert _{L^{\infty }}+\sum_{h\in I_{N}\backslash \left\{
i\right\} }\left\Vert \partial _{x_{h}x_{j}}^{2}\Delta^g_{i,j}\right\Vert
_{L^{\infty }}\Bigg ) \\
&+\frac{\left( L_{y}^{b,\sigma }\right) ^{2}}{N^{2}}\Bigg (\sum_{\ell \in
I_{N}\backslash \left\{ j\right\} ,h\in I_{N}\backslash \left\{ i\right\}
}\left\Vert \partial _{x_{h}x_{\ell }}^{2}\Delta^g_{i,j}\right\Vert
_{L^{\infty }}\Bigg )+\sqrt{\Lambda _{1}}\Bigg(\frac{C^{1,b,\sigma }}{N}+%
\frac{C^{2,b,\sigma }}{N^{2}}\Bigg)\Bigg\}, 
\end{aligned}
\end{eqnarray*}
where
\begin{eqnarray}
C^{1,b,\sigma } & \!\!\!=\!\!\! & L_{y}^{b}+L_{y}^{b,\sigma }L_{y}^{b}+L_{y}^{b}\left(
L_{y}^{b,\sigma }\right) ^{2}+L^{b}\left( L_{y}^{b,\sigma }\right)
^{2}+2L_{y}^{b}\Big (1+L_{y}^{b,\sigma }  \notag \\
&&+\left( L_{y}^{b,\sigma }\right) ^{2}\Big )+2L^{b}L_{y}^{b,\sigma}+L_{y}^{b}+L_{y}^{\sigma }+L_{y}^{b,\sigma }L_{y}^{\sigma }+L_{y}^{\sigma}\left( L_{y}^{b,\sigma }\right) ^{2}  \notag \\
&&+L^{\sigma }\left( L_{y}^{b,\sigma }\right) ^{2}+2L_{y}^{\sigma }\Big(1+L_{y}^{b,\sigma}+\left( L_{y}^{b,\sigma }\right) ^{2}\Big )+2L^{\sigma}L_{y}^{b,\sigma }+L_{y}^{\sigma },  \label{mc1} \\
C^{2,b,\sigma } & \!\!\!=\!\!\! & \left( L_{y}^{b}+L_{y}^{\sigma }\right) \left(
L_{y}^{b,\sigma }\right)^{2}. \label{mc2}
\end{eqnarray}
Then
\begin{eqnarray*} 
\begin{aligned}
\tilde{C}^{i,j} \leq & \, \frac{C\tilde{L}}{N^{2\beta }}+\frac{4L_{y}^{b,\sigma
}\tilde{L}}{N}\cdot \frac{1}{N^{\min \left\{ \beta ,2\beta -1\right\} }}+%
\frac{\left( L_{y}^{b,\sigma }\right) ^{2}\tilde{L}}{N^{2}}\left( N-1\right)
^{2}\frac{\tilde{L}}{N^{2\beta }} \\
& +\frac{2L_{y}^{b,\sigma }\tilde{L}}{N}\left( N-1\right) \frac{\tilde{L}}{%
N^{2\beta }}+\frac{\left( L_{y}^{b,\sigma }\right) ^{2}\tilde{L}}{N^{2}}%
\left( N-1\right) ^{2}\frac{\tilde{L}}{N^{2\beta }} \\
& +C\sqrt{\Lambda _{1}}\left( \frac{C^{1,b,\sigma }}{N}+\frac{C^{2,b,\sigma }%
}{N^{2}}\right),
\end{aligned}
\end{eqnarray*}
while
\begin{eqnarray*}
\begin{aligned}
\Lambda _{1} = & \, C_{1}\mathbb{E}\Bigg [\sum_{\ell \in I_{N}}\left\vert \left(
\partial _{x_{\ell }}\Delta _{i,j}^{g}\right) \left( 0\right) \right\vert
^{2}+\sum_{\ell ,k\in I_{N}}\left\Vert \partial _{x_{\ell }x_{k}}\Delta
_{i,j}^{g}\right\Vert _{L^{\infty }}^{2} \\
& +3\int_{0}^{T}\Bigg (\sum_{\ell }\left\vert \left( \partial _{x_{\ell
}}\Delta _{i,j}^{f}\right) \left( t,0,0\right) \right\vert ^{2}+\sum_{\ell
,k\in I_{N}}\left\Vert \partial _{x_{\ell }x_{k}}\Delta
_{i,j}^{f}\right\Vert _{L^{\infty }}^{2} \\
& +\sum_{\ell ,k\in I_{N}}\left\Vert \partial _{x_{\ell }u_{k}}\Delta
_{i,j}^{f}\right\Vert _{L^{\infty }}^{2}\Bigg )\mathrm{d}t\Bigg ] \\
\leq & \, C_{1}\Bigg [\frac{2}{N^{\beta -1}}+3\left( \frac{1}{N^{\beta }}+\frac{%
1}{N^{2\beta -1}}\right)\Bigg].
\end{aligned}
\end{eqnarray*}
Thus, 
\begin{eqnarray*}
\begin{aligned}
\tilde{C}^{i,j} \leq & \, \frac{C\tilde{L}}{N^{2\beta }}+\frac{4\left(
L_{y}^{b,\sigma }\right) \tilde{L}}{N}\cdot \frac{1}{N^{\min \left\{ \beta
,2\beta -1\right\} }}+\frac{\left( L_{y}^{b,\sigma }\right) ^{2}\tilde{L}}{%
N^{2}}\left( N-1\right) ^{2}\frac{\tilde{L}}{N^{2\beta }} \\
& +\frac{2L_{y}^{b,\sigma }\tilde{L}}{N}\left( N-1\right) \frac{\tilde{L}}{%
N^{2\beta }}+\frac{\left( L_{y}^{b,\sigma }\right) ^{2}\tilde{L}}{N^{2}}%
\left( N-1\right) ^{2}\frac{\tilde{L}}{N^{2\beta }} \\
& +C\sqrt{C_{1}}\max \left\{ C^{1,b,\sigma },C^{2,b,\sigma }\right\} \frac{1%
}{N^{\frac{\beta +1}{2}}}. 
\end{aligned}
\end{eqnarray*}
The proof is complete.
\hfill $\Box $

\paragraph{Concluding remarks.}
\label{sect7}
In this paper, we adopt the BSDE approach to analyze a general class of $\alpha$-potential stochastic differential games. By employing the duality principle, we provide a precise estimation of $\alpha$. 

In contrast to the previous work of \cite{glmsw03}, our BSDE analysis fully utilizes the structure of the open-loop games and has managed to avoid estimating the second-order sensitivity process. 
For the more challenging closed-loop games where players can dynamically adjust their control strategies based on the real-time state of the system,  this BSDE technique  seems promising of  obtaining the precise estimate of $\alpha$ for  non-zero-sum stochastic differential games \cite{GLZ15clo}.

\section{Appendix: Proofs of Technical Lemmas}

\label{APP}


\paragraph{Proof of Lemma \ref{bs2}.}
For any $\left( y,z\right) \in \mathbb{R}^{m}\times \mathbb{R}^{m\times d},$
let
\begin{equation*}
F_{t}\left( y,z\right) =A_{t}y+\sum_{j=1}^{d}B_{t}^{j}z^{j}+f_{t}.
\end{equation*}
Then
\begin{eqnarray*}
\left\vert F_{t}\left( y,z\right) \right\vert 
= \left\vert A_{t}y+\sum_{j=1}^{d}B_{t}^{j}z^{j}\right\vert +\left\vert
f_{t}\right\vert 
\leq C_{1}\left( \left\vert y\right\vert +\sum_{j=1}^{d}\left\vert
z^{j}\right\vert \right) +\left\vert f_{t}\right\vert ,
\end{eqnarray*}
where
$C_{1}=\max \left\{ \left\vert A_{t}\right\vert ,\left\vert
B_{t}^{1}\right\vert ,\ldots ,\left\vert B_{t}^{d}\right\vert \right\}.$
Note that
\begin{eqnarray*} 
\begin{aligned}
\left\vert y_{t}\right\vert 
& = \left\vert \xi +\int_{t}^{T}\left(
A_{s}y_{s}+\sum_{j=1}^{d}B_{s}^{j}z_{s}^{j}+f_{s}\right) \mathrm{d}%
s-\int_{t}^{T}z_{s}\mathrm{d}W_{s}\right\vert \\
& \leq \left\vert \xi \right\vert +\int_{t}^{T}\left( C_{1}\left\vert
y_{s}\right\vert +C_{1}\sum_{j=1}^{d}\left\vert z_{s}^{j}\right\vert
+\left\vert f_{s}\right\vert \right) \mathrm{d}s+\left\vert \int_{t}^{T}z_{s}%
\mathrm{d}W_{s}\right\vert.
\end{aligned}
\end{eqnarray*}
Then letting $y_{t}^{\ast }=\displaystyle\sup_{0\leq s\leq t}\left\vert y_{t}\right\vert$, we have
\begin{equation*}
y_{T}^{\ast }\leq \left\vert \xi \right\vert +\int_{t}^{T}\left(
C_{1}\left\vert y_{s}\right\vert +C_{1}\sum_{j=1}^{d}\left\vert
z_{s}^{j}\right\vert +\left\vert f_{s}\right\vert \right) \mathrm{d}%
s+\sup_{0\leq t\leq T}\left\vert \int_{0}^{t}z_{s}\mathrm{d}W_{s}\right\vert.
\end{equation*}%
By Burkholder-Davis-Gundy inequality
\begin{eqnarray*}
\begin{aligned}
\mathbb{E}\left[ \left\vert y_{T}^{\ast }\right\vert ^{2}\right] 
& \leq  4\mathbb{E}\Bigg [\left\vert \xi \right\vert ^{2}+\left(
\int_{0}^{T}\left\vert f_{s}\right\vert \right) ^{2}+C_{1}^{2}\left(
\int_{0}^{T}\left(\left\vert y_{s}\right\vert +\sum_{j=1}^{d}\left\vert
z_{s}^{j}\right\vert \right) \mathrm{d}s\right) ^{2} 
+\int_{0}^{T}\sum_{j=1}^{d}\left\vert z_{s}^{j}\right\vert ^{2}\mathrm{d}s\Bigg] \\
& \leq 4\mathbb{E}\Bigg [\left\vert \xi \right\vert ^{2}+\left(
\int_{0}^{T}\left\vert f_{s}\right\vert \right) ^{2}+C_{1}^{2}\left(
d+1\right) \int_{0}^{T}\left( \left\vert y_{s}\right\vert
^{2}+\sum_{j=1}^{d}\left\vert z_{s}^{j}\right\vert ^{2}\right) \mathrm{d}s 
+\int_{0}^{T}\sum_{j=1}^{d}\left\vert z_{s}^{j}\right\vert ^{2}\mathrm{d}s\Bigg ] \\
& \leq 4\mathbb{E}\Bigg [\left\vert \xi \right\vert ^{2}+\left(
\int_{0}^{T}\left\vert f_{s}\right\vert \right) ^{2}+C_{1}^{2}\left(
d+1\right) \int_{0}^{T}\left\vert y_{s}\right\vert ^{2}\mathrm{d}s 
+\left( C_{1}^{2}\left( d+1\right) +1\right)
\int_{0}^{T}\sum_{j=1}^{d}\left\vert z_{s}^{j}\right\vert ^{2}\mathrm{d}s\Bigg], 
\end{aligned}
\end{eqnarray*}
which implies
\begin{equation}
\mathbb{E}\left[ \left\vert y_{T}^{\ast }\right\vert ^{2}\right] 
\leq 4\mathbb{E}\Bigg [\left\vert \xi \right\vert ^{2}+\left(
\int_{0}^{T}\left\vert f_{s}\right\vert \right) ^{2}+\left( C_{1}^{2}\left(
d+1\right) +1\right) \int_{0}^{T}\left( \left\vert y_{s}\right\vert
^{2}+\sum_{j=1}^{d}\left\vert z_{s}^{j}\right\vert ^{2}\right) \mathrm{d}s\Bigg].  \label{B0}
\end{equation}
Next, by It\^{o}'s formula  
\begin{eqnarray*}
\mathbb{E}\left[ \left\vert y_{t}\right\vert
^{2}+\int_{t}^{T}\sum_{j=1}^{d}\left\vert z_{s}^{j}\right\vert ^{2}\mathrm{d}s\right] 
&\!\!=\!\!&\mathbb{E}\left[ \left\vert \xi \right\vert
^{2}+2\int_{t}^{T}\left\langle y_{s},F_{s}\left( y_{s},z_{s}\right)
\right\rangle \mathrm{d}s\right]  \notag \\
&\!\!\leq\!\! &\mathbb{E}\left[ \left\vert \xi \right\vert
^{2}+2\int_{t}^{T}\left\vert y_{s}\right\vert \left( C_{1}\left\vert
y_{s}\right\vert +C_{1}\sum_{j=1}^{d}\left\vert z_{s}^{j}\right\vert
+\left\vert f_{s}\right\vert \right) \mathrm{d}s\right]  \notag \\
&\!\!\leq\!\! &\mathbb{E}\left[ \left\vert \xi \right\vert
^{2}+2C_{1}\int_{0}^{T}\left\vert y_{s}\right\vert ^{2}\mathrm{d}s+2C_{1}\int_{0}^{T}\sum_{j=1}^{d}\left\vert y_{s}\right\vert \left\vert
z_{s}^{j}\right\vert \mathrm{d}s+2y_{T}^{\ast }\int_{0}^{T}\left\vert
f_{s}\right\vert \mathrm{d}s\right]  \notag \\
&\!\!\leq\!\! &\mathbb{E}\left[ \left\vert \xi \right\vert ^{2}+\left(
2C_{1}+2C_{1}^{2}\right)\!\int_{0}^{T}\!\!\left\vert y_{s}\right\vert ^{2}\mathrm{d}s+\frac{1}{2}\int_{0}^{T}\!\!\sum_{j=1}^{d}\left\vert z_{s}^{j}\right\vert^{2}\mathrm{d}s+2y_{T}^{\ast }\!\int_{0}^{T}\!\!\left\vert f_{s}\right\vert \mathrm{d}s\right].  \label{B1}
\end{eqnarray*}
Re-arranging the above inequality leads to
\begin{eqnarray}
\mathbb{E}\left[ \left\vert y_{t}\right\vert ^{2}+\frac{1}{2}\int_{t}^{T}\sum_{j=1}^{d}\left\vert z_{s}^{j}\right\vert ^{2}\mathrm{d}s\right]  
\leq \mathbb{E}\left[ \left\vert \xi \right\vert ^{2}+\left(
2C_{1}+2C_{1}^{2}d\right) \int_{0}^{T}\left\vert y_{s}\right\vert ^{2}\mathrm{d}s+2y_{T}^{\ast }\int_{0}^{T}\left\vert f_{s}\right\vert \mathrm{d}s\right].  \label{B3}
\end{eqnarray}
By Fubini Theorem
\begin{equation*}
\mathbb{E}\left[ \left\vert y_{t}\right\vert ^{2}\right] \leq \mathbb{E}%
\left[ \left\vert \xi \right\vert ^{2}+2y_{T}^{\ast }\int_{0}^{T}\left\vert
f_{s}\right\vert \mathrm{d}s\right] +\left( 2C_{1}+2C_{1}^{2}d\right)
\int_{0}^{T}\mathbb{E}\left[ \left\vert y_{s}\right\vert ^{2}\right] \mathrm{d}s. 
\end{equation*} 
Applying Gronwall inequality, we obtain%
\begin{equation}
\mathbb{E}\left[ \left\vert y_{t}\right\vert ^{2}\right] \leq e^{\left(
2C_{1}+2C_{1}^{2}d\right) }\mathbb{E}\left[ \left\vert \xi \right\vert
^{2}+2y_{T}^{\ast }\int_{0}^{T}\left\vert f_{s}\right\vert \mathrm{d}s\right].  \label{B2}
\end{equation}
Plugging (\ref{B2}) into (\ref{B3}) leads to 
\begin{equation}
\mathbb{E}\left[ \int_{0}^{T}\sum_{j=1}^{d}\left\vert z_{s}^{j}\right\vert
^{2}\mathrm{d}s\right] \leq C_{3}\mathbb{E}\left[ \left\vert \xi \right\vert
^{2}+2y_{T}^{\ast }\int_{0}^{T}\left\vert f_{s}\right\vert \mathrm{d}s\right],  \label{B4}
\end{equation}%
where $C_{3}=\left( T\left( 2C_{1}+2C_{1}^{2}d\right) e^{\left(2C_{1}+2C_{1}^{2}d\right) }+1\right)$. From (\ref{B2}) and (\ref{B4}), we derive
\begin{equation}
\sup_{0\leq t\leq T}\mathbb{E}\left[ \left\vert y_{t}\right\vert ^{2}\right]
+\mathbb{E}\left[ \int_{0}^{T}\sum_{j=1}^{d}\left\vert z_{s}^{j}\right\vert
^{2}\mathrm{d}s\right] \leq \left( C_{3}+e^{\left( 2C_{1}+2C_{1}^{2}d\right)
}\right) \mathbb{E}\left[ \left\vert \xi \right\vert ^{2}+\epsilon
\left\vert y_{T}^{\ast }\right\vert ^{2}+\epsilon ^{-1}\left(
\int_{0}^{T}\left\vert f_{s}\right\vert \mathrm{d}s\right) ^{2}\right].
\label{B5}
\end{equation}
Substituting (\ref{B5}) into (\ref{B0}), we obtain 
\begin{eqnarray*}
\begin{aligned}
\mathbb{E}\left[ \left\vert y_{T}^{\ast }\right\vert ^{2}\right] 
\leq & \, 4\mathbb{E}\Bigg [\left\vert \xi \right\vert ^{2}+\left(
\int_{0}^{T}\left\vert f_{s}\right\vert \right) ^{2}+\left( C_{1}^{2}\left(
d+1\right) +1\right) \int_{0}^{T}\left( \left\vert y_{s}\right\vert
^{2}+\sum_{j=1}^{d}\left\vert z_{s}^{j}\right\vert ^{2}\right) \mathrm{d}s%
\Bigg ] \\
\leq & \, 4\mathbb{E}\Bigg [\left\vert \xi \right\vert ^{2}+\left(
\int_{0}^{T}\left\vert f_{s}\right\vert \right) ^{2} \\
&+\left( C_{1}^{2}\left( d+1\right) +1\right) \left( C_{3}+e^{\left(
2C_{1}+2C_{1}^{2}d\right) }\right) \left[ \left\vert \xi \right\vert
^{2}+\epsilon \left\vert y_{T}^{\ast }\right\vert ^{2}+\epsilon ^{-1}\left(
\int_{0}^{T}\left\vert f_{s}\right\vert \mathrm{d}s\right) ^{2}\right] \Bigg].
\end{aligned}
\end{eqnarray*}
Taking
$\epsilon =\frac{1}{2\left( C_{1}^{2}\left( d+1\right) +1\right) \left(
C_{3}+e^{\left( 2C_{1}+2C_{1}^{2}d\right) }\right)},$
%
we have 
\begin{eqnarray*}
\begin{aligned}
\mathbb{E}\left[ \left\vert y_{T}^{\ast }\right\vert ^{2}\right] 
\leq & \, 8\mathbb{E}\Bigg [\left( \left( C_{1}^{2}\left( d+1\right) +1\right) \left(
C_{3}+e^{\left( 2C_{1}+2C_{1}^{2}d\right) }\right) +1\right) \left\vert \xi
\right\vert ^{2} \\
&+\left( 2\left( C_{1}^{2}\left( d+1\right) +1\right) ^{2}\left(
C_{3}+e^{\left( 2C_{1}+2C_{1}^{2}d\right) }\right) ^{2}+1\right) \left(
\int_{0}^{T}\left\vert f_{s}\right\vert \right) ^{2}\Bigg]. 
\end{aligned}
\end{eqnarray*}
Hence, using (\ref{B5}), we get 
\begin{eqnarray*} 
\begin{aligned}
&\mathbb{E}\left[ \int_{0}^{T}\sum_{j=1}^{d}\left\vert z_{s}^{j}\right\vert
^{2}\mathrm{d}s\right] \\
\leq & \, \mathbb{E}\Bigg [\left( C_{3}+e^{\left( 2C_{1}+2C_{1}^{2}d\right)
}\right) +8\left( \left( C_{1}^{2}\left( d+1\right) +1\right) \left(
C_{3}+e^{\left( 2C_{1}+2C_{1}^{2}d\right) }\right) +1\right) \left\vert \xi
\right\vert ^{2} \\
&+\left( C_{3}+e^{\left( 2C_{1}+2C_{1}^{2}d\right) }\right) 8\left( 2\left(
C_{1}^{2}\left( d+1\right) +1\right) ^{2}\left( C_{3}+e^{\left(
2C_{1}+2C_{1}^{2}d\right) }\right) ^{2}+1\right) \left(
\int_{0}^{T}\left\vert f_{s}\right\vert \mathrm{d}s\right) ^{2}\Bigg], 
\end{aligned}
\end{eqnarray*}
which implies the desired result. \hfill $\Box $

\paragraph{Proof of Lemma \protect\ref{l2}.}

Applying It\^{o}'s formula to $\left\vert X_{t,i}\right\vert ^{p}$ yields 
\begin{eqnarray}
\left\vert X_{t,i}\right\vert ^{p} 
&\!\!\!=\!\!\!&\left\vert \xi _{i}\right\vert
^{p}+\int_{0}^{t}\Big (p\left\vert X_{s,i}\right\vert ^{p-1}b_{i}\left(
s,X_{s,i},\mathbf{X}_{s}^{\mathbf{u}},u_{s,i}\right)  
+\frac{p\left( p-1\right) }{2}\left\vert X_{s,i}\right\vert ^{p-2}\sigma
_{i}\left( s,X_{s,i},\mathbf{X}_{s}^{\mathbf{u}},u_{s,i}\right) ^{2}\Big )%
\mathrm{d}s  \notag \\
&&+\int_{0}^{t}p\left\vert X_{s,i}\right\vert ^{p-1}\sigma _{i}\left(
s,X_{s,i},\mathbf{X}_{s}^{\mathbf{u}},u_{s,i}\right) \mathrm{d}W_{s}^{i}.
\label{ineq1}
\end{eqnarray}
Taking the expectation of both sides of (\ref{ineq1}, by the fact
that $\int_{0}^{t}p\left\vert X_{s,i}\right\vert ^{p-1}\sigma _{i}\left(
s,X_{s,i},\mathbf{X}_{s}^{\mathbf{u}},u_{s,i}\right) \mathrm{d}W_{s}^{i}$ is
a martingale (see [29, Problem 2.10.7]), and by Remark \ref{re1}, we get 
\begin{eqnarray*}
\begin{aligned}
\mathbb{E}\left[ \left\vert X_{t,i}\right\vert ^{p}\right] 
=& \, \mathbb{E}\Bigg [\left\vert \xi _{i}\right\vert ^{p}+\int_{0}^{t}\Big (p\left\vert
X_{s,i}\right\vert ^{p-1}b_{i}\left( s,X_{s,i},\mathbf{X}_{s}^{\mathbf{u}%
},u_{s,i}\right) \\
&+\frac{p\left( p-1\right) }{2}\left\vert X_{s,i}\right\vert ^{p-2}\sigma
_{i}\left( s,X_{s,i},\mathbf{X}_{s}^{\mathbf{u}},u_{s,i}\right) ^{2}\Big )%
\mathrm{d}s\Bigg ] \\
\leq & \, \mathbb{E}\Bigg \{\left\vert \xi _{i}\right\vert ^{p}+\int_{0}^{t}%
\Bigg [p\left\vert X_{s,i}\right\vert ^{p-1}\left( L^{b}\left( 1+\left\vert
X_{s,i}\right\vert +\left\vert u_{s,i}\right\vert \right) +\frac{L_{y}^{b}}{N%
}\sum_{k=1}^{N}\left\vert X_{s,k}\right\vert \right) \\
&+\frac{p\left( p-1\right) }{2}\left\vert X_{s,i}\right\vert ^{p-2}\left(
L^{\sigma }\left( 1+\left\vert X_{s,i}\right\vert +\left\vert
u_{s,i}\right\vert \right) +\frac{L_{y}^{\sigma }}{N}\sum_{k=1}^{N}\left%
\vert X_{s,k}\right\vert \right) ^{2}\Bigg ]\mathrm{d}s\Bigg\}. 
\end{aligned}
\end{eqnarray*}
By Young's inequality, we have 
\begin{eqnarray}
\mathbb{E}\left[ \left\vert X_{t,i}\right\vert ^{p}\right]  &\leq &\mathbb{E}%
\Bigg \{\left\vert \xi _{i}\right\vert ^{p}+\int_{0}^{t}\Bigg [L^{b}\left(
\left( 3p-2\right) \left\vert X_{s,i}\right\vert ^{p}+1+\left\vert
u_{s,i}\right\vert ^{p}\right) +L_{y}^{b}\left( p-1\right) \left\vert
X_{s,i}\right\vert ^{p}+\frac{L_{y}^{b}}{N}\sum_{k=1}^{N}\left\vert
X_{s,k}\right\vert ^{p}  \notag \\
&&+\frac{4p\left( p-1\right) }{2}\left\vert X_{s,i}\right\vert ^{p-2}\Bigg (%
\left( L^{\sigma }\right) ^{2}+\left( L^{\sigma }\right) ^{2}\left\vert
X_{s,i}\right\vert ^{2}+\left( L^{\sigma }\right) ^{2}\left\vert
u_{s,i}\right\vert ^{2}  \notag \\
&&+\left( \frac{L_{y}^{\sigma }}{N}\right) ^{2}\left(
\sum_{k=1}^{N}\left\vert X_{s,k}\right\vert \right) ^{2}\Bigg )\Bigg ]%
\mathrm{d}s\Bigg \}  \notag \\
&\leq &\mathbb{E}\Bigg \{\left\vert \xi _{i}\right\vert ^{p}+\int_{0}^{t}%
\Bigg [L^{b}\left( \left( 3p-2\right) \left\vert X_{s,i}\right\vert
^{p}+1+\left\vert u_{s,i}\right\vert ^{p}\right) +L_{y}^{b}\left( p-1\right)
\left\vert X_{s,i}\right\vert ^{p}+\frac{L_{y}^{b}}{N}\sum_{k=1}^{N}\left%
\vert X_{s,k}\right\vert ^{p}  \notag \\
&&+\frac{4p\left( p-1\right) }{2}\left\vert X_{s,i}\right\vert ^{p-2}\Bigg (%
\left( L^{\sigma }\right) ^{2}+\left( L^{\sigma }\right) ^{2}\left\vert
X_{s,i}\right\vert ^{2}+\left( L^{\sigma }\right) ^{2}\left\vert
u_{s,i}\right\vert ^{2}  \notag \\
&&+N\left( \frac{L_{y}^{\sigma }}{N}\right) ^{2}\sum_{k=1}^{N}\left\vert
X_{s,k}\right\vert ^{2}\Bigg )\Bigg ]\mathrm{d}s\Bigg \}  \notag \\
&\leq &\mathbb{E}\Bigg \{\left\vert \xi _{i}\right\vert ^{p}+\int_{0}^{t}%
\Bigg [L^{b}\left( \left( 3p-2\right) \left\vert X_{s,i}\right\vert
^{p}+1+\left\vert u_{s,i}\right\vert ^{p}\right) +L_{y}^{b}\left( p-1\right)
\left\vert X_{s,i}\right\vert ^{p}+\frac{L_{y}^{b}}{N}\sum_{k=1}^{N}\left%
\vert X_{s,k}\right\vert ^{p}  \notag \\
&&+2p\left( p-1\right) \left\vert X_{s,i}\right\vert ^{p-2}\Bigg (\left(
L^{\sigma }\right) ^{2}+\left( L^{\sigma }\right) ^{2}\left\vert
X_{s,i}\right\vert ^{2}+\left( L^{\sigma }\right) ^{2}\left\vert
u_{s,i}\right\vert ^{2}  \notag \\
&&+\frac{\left( L_{y}^{\sigma }\right) ^{2}}{N}\sum_{k=1}^{N}\left\vert
X_{s,k}\right\vert ^{2}\Bigg )\Bigg ]\mathrm{d}s\Bigg \}  \notag \\
&\leq &\mathbb{E}\Bigg \{\left\vert \xi _{i}\right\vert ^{p}+\int_{0}^{t}%
\Bigg [L^{b}\left( \left( 3p-2\right) \left\vert X_{s,i}\right\vert
^{p}+1+\left\vert u_{s,i}\right\vert ^{p}\right) +L_{y}^{b}\left( p-1\right)
\left\vert X_{s,i}\right\vert ^{p}+\frac{L_{y}^{b}}{N}\sum_{k=1}^{N}\left%
\vert X_{s,k}\right\vert ^{p}  \notag \\
&&+2p\left( p-1\right) \left( L^{\sigma }\right) ^{2}\left\vert
X_{s,i}\right\vert ^{p}+2\left( p-1\right) \left( p-2\right) \left(
L^{\sigma }\right) ^{2}\left\vert X_{s,i}\right\vert ^{p}+4\left( p-1\right)
\left( L^{\sigma }\right) ^{2}  \notag \\
&&+2\left( p-1\right) \left( p-2\right) \left( L^{\sigma }\right)
^{2}\left\vert X_{s,i}\right\vert ^{p}+4\left( p-1\right) \left( L^{\sigma
}\right) ^{2}\left\vert u_{s,i}\right\vert ^{p}  \notag \\
&&+2\left( p-1\right) \left( p-2\right) \left( L_{y}^{\sigma }\right)
^{2}\left\vert X_{s,i}\right\vert ^{p}+\frac{4\left( p-1\right) \left(
L_{y}^{\sigma }\right) ^{2}}{N}\sum_{k=1}^{N}\left\vert X_{s,k}\right\vert
^{p}\Bigg ]\mathrm{d}s\Bigg \}  \notag \\
&\leq &\mathbb{E}\Bigg \{\left\vert \xi _{i}\right\vert ^{p}+\int_{0}^{t}%
\Bigg [\big (L^{b}\left( 3p-2\right) +L_{y}^{b}\left( p-1\right) +2p\left(
p-1\right) \left( L^{\sigma }\right) ^{2}  \notag \\
&&+4\left( p-1\right) \left( p-2\right) \left( L^{\sigma }\right)
^{2}+2\left( p-1\right) \left( p-2\right) \left( L_{y}^{\sigma }\right) ^{2}%
\big )\left\vert X_{s,i}\right\vert ^{p}  \notag \\
&&+L^{b}+4\left( p-1\right) \left( L^{\sigma }\right) ^{2}+\left(
L^{b}+4\left( p-1\right) \left( L^{\sigma }\right) ^{2}\right) \left\vert
u_{s,i}\right\vert ^{p}  \notag \\
&&+\frac{L_{y}^{b}+4\left( p-1\right) \left( L_{y}^{\sigma }\right) ^{2}}{N}%
\sum_{k=1}^{N}\left\vert X_{s,k}\right\vert ^{p}\Bigg ]\mathrm{d}s\Bigg \}
\notag \\
&\leq &\mathbb{E}\left[ \left\vert \xi _{i}\right\vert ^{p}\right] +\left(
L^{b}+4\left( p-1\right) \left( L^{\sigma }\right) ^{2}\right) T+\left(
L^{b}+4\left( p-1\right) \left( L^{\sigma }\right) ^{2}\right) \left\Vert
u_{i}\right\Vert _{\mathcal{H}^{p}\left( \mathbb{R}\right) }^{p}  \notag \\
&&+\int_{0}^{t}\Bigg [\big (L^{b}\left( 3p-2\right) +L_{y}^{b}\left(
p-1\right) +2p\left( p-1\right) \left( L^{\sigma }\right) ^{2}  \notag \\
&&+4\left( p-1\right) \left( p-2\right) \left( L^{\sigma }\right)
^{2}+2\left( p-1\right) \left( p-2\right) \left( L_{y}^{\sigma }\right) ^{2}%
\big )\mathbb{E}\left[ \left\vert X_{s,i}\right\vert ^{p}\right]   \notag \\
&&+\frac{L_{y}^{b}+4\left( p-1\right) \left( L_{y}^{\sigma }\right) ^{2}}{N}%
\sum_{k=1}^{N}\mathbb{E}\left[ \left\vert X_{s,k}\right\vert ^{p}\right] %
\Bigg ]\mathrm{d}s.  \label{ineq2}
\end{eqnarray}
Summing up the above equation over the index $i\in I_{N}$ yields for all $t\in [0,T]$,
\begin{eqnarray*}
\begin{aligned} 
\sum_{i=1}^{N}\mathbb{E}\left[ \left\vert X_{t,i}\right\vert ^{p}\right]
\leq &\sum_{i=1}^{N}I_{\xi _{i},b,\sigma ,p,u_{i}}^{0}+\int_{0}^{t}\Bigg [%
\big (L^{b}\left( 3p-2\right) 
+L_{y}^{b}p+\left( 6p^{2}+14p+8\right) \left( L^{\sigma }\right) ^{2} \\
&+2\left( p-1\right) p\left( L_{y}^{\sigma }\right) ^{2}\big )\sum_{k=1}^{N}%
\mathbb{E}\left[ \left\vert X_{s,k}\right\vert ^{p}\right] \Bigg ]\mathrm{d}s.
\end{aligned}
\end{eqnarray*}
Immediately, Gronwall's inequality gives
\begin{equation*}
\sum_{k=1}^{N}\mathbb{E}\left[ \left\vert X_{t,k}\right\vert ^{p}\right]
\leq \sum_{k=1}^{N}I_{\xi _{k},b,\sigma ,p,u_{k}}^{0}\cdot e^{I_{b,\sigma,p}^{1}\cdot t}.
\end{equation*}
Substituting the above inequality into (\ref{ineq2}) and applying Gronwall's
inequality again yield
\begin{eqnarray*}
\begin{aligned}
\mathbb{E}\left[ \left\vert X_{t,i}\right\vert ^{p}\right] 
&\leq I_{\xi_{i},b,\sigma ,p,u_{i}}^{0}+\int_{0}^{t}\Bigg [\big (I_{b,\sigma ,p}^{2}\big
)\mathbb{E}\left[ \left\vert X_{s,i}\right\vert ^{p}\right] 
+\frac{L_{y}^{b}+4\left( p-1\right) \left( L_{y}^{\sigma }\right) ^{2}}{N}%
\sum_{k=1}^{N}\mathbb{E}\left[ \left\vert X_{s,k}\right\vert ^{p}\right] %
\Bigg ]\mathrm{d}s \\
&\leq I_{\xi _{i},b,\sigma ,p,u_{i}}^{0}+\int_{0}^{t}\Bigg [\big (%
I_{b,\sigma ,p}^{2}\big )\mathbb{E}\left[ \left\vert X_{s,i}\right\vert ^{p}%
\right] 
+\frac{L_{y}^{b}+4\left( p-1\right) \left( L_{y}^{\sigma }\right) ^{2}}{N}%
\sum_{k=1}^{N}I_{\xi _{k},b,\sigma ,p,u_{k}}^{0}\cdot e^{I_{b,\sigma
,p}^{1}\cdot t},
\end{aligned}
\end{eqnarray*}
and then
\begin{equation*}
\mathbb{E}\left[ \left\vert X_{t,i}\right\vert ^{p}\right] \leq \left[
I_{\xi _{i},b,\sigma ,p,u_{i}}^{0}+\frac{L_{y}^{b}+4\left( p-1\right) \left(
L_{y}^{\sigma }\right) ^{2}}{N}\sum_{k=1}^{N}I_{\xi _{k},b,\sigma
,p,u_{k}}^{0}\cdot e^{I_{b,\sigma ,p}^{1}\cdot t}\right] \cdot
e^{I_{b,\sigma ,p}^{2}\cdot t}.
\end{equation*}%
We thus complete the proof. \hfill $\Box $
\begin{lemma}
\label{l3}Fix $p\geq 2$ and for each $i,j\in I_{N}$, let $B_{i},\bar{B}%
_{i,j},D_{i},\bar{D}_{i,j}:\Omega \times \left[ 0,T\right] \rightarrow
\mathbb{R}$ be bounded adapted processes, and $f_{i},\bar{f}_{i}\in \mathcal{%
H}^{p}\left( \mathbb{R}\right) $. Let $\mathbf{S}=(S_{i})_{i=1}^{N}\in
\mathcal{S}^{p}\left( \mathbb{R}^{N}\right) $ satisfy the following
dynamics: for all $t\in [0,T]$,
\begin{eqnarray*}
\begin{aligned}
\mathrm{d}S_{t,i} 
&=\left( B_{i}\left( t\right) S_{t,i}+\sum_{j=1}^{N}\bar{B}_{i,j}\left( t\right) S_{t,j}+f_{t,i}\right) \mathrm{d}t 
+\left( D_{i}\left( t\right) S_{t,i}+\sum_{j=1}^{N}\bar{D}_{i,j}\left(
t\right) S_{t,j}+\bar{f}_{t,i}\right) \mathrm{d}W_{t}^{i}, \\
S_{0,i} &=0,\qquad \forall i\in I_{N}.
\end{aligned}
\end{eqnarray*}
Then for all $i\in I_{N}$, we have
\begin{eqnarray*}
\begin{aligned}
\mathbb{E}\left[ \left\vert S_{t,i}\right\vert ^{p}\right] 
\leq & \, \Bigg[T\left( \left\Vert \bar{B}\right\Vert _{\infty }+3\left( p-1\right)
N\left\Vert \bar{D}\right\Vert _{\infty }^{2}\right) 
\left( \sum_{k=1}^{N}\left( \left\Vert f_{k}\right\Vert _{\mathcal{H}^{p}\left( \mathbb{R}\right) }^{p}+3\left( p-1\right) \left\Vert \bar{f}_{k}\right\Vert _{\mathcal{H}^{p}\left( \mathbb{R}\right) }^{p}\right)\right) e^{I_{B,D,\bar{B},\bar{D},p}^{3}\cdot T} \\
&+\left\Vert f_{i}\right\Vert _{\mathcal{H}^{p}\left( \mathbb{R}\right)}^{p}+3\left( p-1\right) \left\Vert \bar{f}_{i}\right\Vert _{\mathcal{H}^{p}\left( \mathbb{R}\right) }^{p}\Bigg ]e^{I_{B,D,\bar{B},\bar{D},p}^{4}\cdot T},
\end{aligned}
\end{eqnarray*}
where
\begin{eqnarray*}
\begin{aligned}
I_{B,D,\bar{B},\bar{D},p}^{3} 
= & \, p\left\Vert B\right\Vert _{\infty}+N\left\Vert \bar{B}\right\Vert _{\infty }p+\left( p-1\right) +\frac{3}{2}\left( p-1\right) p\left\Vert D\right\Vert _{\infty }^{2} \\
&+\frac{3}{2}\left( p-1\right) \left( p-2\right) N^{2}\left\Vert \bar{D}\right\Vert _{\infty}^{2}+\frac{3}{2}\left( p-1\right) \left( p-2\right)
+3\left( p-1\right) N^{2}\left\Vert \bar{D}\right\Vert _{\infty }^{2}, \\
I_{B,D,\bar{B},\bar{D},p}^{4} 
= & \, p\left\Vert B\right\Vert _{\infty}+N\left\Vert \bar{B}\right\Vert _{\infty }\left( p-1\right) +\left(p-1\right) +\frac{3}{2}\left( p-1\right) p\left\Vert D\right\Vert _{\infty}^{2} \\
&+\frac{3}{2}\left( p-1\right) N^{2}\left\Vert \bar{D}\right\Vert_{\infty}^{2}\left( p-2\right) +\frac{3}{2}\left( p-1\right) \left( p-2\right), 
\end{aligned}
\end{eqnarray*}
where $\left\Vert B\right\Vert _{\infty }=\max_{i\in I_{N}}\left\Vert
B_{i}\right\Vert _{L^{\infty }},$ $\left\Vert D\right\Vert _{\infty
}=\max_{i\in I_{N}}\left\Vert D_{i}\right\Vert _{L^{\infty }}$, $\left\Vert
\bar{B}\right\Vert _{\infty }=\max_{i,j\in I_{N}}\left\Vert \bar{B}%
_{i,j}\right\Vert _{L^{\infty }}$, and $\left\Vert \bar{D}\right\Vert
_{\infty}=\max_{i,j\in I_{N}}\left\Vert \bar{D}_{i,j}\right\Vert_{L^{\infty }}$.
\end{lemma}

\paragraph{Proof.}

Applying It\^{o}'s formula to $\left\vert S_{t,i}\right\vert ^{p}$, we obtain 
\begin{eqnarray*} 
\begin{aligned}
\left\vert S_{t,i}\right\vert ^{p} 
=& \, \int_{0}^{t}p\left\vert
S_{r,i}\right\vert ^{p-1}\left( B_{i}\left( r\right) S_{r,i}+\sum_{k=1}^{N}%
\bar{B}_{i,k}\left( r\right) S_{r,k}+f_{r,i}\right) \mathrm{d}r \\
&+\int_{0}^{t}p\left\vert S_{r,i}\right\vert ^{p-1}\left( D_{i}\left(t\right) S_{t,i}+\sum_{k=1}^{N}\bar{D}_{i,k}\left( t\right) S_{t,k}+\bar{f}_{t,i}\right) \mathrm{d}W_{t}^{i} \\
&+\frac{1}{2}p\left( p-1\right) \int_{0}^{t}\left\vert S_{r,i}\right\vert^{p-2}\left( D_{i}\left( r\right) S_{r,i}+\sum_{k=1}^{N}\bar{D}_{i,k}\left(r\right) S_{r,k}+\bar{f}_{r,i}\right) ^{2}\mathrm{d}r.
\end{aligned}
\end{eqnarray*}
Taking the expectation leads to 
\begin{eqnarray*}
\begin{aligned}
\mathbb{E}\left[ \left\vert S_{t,i}\right\vert ^{p}\right] 
= & \, \mathbb{E}\Bigg [\displaystyle\int_{0}^{t}p\left\vert S_{r,i}\right\vert ^{p-1}\left(
B_{i}\left( r\right) S_{r,i}+\sum_{k=1}^{N}\bar{B}_{i,k}\left( r\right)
S_{r,k}+f_{r,i}\right) \mathrm{d}r \\
&+\frac{1}{2}p\left( p-1\right) \int_{0}^{t}\left\vert S_{r,i}\right\vert^{p-2}\left( D_{i}\left( r\right) S_{r,i}+\sum_{k=1}^{N}\bar{D}_{i,k}\left(r\right) S_{r,k}+\bar{f}_{r,i}\right) ^{2}\mathrm{d}r\Bigg].
\end{aligned}
\end{eqnarray*}
Applying Young's inequality yields 
\begin{eqnarray}
\mathbb{E}\left[ \left\vert S_{t,i}\right\vert ^{p}\right] 
&\leq &\mathbb{E}\Bigg [\displaystyle\int_{0}^{t}\Big (p\left\Vert B\right\Vert _{\infty
}\left\vert S_{r,i}\right\vert ^{p}+p\left\Vert \bar{B}\right\Vert _{\infty
}\sum_{k=1}^{N}\left\vert S_{r,i}\right\vert ^{p-1}\left\vert
S_{r,k}\right\vert +p\left\vert S_{r,i}\right\vert ^{p-1}\left\vert
f_{r,i}\right\vert \Big )\mathrm{d}r  \notag \\
&&+\frac{1}{2}p\left( p-1\right) \int_{0}^{t}\left\vert S_{r,i}\right\vert
^{p-2}\left( \left\Vert D\right\Vert _{\infty }\left\vert S_{r,i}\right\vert
+\left\Vert \bar{D}\right\Vert _{\infty }\sum_{k=1}^{N}\left\vert
S_{r,k}\right\vert +\left\vert \bar{f}_{r,i}\right\vert\right) ^{2}\mathrm{d}r\Bigg ]  \notag \\
&\leq &\mathbb{E}\Bigg [\displaystyle\int_{0}^{t}\Big (p\left\Vert B\right\Vert _{\infty }\left\vert S_{r,i}\right\vert ^{p}+p\left\Vert \bar{B}\right\Vert _{\infty }\sum_{k=1}^{N}\left( \frac{p-1}{p}\left\vert
S_{r,i}\right\vert ^{p}+\frac{1}{p}\left\vert S_{r,k}\right\vert ^{p}\right)
+p\left( \frac{p-1}{p}\left\vert S_{r,i}\right\vert^{p}+\frac{1}{p}\left\vert f_{r,i}\right\vert ^{p}\right) \Big )\mathrm{d}r  \notag \\
&&+\frac{3}{2}p\left( p-1\right) \int_{0}^{t}\left\vert S_{r,i}\right\vert
^{p-2}\left( \left\Vert D\right\Vert _{\infty }^{2}\left\vert
S_{r,i}\right\vert ^{2}+N\left\Vert \bar{D}\right\Vert _{\infty
}^{2}\sum_{k=1}^{N}\left\vert S_{r,k}\right\vert ^{2}+\left\vert \bar{f}%
_{r,i}\right\vert ^{2}\right) \mathrm{d}r\Bigg ]  \notag \\
&\leq &\mathbb{E}\Bigg [\displaystyle\int_{0}^{t}\Big (p\left\Vert
B\right\Vert _{\infty }\left\vert S_{r,i}\right\vert ^{p}+\left\Vert \bar{B}\right\Vert_{\infty }\sum_{k=1}^{N}\left( \left( p-1\right) \left\vert
S_{r,i}\right\vert ^{p}+\left\vert S_{r,k}\right\vert^{p}\right)  +\left( p-1\right) \left\vert S_{r,i}\right\vert^{p}+\left\vert
f_{r,i}\right\vert ^{p}\Big )\mathrm{d}r  \notag \\
&&+\frac{3}{2}p\left( p-1\right) \int_{0}^{t}\Big (\left\Vert D\right\Vert_{\infty }^{2}\left\vert S_{r,i}\right\vert^{p}+N\left\Vert \bar{D}\right\Vert _{\infty}^{2}\sum_{k=1}^{N}\left\vert S_{r,i}\right\vert
^{p-2}\left\vert S_{r,k}\right\vert ^{2}
+\left\vert S_{r,i}\right\vert ^{p-2}\left\vert \bar{f}_{r,i}\right\vert
^{2}\Big )\mathrm{d}r\Bigg ]  \notag \\
&\leq &\mathbb{E}\Bigg [\displaystyle\int_{0}^{t}\Big (p\left\Vert
B\right\Vert _{\infty }\left\vert S_{r,i}\right\vert ^{p}+\left\Vert \bar{B}%
\right\Vert _{\infty }\sum_{k=1}^{N}\left( \left( p-1\right) \left\vert
S_{r,i}\right\vert ^{p}+\left\vert S_{r,k}\right\vert ^{p}\right)  +\left( p-1\right) \left\vert S_{r,i}\right\vert ^{p}+\left\vert
f_{r,i}\right\vert ^{p}\Big )\mathrm{d}r  \notag \\
&&+\frac{3}{2}\left( p-1\right) \int_{0}^{t}\Big (p\left\Vert D\right\Vert
_{\infty }^{2}\left\vert S_{r,i}\right\vert ^{p}+N\left\Vert \bar{D}%
\right\Vert _{\infty }^{2}\sum_{k=1}^{N}\left( \left( p-2\right) \left\vert S_{r,i}\right\vert ^{p}+2\left\vert S_{r,k}\right\vert^{p}\right)  \notag \\
&&+\left( p-2\right) \left\vert S_{r,i}\right\vert^{p}+2\left\vert \bar{f}_{r,i}\right\vert ^{p}\Big )\mathrm{d}r\Bigg ]  \notag \\
&\leq &\mathbb{E}\Bigg \{\displaystyle\int_{0}^{t}\Bigg [\Big (p\left\Vert
B\right\Vert _{\infty }+N\left\Vert \bar{B}\right\Vert _{\infty }\left(
p-1\right) +\left( p-1\right) +\frac{3}{2}\left( p-1\right) p\left\Vert
D\right\Vert _{\infty }^{2}  \notag \\
&&+\frac{3}{2}\left( p-1\right) N^{2}\left\Vert \bar{D}\right\Vert _{\infty
}^{2}\left( p-2\right) +\frac{3}{2}\left( p-1\right) \left( p-2\right) \Big )%
\left\vert S_{r,i}\right\vert ^{p}  \notag \\
&&+\left( \left\Vert \bar{B}\right\Vert _{\infty }+3\left( p-1\right)
N\left\Vert \bar{D}\right\Vert _{\infty }^{2}\right)
\sum_{k=1}^{N}\left\vert S_{r,k}\right\vert^p)
+\left\vert f_{r,i}\right\vert ^{p}+3\left( p-1\right) \left\vert \bar{f}%
_{r,i}\right\vert ^{p}\Big )\Bigg ]\mathrm{d}r\Bigg\}.  \label{ineq3}
\end{eqnarray}
Summarizing (\ref{ineq3}) over the index $i\in I_{N}$ yields for all $t\in[0,T]$,
\begin{eqnarray*} 
\begin{aligned}
\sum_{k=1}^{N}\mathbb{E}\left[ \left\vert S_{t,k}\right\vert ^{p}\right]
\leq & \,\mathbb{E}\Bigg \{\int_{0}^{t}\Bigg [\Big (p\left\Vert B\right\Vert_{\infty }+N\left\Vert \bar{B}\right\Vert _{\infty }\left( p-1\right)
+\left( p-1\right) +\frac{3}{2}\left( p-1\right) p\left\Vert D\right\Vert_{\infty }^{2} \\
&+\frac{3}{2}\left( p-1\right) N^{2}\left\Vert \bar{D}\right\Vert _{\infty
}^{2}\left( p-2\right) +\frac{3}{2}\left( p-1\right) \left( p-2\right)\Big)\sum_{k=1}^{N}\left\vert S_{r,k}\right\vert ^{p} \\
&+N\left( \left\Vert \bar{B}\right\Vert _{\infty }+3\left( p-1\right)N\left\Vert \bar{D}\right\Vert _{\infty }^{2}\right)
\sum_{k=1}^{N}\left\vert S_{r,k}\right\vert ^{p} +\sum_{k=1}^{N}\left( \left\vert f_{r,k}\right\vert ^{p}+3\left(
p-1\right) \left\vert \bar{f}_{r,k}\right\vert ^{p}\right) \Big )\Bigg ]%
\mathrm{d}r\Bigg \} \\
=& \, \mathbb{E}\Bigg \{\int_{0}^{t}\Bigg [\Big (p\left\Vert B\right\Vert
_{\infty }+N\left\Vert \bar{B}\right\Vert _{\infty }p+\left( p-1\right) +%
\frac{3}{2}\left( p-1\right) p\left\Vert D\right\Vert _{\infty }^{2} \\
&+\frac{3}{2}\left( p-1\right) \left( p-2\right) N^{2}\left\Vert \bar{D}\right\Vert _{\infty }^{2}+\frac{3}{2}\left( p-1\right) \left( p-2\right) \\
&+3\left( p-1\right) N^{2}\left\Vert \bar{D}\right\Vert _{\infty }^{2}\Big )\sum_{k=1}^{N}\left\vert S_{r,k}\right\vert ^{p} 
+\sum_{k=1}^{N}\left( \left\vert f_{r,k}\right\vert ^{p}+3\left(
p-1\right) \left\vert \bar{f}_{r,k}\right\vert ^{p}\right) \Big )\Bigg ]\mathrm{d}r\Bigg\}.
\end{aligned}
\end{eqnarray*}
Gronwall's inequality implies that
\begin{eqnarray*}
\begin{aligned}
\sum_{k=1}^{N}\mathbb{E}\left[ \left\vert S_{t,k}\right\vert ^{p}\right]
&\leq \mathbb{E}\left[ \left( \int_{0}^{T}\sum_{k=1}^{N}\left( \left\vert
f_{r,k}\right\vert ^{p}+3\left( p-1\right) \left\vert \bar{f}%
_{r,k}\right\vert ^{p}\right) \mathrm{d}r\right) e^{I_{B,D,\bar{B},\bar{D}%
,p}^{3}\cdot T}\right] \\
&=\mathbb{E}\left[ \left( \sum_{k=1}^{N}\int_{0}^{T}\left\vert
f_{r,k}\right\vert ^{p}+3\left( p-1\right)
\sum_{k=1}^{N}\int_{0}^{T}\left\vert \bar{f}_{r,k}\right\vert ^{p}\mathrm{d}%
r\right) e^{I_{B,D,\bar{B},\bar{D},p}^{3}\cdot T}\right] \\
&=\left( \sum_{k=1}^{N}\mathbb{E}\left[ \int_{0}^{T}\left\vert
f_{r,k}\right\vert ^{p}\mathrm{d}r\right] +3\left( p-1\right) \sum_{k=1}^{N}%
\mathbb{E}\left[ \int_{0}^{T}\left\vert \bar{f}_{r,k}\right\vert ^{p}\mathrm{%
d}r\right] \right) e^{I_{B,D,\bar{B},\bar{D},p}^{3}\cdot T} \\
&=\left[ \sum_{k=1}^{N}\left( \left\Vert f_{k}\right\Vert _{\mathcal{H}^{p}\left( \mathbb{R}\right) }^{p}+3\left( p-1\right) \left\Vert \bar{f}_{k}\right\Vert _{\mathcal{H}^{p}\left( \mathbb{R}\right) }^{p}\right)\right] e^{I_{B,D,\bar{B},\bar{D},p}^{3}\cdot T}. 
\end{aligned}
\end{eqnarray*}
Substituting the above inequality into (\ref{ineq3}) yields for all $t>0$,
\begin{eqnarray*}
\begin{aligned}
\mathbb{E}\left[ \left\vert S_{t,i}\right\vert ^{p}\right] \leq & \, \mathbb{E}\Bigg [\int_{0}^{t}\Bigg [\Big (p\left\Vert B\right\Vert _{\infty
}+N\left\Vert \bar{B}\right\Vert _{\infty }\left( p-1\right) +\left(p-1\right) +\frac{3}{2}\left( p-1\right) p\left\Vert D\right\Vert _{\infty}^{2} \\
&+\frac{3}{2}\left( p-1\right) N^{2}\left\Vert \bar{D}\right\Vert _{\infty
}^{2}\left( p-2\right) +\frac{3}{2}\left( p-1\right) \left( p-2\right) \Big)\mathbb{E}\left[ \left\vert S_{r,i}\right\vert ^{p}\right] \\
&+\left( \left\Vert \bar{B}\right\Vert _{\infty }+3\left( p-1\right)
N\left\Vert \bar{D}\right\Vert _{\infty }^{2}\right) \mathbb{E}\left[
\sum_{k=1}^{N}\left\vert S_{r,k}\right\vert ^{p}\right] +\mathbb{E}\left[ \left\vert f_{r,i}\right\vert ^{p}+3\left( p-1\right)\left\vert \bar{f}_{r,i}\right\vert ^{p}\right] \Big )\Bigg ]\mathrm{d}r \\
\leq & \, \int_{0}^{t}\Big (p\left\Vert B\right\Vert _{\infty }+N\left\Vert\bar{B}\right\Vert _{\infty }\left( p-1\right) +\left( p-1\right) +\frac{3}{2}\left( p-1\right) p\left\Vert D\right\Vert _{\infty }^{2} \\
&+\frac{3}{2}\left( p-1\right) N^{2}\left\Vert \bar{D}\right\Vert _{\infty
}^{2}\left( p-2\right) +\frac{3}{2}\left( p-1\right) \left( p-2\right)\Big)\mathbb{E}\left[ \left\vert S_{r,i}\right\vert ^{p}\right] \mathrm{d}r \\
&+\int_{0}^{t}\Bigg \{\left( \left\Vert \bar{B}\right\Vert _{\infty
}+3\left( p-1\right) N\left\Vert \bar{D}\right\Vert _{\infty }^{2}\right) \cdot \left[ \sum_{k=1}^{N}\left( \left\Vert f_{k}\right\Vert _{\mathcal{H}%
^{p}\left( \mathbb{R}\right) }^{p}+3\left( p-1\right) \left\Vert \bar{f}_{k}\right\Vert _{\mathcal{H}^{p}\left( \mathbb{R}\right) }^{p}\right)\right] e^{I_{B,D,\bar{B},\bar{D},p}^{3}\cdot T} \\
&+\left\vert f_{r,i}\right\vert ^{p}+3\left( p-1\right) \left\vert \bar{f}_{r,i}\right\vert ^{p}\Bigg \}\mathrm{d}r, 
\end{aligned}
\end{eqnarray*}
which implies the desired result.\hfill $\Box $

\paragraph{Proof of Lemma \protect\ref{l4}.}

First of all, we have 
\begin{eqnarray*} 
\begin{aligned}
\mathbb{E}\left[ \left\vert S_{t,i}\right\vert ^{p}\right] \leq &\Bigg[T\left( \left\Vert \bar{B}\right\Vert _{\infty }+3\left( p-1\right)
N\left\Vert \bar{D}\right\Vert _{\infty }^{2}\right) 
\left( \sum_{k=1}^{N}\left( \left\Vert f_{k}\right\Vert _{\mathcal{H}%
^{p}\left( \mathbb{R}\right) }^{p}+3\left( p-1\right) \left\Vert \bar{f}%
_{k}\right\Vert _{\mathcal{H}^{p}\left( \mathbb{R}\right) }^{p}\right)
\right) e^{I_{B,D,\bar{B},\bar{D},p}^{3}\cdot T} \\
&+\left\Vert f_{i}\right\Vert _{\mathcal{H}^{p}\left( \mathbb{R}\right)
}^{p}+3\left( p-1\right) \left\Vert \bar{f}_{i}\right\Vert _{\mathcal{H}%
^{p}\left( \mathbb{R}\right) }^{p}\Bigg ]e^{I_{B,D,\bar{B},\bar{D},p}^{4}\cdot T}. 
\end{aligned}
\end{eqnarray*}
Employing Lemma \ref{l3} with
\begin{eqnarray*}
\begin{aligned}
\mathbf{S} =&\mathbf{Y}^{\mathbf{u},u_{h}^{\prime }},B_{i}\left( t\right)
=\left( \partial _{x}b_{i}\right) \left(t,X_{t,i}^{\mathbf{u}},\mathbf{X}_{t}^{\mathbf{u}},u_{t,i}\right), \\
\bar{B}_{i,j}\left( t\right) =&\left( \partial _{y_{j}}b_{i}\right) \left(
t,X_{t,i}^{\mathbf{u}},\mathbf{X}_{t}^{\mathbf{u}},u_{t,i}\right)
,f_{t,i}=\left( \partial _{u}b_{i}\right) \left( t,X_{t,i}^{\mathbf{u}},%
\mathbf{X}_{t}^{\mathbf{u}},u_{t,i}\right) \delta _{h,i}u_{t,h}^{\prime }, \\
D_{i}\left( t\right) =&\left( \partial _{x}\sigma _{i}\right) \left(t,X_{t,i}^{\mathbf{u}},\mathbf{X}_{t}^{\mathbf{u}},u_{t,i}\right) ,\bar{D}_{i,j}\left( t\right) =\left( \partial_{y_{j}}\sigma _{i}\right)\left(t,X_{t,i}^{\mathbf{u}},\mathbf{X}_{t}^{\mathbf{u}},u_{t,i}\right), \\
\bar{f}_{t,i} =&\left( \partial _{u}\sigma _{i}\right) \left( t,X_{t,i}^{\mathbf{u}},\mathbf{X}_{t}^{\mathbf{u}},u_{t,i}\right) \delta_{h,i}u_{t,h}^{\prime} 
\end{aligned}
\end{eqnarray*}
yields that for any $i\in I_{N}$,
\begin{eqnarray*}
\mathbb{E}\left[ \left\vert Y_{t,i}^{\mathbf{u},u_{h}^{\prime
}}\right\vert ^{p}\right]
\leq \left[ T\left( \frac{L_{y}^{b}}{N}+3\left( p-1\right) \frac{\left(
L_{y}^{\sigma }\right) ^{2}}{N}\right) e^{\bar{I}_{b,\sigma ,p}^{3}\cdot
T}\left( L^{b}+3\left( p-1\right) L^{\sigma }\right) +\left( 3p-2\right)
\delta _{h,i}\right] \left\Vert u_{h}^{\prime }\right\Vert _{\mathcal{H}%
^{p}\left( \mathbb{R}\right) }^{p}e^{\bar{I}_{b,\sigma ,p}^{4}\cdot T},
\end{eqnarray*}
with
\begin{eqnarray*}
\begin{aligned}
\bar{I}_{b,\sigma ,p}^{3} = & \, pL^{b}+L_{y}^{b}p+\left( p-1\right) +\frac{3}{2}%
\left( p-1\right) p\left( L^{\sigma }\right) ^{2} \\
&+\frac{3}{2}\left( p-1\right) \left( p-2\right) \left( L_{y}^{\sigma
}\right) ^{2}+\frac{3}{2}\left( p-1\right) \left( p-2\right) +3\left( p-1\right) \left( L_{y}^{\sigma }\right) ^{2}, \\
\bar{I}_{b,\sigma ,p}^{4} = & \, pL^{b}+L_{y}^{b}\left( p-1\right) +\left(p-1\right) +\frac{3}{2}\left( p-1\right) p\left( L^{\sigma }\right) ^{2} \\
&+\frac{3}{2}\left( p-1\right) \left( p-2\right) \left( L_{y}^{\sigma
}\right) ^{2}+\frac{3}{2}\left( p-1\right) \left( p-2\right), 
\end{aligned}
\end{eqnarray*}
where we adopt the fact that $\left\Vert B\right\Vert _{\infty }\leq L^{b},$ $\left\Vert D\right\Vert
_{\infty }\leq L^{\sigma },$ $\left\Vert \bar{B}\right\Vert _{\infty }\leq
\frac{L_{y}^{b}}{N}$, and $\left\Vert \bar{D}\right\Vert _{\infty }\leq \frac{L_{y}^{\sigma }}{N}$.  
\hfill $\Box$

\paragraph{Proof of Corollary \ref{comthe}.}
We formulate SDE (\ref{comnoise}) as
\begin{equation}
\left\{
\begin{array}{rcl}
\mathrm{d}X_{t,i} & = & \mathsf{b}_{t,i}u_{t,i}\mathrm{d}t+\left( \sigma
_{t,i},1\right) \mathrm{d}\tilde{W}_{t,i}, \\
X_{0,i} & = & \xi _{i},%
\end{array}%
\right.   \label{comnoise2}
\end{equation}%
where $\tilde{W}_{t,i}=\left(
\begin{array}{c}
W_{t,i} \\
W_{t}^{0}%
\end{array}%
\right) $, for $i\in I_{N},$ which is a two dimensional Brownian motion.
Introduce the following notations
\begin{equation*}
\begin{aligned} \mathbf{B}_{0,x,y}\left(
t,\mathbf{X}_{t}^{\mathbf{u}},\mathbf{u}_{t}\right) =& \,
\mathbf{\bar{B}}_{x}\left(
t,\mathbf{X}_{t}^{\mathbf{u}},\mathbf{u}_{t}\right) +\mathbf{B}_{0,y}\left(
t,\mathbf{X}_{t}^{\mathbf{u}},\mathbf{u}_{t}\right), \\ \mathbf{\Pi
}_{0,x,y}^{j}\left( t,\mathbf{X}_{t}^{\mathbf{u}},\mathbf{u}_{t}\right) =&
\, \mathbf{\bar{\Pi}}_{x}^{j}\left(
t,\mathbf{X}_{t}^{\mathbf{u}},\mathbf{u}_{t}\right) +\mathbf{\Pi
}_{0,y}^{j}\left( t,\mathbf{X}_{t}^{\mathbf{u}},\mathbf{u}_{t}\right) ,\quad
\quad j\in I_{N}, 
\end{aligned}
\end{equation*}
with

\begin{equation*}
\begin{aligned}
\mathbf{\bar{B}}_{x}\left(
t,\mathbf{X}_{t}^{\mathbf{u}},\mathbf{u}_{t}\right) =& \, 0_{N\times
N},\text{ }\mathbf{B}_{0,y}\left(
t,\mathbf{X}_{t}^{\mathbf{u}},\mathbf{u}_{t}\right) =0_{N\times N}, \\
\mathbf{B}_{1,u}^{i}\left(
t,\mathbf{X}_{t}^{\mathbf{u}},\mathbf{u}_{t}\right) =& \, \left[ 0,\ldots
,\mathsf{b}_{t,i},\ldots ,0\right] ^{\top }{}_{N\times 1},\text{ }i\in I_{N},
\end{aligned}
\end{equation*}%
and
\begin{equation*}
\mathbf{\bar{\Pi}}_{x}^{i}\left( t,\mathbf{X}_{t}^{\mathbf{u}},\mathbf{u}%
_{t}\right) =0_{N\times N},\text{ }\mathbf{\Pi }_{0,y}^{i}\left( t,\mathbf{X}%
_{t}^{\mathbf{u}},\mathbf{u}_{t}\right) =0_{N\times N},\text{ }\mathbf{\Pi }%
_{1,u}^{j}\left( t,\mathbf{X}_{t}^{\mathbf{u}},\mathbf{u}_{t}\right)
=0_{N\times 1},\text{ }j\in I_{N}.
\end{equation*}%
Then the sensitive process is
\begin{equation}
\left\{
\begin{array}{rcl}
\mathrm{d}\mathbf{Y}_{t}^{\mathbf{u},u_{h}^{\prime }} & = & \mathbf{B}%
_{1,u}^{h}\left( t,\mathbf{X}_{t}^{\mathbf{u}},\mathbf{u}_{t}\right)
u_{t,h}^{\prime }\mathrm{d}t, \\
\mathbf{Y}_{0}^{\mathbf{u},u_{h}^{\prime }} & = & 0.%
\end{array}%
\right.
\end{equation}%
Now consider the following BSDE
\begin{equation}
\left\{
\begin{array}{rcl}
-\mathrm{d}P_{t,i} & = & \left( \partial _{x}f_{i}\right) \left( t,\mathbf{X}%
_{t}^{\mathbf{u}},\mathbf{u}_{t}\right) \mathrm{d}t-\sum_{j=1}^{N}\left(
Q_{t,i,j},Q_{t,i}^{0}\right) \mathrm{d}\tilde{W}_{t}^{j}, \\
P_{T,i} & = & \partial _{x}g_{i}\left( \mathbf{X}_{T}^{\mathbf{u}}\right) ,%
\text{ }\forall i\in I_{N},%
\end{array}%
\right.   \label{adjcom}
\end{equation}%
where
\begin{eqnarray*}
\partial _{x}g_{i}\left( \mathbf{X}_{T}^{\mathbf{u}}\right)  &=&G_{i}\left(
X_{T,i}-X_{T}^{\left( N\right) }\right) \left[ -\frac{1}{N},\ldots \frac{%
\left( N-1\right) }{N},\ldots ,-\frac{1}{N}\right] ^{\top }, \\
\left( \partial _{x}f_{i}\right) \left( t,\mathbf{X}_{t}^{\mathbf{u}},%
\mathbf{u}_{t}\right)  &=&\hat{Q}_{t,i}\left( X_{i}-X_{T}^{\left( N\right)
}\right) \left[ -\frac{1}{N},\ldots \frac{\left( N-1\right) }{N},\ldots ,-%
\frac{1}{N}\right] ^{\top }.
\end{eqnarray*}%
By Lemma \ref{bs2}, there exists a unique adapted solution $\left(
P_{t,i},\left\{ Q_{t,i,j}\right\} _{j=1,\ldots ,N},Q_{t,i}^{0}\right)
_{0\leq t\leq T}$ to BSDE (\ref{adjcom}). From Theorem \ref{the1}, we get
the desired result. \hfill $\Box $

\bibliographystyle{plain}
\bibliography{reference}

\end{document}